\documentclass[10pt,reqno]{amsart}
\usepackage{amsmath,mathrsfs}
\usepackage{amsfonts}
\usepackage{amssymb}
\usepackage{graphicx}
\usepackage{amsthm,color,yfonts}
\usepackage{pdfsync}
\usepackage{epstopdf}%
\usepackage{pifont}
\usepackage{bbm}
\usepackage[margin=4.1cm]{geometry}
\usepackage[colorlinks=true]{hyperref}
\hypersetup{linkcolor=red,citecolor=blue,filecolor=dullmagenta,urlcolor=blue} 

\usepackage{wrapfig}
\usepackage{tikz}
\usetikzlibrary{arrows,calc,decorations.pathreplacing}
\definecolor{light-gray1}{gray}{0.90}
\definecolor{light-gray2}{gray}{0.80}
\definecolor{light-gray3}{gray}{0.60}

\newcommand{\wt}{\widetilde}

\newcommand{\C}{{\mathbb C}}
\newcommand{\R} {\mathbb R}
\newcommand{\cuad}{{\sqcap\kern-.68em\sqcup}}

\definecolor{darkgreen}{rgb}{0.2,0.7,0.1}

\def\wh#1{\widehat{#1}}

\newcommand{\Z}{\mathbb{Z}}
\newcommand{\T}{\mathbb{T}}

\newcommand{\med}{\operatorname{med}}
\newcommand{\vu}{\mathbf{u}}

\newcommand{\vH}{\mathbb{H}}

\def\bm{\left( \begin{array}{cc}}
\def\endm{\end{array}\right)}

\providecommand{\norm}[1]{\left\| #1 \right\|}

\newcommand{\be}{\begin{equation}}
\newcommand{\ee}{\end{equation}}
\newcommand{\ba}{\begin{aligned}}
\newcommand{\ea}{\end{aligned}}
\newcommand{\bl}{\begin{lemma}}
\newcommand{\el}{\end{lemma}}
\newcommand{\bpm}{\begin{pmatrix}}
\newcommand{\epm}{\end{pmatrix}}

\numberwithin{equation}{section}

\newtheorem{theorem}{Theorem}[section]
\newtheorem{theoremA}{Theorem}

\newtheorem{lemmaB}{Lemma}

\newtheorem*{theorem*}{Theorem}
\newtheorem{proposition}{Proposition}[section]
\newtheorem{corollary}{Corollary}[section]
\newtheorem{lemma}{Lemma}[section]

\theoremstyle{remark}
\newtheorem{remark}{Remark}[section]

\title[2D \emph{abcd} decay]{Local energy decay for small radial solutions of the two-dimensional Hamiltonian $abcd$ Boussinesq system}

\author[C. Kwak]{Chulkwang Kwak}
\address{Department of Mathematics, Ewha Womans University, Seoul 03760 and Korea Institute for Advanced Study, Seoul 02455, Republic of Korea}
\email{ckkwak@ewha.ac.kr}

\author[C. Mu\~noz]{Claudio Mu\~noz}
\address{CNRS and Departamento de Ingenier\'{\i}a Matem\'atica and Centro
de Modelamiento Matem\'atico (UMI 2807 CNRS), Universidad de Chile, Casilla
170 Correo 3, Santiago, Chile.}
\email{cmunoz@dim.uchile.cl}
\thanks{C. K. was partially supported by Young Research Program of the National Research Foundation of Korea(NRF) grant funded by the Korea government(MSIT) (No. RS-2023-00210210) and Global - Learning \& Academic research institution for Master’s·PhD students, and Postdocs(G-LAMP) Program of the National Research Foundation of Korea(NRF) grant funded by the Ministry of Education(No. RS-2025-25442252). C.~M.'s work was partly funded by Inria Lille PANDA and Chilean research grants FONDECYT 1231250 and Centro de Modelamiento Matem\'atico (CMM) BASAL fund FB210005 for center of excellence from ANID-Chile.}


\begin{document}

\begin{abstract}
We consider the $abcd$ Boussinesq system in two space dimensions, a four-parameter family of models for the two-way propagation of small-amplitude long surface waves, in the generic regime $a,c<0$, $b,d>0$. We first revisit the local well-posedness theory in Sobolev spaces in dimensions $2\le n\le4$ and show that for the Hamiltonian case $b=d$ the Cauchy problem is globally well-posed for small data in the energy space $\mathbb{H}^1(\mathbb R^n)$. Then we consider the asymptotic behavior of these global solutions in the $n=2$ case. Under the condition $-a,-c>b/9$, we prove that every small solution with radial symmetry exhibits local energy decay. As a consequence, there are no small radial time-periodic (in particular, stationary) solutions. Decay also holds, along sequences of times, on balls whose radius grows almost linearly in time. The proof is based on new virial identities adapted to the radial structure of the system. In the proof we introduce two correction functionals such that the resulting quadratic form is coercive precisely when $-a,-c>b/9$, a threshold which is sharp for the functional. To our knowledge this is the first result describing the long-time dynamics of $abcd$ systems in dimension two. 
\end{abstract}

\maketitle
\tableofcontents

\section{Introduction}

\subsection{General setting}
The mathematical study of surface water waves starts from the free-surface \emph{Euler equations}, which describe the motion of an inviscid and incompressible fluid in two or three spatial dimensions. Assuming further that the flow is irrotational and that the bottom of the channel is flat and horizontal, one may focus on specific regimes of wave amplitude and wavelength and derive asymptotic models that capture the dynamics of the solutions in those regimes. Of particular interest is the small-amplitude, long-wave regime, known as the \emph{Boussinesq regime}, characterized by
\[\varepsilon = \frac{A}{h} \sim \left(\frac{h}{\lambda}\right)^2 \ll 1,\]
where $h$ denotes the mean depth of the fluid, $A$ a typical wave amplitude and $\lambda$ a typical horizontal wavelength. In this regime, nonlinear and dispersive effects are of the same order, and their balance gives rise to a rich family of mathematically interesting models. Among them are the \emph{abcd Boussinesq systems}, introduced by Bona, Chen and Saut \cite{BCS1} as a generalization of the classical Boussinesq model \cite{Bous}. They are obtained through a formal asymptotic expansion of the Euler equations in suitably scaled variables, which yields 
\begin{equation}\label{2D ABCD_e}
\begin{aligned}
&\partial_t \eta + \nabla \cdot u + \varepsilon \left( a\, \Delta \nabla \cdot u + \nabla \cdot (\eta u) - b\, \Delta \partial_t \eta \right) = O(\varepsilon^2), \\
&\partial_t u + \nabla \eta + \varepsilon \left( c\, \Delta \nabla \eta + \tfrac12 \nabla |u|^2 - d\, \Delta \partial_t u \right) = O(\varepsilon^2),
\end{aligned}
\end{equation}
where the unknowns $\eta$ and $u$ denote, respectively, the deviation of the free surface from its rest position and the (scaled) horizontal velocity at height $\theta h$ above the bottom, $\theta\in[0,1]$ being the parameter defined in \eqref{Conds}. Neglecting the $O(\varepsilon^2)$ terms leads to the \emph{abcd} systems. 

The parameters $a,b,c$ and $d$ are not arbitrary, and must obey the following conditions (at least in the regime without surface tension):\footnote{In the regime with surface tension, $a$ will be replaced by $a - \tau$ in \eqref{Conds}, where $\tau \ge 0$ is the Bond number which measures surface tension effects, see \cite{DD}.}
\begin{equation}\label{Conds}
a +b =\frac1{2}\left(\theta^2-\frac13\right),\quad c+d =\frac1{2} (1-\theta^2)\geq 0, \quad \mbox{where} \quad \theta \in [0,1].
\end{equation}
The three degrees of freedom in the relation (immediately follows from \eqref{Conds}) $a+b+c+d = \frac13$ arise from the height at which the horizontal velocity is taken and from a double use of the Benjamin-Bona-Mahony (BBM) trick.  On the other hand, the parameters satisfying \eqref{Conds} have been classified in the sense of  well-posedness. In \cite{BCS1, BCS2}, the following particular cases are referred to as \emph{linearly well-posed regimes}:\footnote{The case when $b=d<0$ and $a=c>0$ has been dealt with in \cite{BCS1} as ``linearly" well-posed system, while it was excluded in \cite{BCS2} due to the singularity of operators $I +a\partial_x^2$ and $I-b\partial_x^2$.}
\begin{align}
& a<0, ~c< 0, ~ b >0,~d > 0 \label{generic}\\
& a=c > 0, ~ b \ge0,~d \ge 0 \label{strong}.
\end{align}  
The first one \eqref{generic} is usually referred to as \emph{generic} case. Case \eqref{strong} will not be considered in this work.

A natural question is whether the corresponding Boussinesq systems provide a rigorous approximation to the full water waves system. This leads to the problem of the rigorous \emph{justification} of the asymptotic models, namely, showing that solutions to the asymptotic models converge strongly to solutions to the full water waves system. In the generic regimes without surface tension, this question has been completely settled by combining several independent results \cite{BCL, ASL, MSZ2012, SX, SWX}. We briefly recall here the notion of justification and the results available so far. The rigorous justification of the Boussinesq systems can be carried out in the following four steps. First, \textit{consistency}: one shows that the full water waves system is consistent with the Boussinesq systems; for any family of solutions $(\phi^{\varepsilon}, \eta^{\varepsilon})$ to the Euler system, the corresponding pair $(\vu^{\varepsilon}, \eta^{\varepsilon})$\footnote{The velocity $\vu^{\varepsilon}$ is determined by $\phi^{\varepsilon}$; namely, $\vu^{\varepsilon} := \nabla \psi^{\varepsilon}$ with $\psi^{\varepsilon}(t,X) = \phi^{\varepsilon}\big(t, X, 1 + \varepsilon \eta^{\varepsilon}(t,X)\big)$.} satisfies \eqref{2D ABCD_e} up to a residual of order $O(\varepsilon^2)$, bounded in a suitable functional class. The second step is the \textit{well-posedness of the full water waves system}: one establishes the well-posedness of the water waves system on the relevant time scale, of order $O(1/\varepsilon)$. The third step concerns the \textit{long-time existence for the Boussinesq systems}. Here one proves that solutions to the Boussinesq systems exist at least on time scales of order $O(1/\varepsilon)$, together with the appropriate uniform bounds. Finally, the fourth step is to prove \textit{error estimates}. One proves that solutions to the Boussinesq systems converge strongly to those of the full water waves system, with an error of order $O(\varepsilon^2 t)$; that is, $\| \mathcal{U}_{\mathrm{Euler}} - \mathcal{U}_{\mathrm{Bous}} \| = O(\varepsilon^2 t)$, with $t \in [0, T/\varepsilon]$, and where $\mathcal{U}_{\mathrm{Euler}}$ and $\mathcal{U}_{\mathrm{Bous}}$ denote solutions to the full water waves system and to the Boussinesq systems, respectively. The first and last steps were achieved by Bona, Colin and Lannes \cite{BCL} in 2005, see also \cite{BLS}. The second step was established by Alvarez-Samaniego and Lannes \cite{ASL} in 2008 in the absence of surface tension, and by Ming, Zhang and Zhang \cite{MZZ1, MZZ2} in 2012 in the presence of surface tension. The third step was addressed in \cite{MSZ2012, SX, SWX}. In \cite{MSZ2012}, the authors obtained the first long-time existence result for the Boussinesq systems in two particular regimes: the generic case ($a, c < 0$ and $b, d > 0$) and the BBM--BBM case ($a = c = 0$ and $b, d > 0$). All the remaining admissible, linearly well-posed systems were treated in \cite{SX}, with the exception of the \emph{strongly dispersive} systems, corresponding to $b = d = 0$ with either $a = 0,\ c < 0$ or $a < 0,\ c = 0$. An alternative proof of long-time existence on time scales of order $O(1/\varepsilon)$ was given in \cite{SWX} for most of the Boussinesq systems ($a, c \le 0$ and $b, d \ge 0$), with the exception of the two cases $a = d = 0$, $c < 0$, $b > 0$ and $a = b = 0$, $c < 0$, $d > 0$, which in particular involve the strongly dispersive case. We also remark that all proofs in \cite{MSZ2012, SX, SWX} relax the noncavitation condition $1+\eta\geq 0$. The KdV-KdV case was treated in depth in \cite{LPS}, but it is not considered in this work. For the symmetric case, see \cite{BCL}. For ill-posedness results, see \cite{CL2013,KwakMaulen}.

For a comprehensive account of the derivation and the rigorous justification of shallow-water asymptotic models, including the $abcd$ family, we refer to the monograph \cite{Lannes2013}. The particular case $a=c=0$, $b=d>0$ is the BBM--BBM system, whose theoretical and numerical study goes back to \cite{BonaChen1998}. Long-time existence on the Boussinesq time scale $O(1/\varepsilon)$ for most $abcd$ systems, including the case of non-flat bottoms, was also established in \cite{Burtea2016}. Two-dimensional $abcd$ systems on bounded domains, as well as their numerical approximation, were studied in \cite{DMS2009}. 

\medskip

We now focus on the two-dimensional $abcd$ system in the generic regime. Neglecting the $O(\varepsilon^2)$ terms in \eqref{2D ABCD_e} and rescaling the variables, we obtain
\begin{equation}\label{2D ABCD}
\begin{aligned}
&(1-b\Delta)\partial_t\eta+\nabla\cdot\left(u+a\Delta u+\eta u\right)=0,\\
&(1-d\Delta)\partial_tu+\nabla\left(\eta+c\Delta\eta+\frac12|u|^2\right)=0,
\end{aligned}
\end{equation}
for $(t,x)\in\R\times\R^2$, where $\eta=\eta(t,x)$ is real-valued and $u=u(t,x)\in\R^2$.

\subsection{Main results}

The main purpose of this work is to establish local energy decay for small radial solutions to the two-dimensional $abcd$ systems \eqref{2D ABCD}. As a first step, we develop the local and global well-posedness theory in the energy space. Although the long-time dynamics studied in this paper are two-dimensional, the well-posedness argument extends naturally to higher dimensions and yields small data global well-posedness in $H^1(\R^n)$ for $2\le n\le4$. The one-dimensional well-posedness theory was already developed by Bona, Chen and Saut \cite{BCS1, BCS2}.

\subsubsection{Well-posedness theory}\label{sec:WP theory}

The classical notion of well-posedness was introduced by Hadamard \cite{Hadamard}, that is to say, the mathematical models of physical phenomena should have the following properties:\footnote{A precise formulation can be found in \cite{MS}.} $(i)$ there exists a solution; $(ii)$ the solution is unique; and $(iii)$ the solution behaves continuously with the initial condition.

For the well-posedness theory, we consider the $abcd$ system in $\R^n$ for $n\ge2$:
\begin{equation}\label{nD ABCD}
\begin{aligned}
&(1-b\Delta)\partial_t\eta+\nabla\cdot\left(a\Delta u+u+\eta u\right)=0,\\
&(1-d\Delta)\partial_tu+\nabla\left(c\Delta\eta+\eta+\frac12|u|^2\right)=0,\\
&\eta(0)=\eta_0,\qquad u(0)=u_0,
\end{aligned}
\end{equation}
where $\eta=\eta(t,x) \in \R$, $u=u(t,x) \in \R^n$, and $(t,x)\in\R\times\R^n$. In higher dimensions, we impose the irrotational condition
\begin{equation}\label{eq:irrotational condition}
\partial_{x_j}u_k=\partial_{x_k}u_j,\qquad 1\le j,k\le n,
\end{equation}
which relates the curl-free condition arising in the derivation from the irrotational water waves problem. Therefore, on $\R^n$, one may write $u=\nabla\psi$ for a scalar potential $\psi$.

\begin{remark}\label{rem:propagation}
Notice that \eqref{eq:irrotational condition} is formally preserved by the flow. Indeed, applying the operator $\partial_{x_j}$ to the $k$-th component of the second equation in \eqref{nD ABCD}, and subtracting the symmetric expression obtained by interchanging $j$ and $k$, we obtain
\[(1-d\Delta)\,\partial_t\left(\partial_{x_j}u_k-\partial_{x_k}u_j\right)=0,\qquad 1\le j,k\le n,\]
so that $\partial_{x_j}u_k-\partial_{x_k}u_j$ is formally constant in time. In particular it vanishes identically if it vanishes at $t=0$. Thus it suffices to impose \eqref{eq:irrotational condition} on the initial data. Similarly, \eqref{nD ABCD} is invariant under rotations $x\mapsto Rx$, with $R\in {\rm O}(n)$, acting as $(\eta,u)(t,x)\mapsto(\eta(t,Rx),R^{T}u(t,Rx))$. Hence, by uniqueness, if $\eta_0$ is radial and $u_0(x)=q_0(|x|)\,x/|x|$, the solution keeps this form for all times: $\eta(t,\cdot)$ is radial and $u(t,x)=q(t,|x|)\,x/|x|$. Note that such a vector field is automatically irrotational, being the gradient of the radial function $\int_0^{|x|}q(t,s)\,ds$.
\end{remark}

The system is called \emph{Hamiltonian} if $b=d$. In this case, sufficiently regular solutions to \eqref{nD ABCD} conserve the energy
\begin{equation}\label{ENERGY}
E(u,\eta)=\frac12\int_{\R^n}\left(-a|Du|^2-c|\nabla\eta|^2+|u|^2+\eta^2+\eta|u|^2\right)\;dx=E(u_0,\eta_0).
\end{equation}
Here $Du$ denotes the Jacobian matrix of $u$, so that
\[|Du|^2=\sum_{j,k=1}^n|\partial_{x_k}u_j|^2.\]
Indeed, for sufficiently regular solutions, a direct differentiation of the energy gives
\[\frac{d}{dt}E(u,\eta)=\int_{\R^n}\left(\eta+c\Delta\eta+\frac12|u|^2\right)\eta_t+\left(u+a\Delta u+\eta u\right)\cdot u_t\;dx.\]
Using \eqref{nD ABCD}, integration by parts, the self-adjointness of $(1-b\Delta)^{-1}$ and $(1-d\Delta)^{-1}$, and the resolvent identity
\[(1-b\Delta)^{-1}-(1-d\Delta)^{-1}=(b-d)(1-b\Delta)^{-1}(1-d\Delta)^{-1}\Delta,\]
one obtains
\[\begin{aligned}
\frac{d}{dt}E(u,\eta)
=&~{}\int_{\R^n}\left(u+a\Delta u+\eta u\right)\cdot(1-b\Delta)^{-1}\nabla\left(\eta+c\Delta\eta+\frac12|u|^2\right)\;dx\\
&~{}-\int_{\R^n}\left(u+a\Delta u+\eta u\right)\cdot(1-d\Delta)^{-1}\nabla\left(\eta+c\Delta\eta+\frac12|u|^2\right)\;dx\\
=&~{}(b-d)\int_{\R^n}\left(u+a\Delta u+\eta u\right)\cdot(1-b\Delta)^{-1}(1-d\Delta)^{-1}\Delta\nabla\left(\eta+c\Delta\eta+\frac12|u|^2\right)\;dx.
\end{aligned}\]
Hence $\frac{d}{dt}E(u,\eta)=0$ when $b=d$, which gives the conservation of \eqref{ENERGY}.

The system also satisfies the conservation laws
\[M_j(u_j,\eta)=\int_{\R^n}u_j\eta+b\nabla u_j\cdot\nabla\eta\;dx=M_j(u_{j,0},\eta_0),\qquad j=1,\ldots,n,\]
and
\[\int_{\R^n}u_j\;dx=\int_{\R^n}u_{j,0}\;dx,\qquad j=1,\ldots,n,\qquad \int_{\R^n}\eta\;dx=\int_{\R^n}\eta_0\;dx.\]
The conservation of $M_j$ follows similarly from \eqref{nD ABCD} and the irrotational condition \eqref{eq:irrotational condition}, while the conservation of the spatial means follows directly by integrating \eqref{nD ABCD} over $\R^n$. Although $M_j$ will not be used in this paper, it will define virial functionals by suitable modifications.

As concerns the local theory, we first recall the following two-dimensional result by Dougalis, Mitsotakis and Saut \cite{DMS2007}.
\begin{theoremA}[\cite{DMS2007}]\label{thm:LWP0}
Let $s>0$ and assume that $a,c<0$ and $b,d>0$. Then, for every initial data $(u_0, \eta_0)\in H^s(\R^2)^2 \times H^s(\R^2)$, there exist $T>0$ and a unique solution $(u, \eta)\in C([0,T];H^s(\R^2)^2 \times H^s(\R^2))$ to the two-dimensional $abcd$ system \eqref{2D ABCD}.
\end{theoremA}

We extend this local theory to higher dimensions and, in particular, include the endpoint regularity in dimensions $n\ge3$.

\begin{theorem}\label{thm:LWP}
Let $n\ge2$, and assume that $s>0$ if $n=2$, and $s\ge\frac{n-2}{2}$ if $n\ge3$. Then the \emph{generic} $abcd$ system \eqref{nD ABCD}, under the irrotational condition \eqref{eq:irrotational condition}, is locally well-posed in $\vH^s(\R^n)$. The precise definition of $\vH^s$ will be given in Section \ref{sec:notations}.
\end{theorem}

The proof of Theorem \ref{thm:LWP0} is based on the following classical bilinear estimate of Grisvard \cite{Grisvard}.

\begin{lemmaB}[Grisvard \cite{Grisvard}]\label{lem:bi0}
Let $s_1,s_2,s_3\in\R$ be such that $s_1\ge s_3$, $s_2\ge s_3$, $s_1+s_2\ge0$, and $s_1+s_2-s_3>\frac n2$. Then, the bilinear map $(f,g)\mapsto fg$ is continuous from $H^{s_1}(\R^n)\times H^{s_2}(\R^n)$ into $H^{s_3}(\R^n)$.
\end{lemmaB}

Taking $s_1=s_2=s$ and $s_3=s-1$, Lemma \ref{lem:bi0} yields
\[\norm{\partial_{x_j}(1-\Delta)^{-1}(fg)}_{H^s(\R^n)}\lesssim\norm{fg}_{H^{s-1}(\R^n)}\lesssim\norm{f}_{H^s(\R^n)}\norm{g}_{H^s(\R^n)},\]
which, in dimension two, is the form used in the local well-posedness argument of \cite{DMS2007}. This estimate holds provided that $s>\frac{n-2}{2}$ for $n\ge2$. We improve the endpoint regularity in dimensions $n\ge3$ as follows.

\begin{proposition}\label{prop:refined bilinear}
Let $n \ge 2$, and assume that $s>0$ if $n=2$, and $s\ge\frac{n-2}{2}$ if $n\ge3$. Then, for all $f,g\in H^s(\R^n)$,
\[\norm{(1-\Delta)^{-1}|\nabla|(fg)}_{H^s(\R^n)}\lesssim\norm{f}_{H^s(\R^n)}\norm{g}_{H^s(\R^n)}.\]
\end{proposition}

Theorem \ref{thm:LWP} follows from Proposition \ref{prop:refined bilinear} together with the contraction mapping principle.

\begin{remark}[Unconditional uniqueness]
Proposition \ref{prop:refined bilinear} also yields unconditional uniqueness for distributional solutions in $L^\infty((0,T);\vH^s(\R^n))$ by applying the argument of Kim and the first author \cite{KimKwak2026}. In particular, any two distributional solutions in this class with the same initial data coincide on $[0,T]$.
\end{remark}

\begin{remark}\label{rem:SLWP}
Theorem \ref{thm:LWP} privides local well-posedness in the energy space $\vH^1(\R^n)$ for $2 \le n \le 4$. Notice that the case $n=4$ is precisely the one provided by Proposition \ref{prop:refined bilinear}, since Lemma \ref{lem:bi0} with $s_1=s_2=1$ and $s_3=0$ only gives the energy estimate when $n \le 3$. Thus, Proposition \ref{prop:refined bilinear} extends the local theory in the energy space up to dimension $4$. Moreover, for every parameter $(a,b,c,d)$ for which the system \eqref{nD ABCD} is \emph{generic}, there exists $\epsilon_0=\epsilon_0(a,b,c,d,n)>0$ such that if
\[\norm{(u_0,\eta_0)}_{\vH^1}\le\epsilon_0,\]
then there exists a unique solution $(u,\eta)\in C([0,1];\vH^1)$ in the sense of Theorem \ref{thm:LWP}, provided $2\le n\le4$. This uniform local existence time for small initial data will be used below, together with the small data a priori energy estimate, to establish small data global well-posedness in the energy space.
\end{remark}


\begin{theorem}\label{thm:GWP}
Let $2\le n\le4$ and $s\ge1$. Assume that the parameters $(a,b,c,d)$ satisfy the generic condition \eqref{generic} and the Hamiltonian condition $b=d$. Then the $abcd$ system \eqref{nD ABCD}, under the irrotational condition \eqref{eq:irrotational condition}, is globally well-posed in $\vH^s(\R^n)$ for initial data in $\vH^s(\R^n)$ whose $\vH^1(\R^n)$ norm is sufficiently small. Moreover, the solution satisfies the uniform bound
\[\sup_{t\in\R}\|(u(t),\eta(t))\|_{\vH^1(\R^n)}\lesssim \|(u_0,\eta_0)\|_{\vH^1(\R^n)}.
\]
\end{theorem}

The proof of Theorem \ref{thm:GWP} rests on two ingredients. The first is a diagonalization of the linear part of \eqref{nD ABCD} in Fourier variables. The corresponding linear propagator is unitary on $H^s$, while the nonlinearity gains one full derivative. A fixed-point argument based on a bilinear estimate gives local well-posedness as stated in Theorem \ref{thm:LWP}. The second ingredient is the coercivity of the Hamiltonian \eqref{ENERGY} on small sets of the energy space. The cubic term
\[\int_{\R^n}\eta|u|^2\,dx\]
is controlled through the Sobolev embedding $H^1(\R^n)\hookrightarrow L^3(\R^n)$, which yields a small multiple of the $\vH^1$-norm together with a quartic error. Consequently, for sufficiently small initial data, the conserved energy controls the $\vH^1$-norm uniformly in time. This a priori argument is valid for $n\le6$, while the local well-posedness theory obtained from Theorem \ref{thm:LWP} is available for $2 \le n \le 4$. 

\begin{remark}
To the best of our knowledge, Theorem \ref{thm:GWP} gives the first global well-posedness result for the $abcd$ system valid in higher dimensions. A related two-dimensional global existence result was obtained by Saut and Xu for the Hamiltonian ``Boussinesq-Full dispersion" systems \cite{SautXu}. Their argument is also based on the coercivity of the Hamiltonian. In the present setting, we use the exact Sobolev embedding $H^1(\R^n)\hookrightarrow L^3(\R^n)$ directly, which gives an a priori argument valid up to dimension $6$. 
\end{remark}

\begin{remark}
The local and global well-posedness results above stated also hold on the periodic domain $\T^n$. Indeed, the bilinear estimate does not rely on any dispersive property of the linear flow and remains valid in the periodic setting. Therefore the iteration argument for Theorem \ref{thm:LWP} holds in $\T^n$, and the same energy argument and persistence of regularity provides the corresponding global well-posedness result of Theorem \ref{thm:GWP}.
\end{remark}

\subsubsection{Local energy decay}
We now consider the long-time behavior of small solutions to the $abcd$ systems. Since $abcd$ systems possess solitary wave solutions \cite{Chen1998,ChenNguyenSun2010,ChenNguyenSun2011}, no decay is expected for large solutions in the full $abcd$ regime of parameters. Therefore, the smallness assumption in Theorem~\ref{thm:MT} is natural and sharp. In one space dimension, the long-time behavior of small solutions to the Hamiltonian $abcd$ system was considered in \cite{KMPP2019} and in \cite{KM2020}. The present paper is the two-dimensional counterpart of these works, proposing a road path towards the study of local energy decay in higher dimensions. The new difficulty is that the virial functionals used in previous works \cite{KMPP2019,KM2020} are no longer coercive in dimension two or higher. 

The introduction and development of virial identities to describe the energy based long-time dynamics of dispersive models of generalized KdV type goes back to the works by Martel and Merle \cite{MartelMerle2000,MM0,MM1,MMnon,MMas1,MMas2}. Recently, iterated virial estimates in the energy space for scalar field models were developed in \cite{KMM2017,KMM2017b,KMMV2021,CM2023,CLY2026} in the context of kinks and breathers. Local energy decay results for arbitrarily large data have been obtained for the KdV equation \cite{MunozPonce2019}, for the intermediate long wave equation \cite{MPS2021}, for the BBM equation \cite{KwakMunoz2019}, for Boussinesq-type equations \cite{MPP2018,MaMu0,Mau,MaMu,MaMu2,BCSM2026} and other related models \cite{MMP,BCM2026,CM2027,CSM2026}. In two space dimensions, virial identities have been used to prove local energy decay for the Zakharov--Kuznetsov and Kadomtsev-Petviashvilii equations \cite{MMPP2021,MMPP2}, as well as asymptotic stability of its solitons \cite{CMPS2016}. In those works the anisotropic structure of the equation (a preferred direction of movement) plays an essential role. Here we rely on radial symmetry as the main symmetry assumption. 

We consider the two-dimensional system \eqref{2D ABCD} and restrict our attention to small radial solutions. Our main result is the following.

\begin{theorem}\label{thm:MT}
Assume $n=2$ and the conditions
\[b=d, \qquad -a>\frac{b}9,\qquad -c>\frac{b}9.\]
There exist $\varepsilon_0>0$ depending only on $a,b,c$, such that the following holds. Let $(\eta, u)$ be a global radial solution of \emph{abcd} obtained from Theorem~\ref{thm:GWP}, in the sense that $\eta(t,\cdot)$ is radial and $u(t,x)=q(t,|x|)\,x/|x|$ as a vector field, and satisfying
\[\sup_{t\in\R}\left(\|u(t)\|_{H^1(\R^2)}+\|\eta(t)\|_{H^1(\R^2)}\right)\le\varepsilon\le\varepsilon_0.\]
Then, for every $R >0$, with $B_R=\{x\in\R^2:|x|<R\}$,
\[
\lim_{t\to\pm\infty}\left(\|u(t)\|_{H^1(B_R)}+\|\eta(t)\|_{H^1(B_R)}\right)=0.
\]
\end{theorem}

Theorem~\ref{thm:MT} has two direct consequences. The first one is a rigidity statement: no small radial solution can be periodic in time.

\begin{corollary}\label{cor:nonexistence}
Under the parameter assumptions of Theorem~\ref{thm:MT}, let $(u,\eta)$ be a global radial solution of \eqref{2D ABCD} with the same smallness condition from Theorem~\ref{thm:MT}. Assume that this solution is periodic in time, i.e.\ $(\eta,u)(t+T,\cdot)=(\eta,u)(t,\cdot)$ for some $T>0$ and all $t\in\R$. Then $(\eta,u)\equiv0$. In particular, \eqref{2D ABCD} has no nontrivial small radial stationary solutions, and no small radial breathers.
\end{corollary}

\begin{proof}
Fix $R>0$ and $t\in\R$. By hypothesis on periodicity and \eqref{eq:local H1 decay},
\[
\begin{aligned}
& \|u(t)\|_{H^1(B_R)}+\|\eta(t)\|_{H^1(B_R)} \\
&\quad =\lim_{k\to\infty}\left(\|u(t+kT)\|_{H^1(B_R)}+\|\eta(t+kT)\|_{H^1(B_R)}\right)=0.
\end{aligned}
\]
Since $R$ and $t$ are arbitrary, we have $(\eta,u)\equiv0$.
\end{proof}

Note that a solitary wave $Q(x-ct)$ with $c\neq0$ is never radial, so radial symmetry already excludes the solitary waves of \cite{Chen1998,ChenNguyenSun2010,ChenNguyenSun2011}. Consequently, Corollary~\ref{cor:nonexistence} shows that, in the small radial class, there are no other coherent structures either.

The final result in this paper extends the decay on fixed balls proved in Theorem~\ref{thm:MT} to decay, along sequences of times, on balls whose radius grows almost linearly in time.

\begin{corollary}\label{cor:expanding}
Under the assumptions of Theorem~\ref{thm:MT}, there exists $L_0\ge1$ depending only on $a,b,c$ such that the following holds. Let $T_0>0$ and let $L:[T_0,\infty)\to[L_0,\infty)$ be a nondecreasing $C^1$ function with $L(t)\to\infty$ as $t\to\infty$ and
\begin{equation}\label{eq:expanding condition}
\lim_{T\to\infty} \dfrac{\log L(T)}{\displaystyle \int_{T_0}^T L^{-1}(t) dt }=0.
\end{equation}
Then
\[\liminf_{t\to+\infty}\left(\|u(t)\|_{H^1(B_{L(t)})}+\|\eta(t)\|_{H^1(B_{L(t)})}\right)=0,\]
that is, there exists an increasing sequence of times $t_k\to+\infty$ such that 
\[
\lim_{k\to +\infty} \left( \|u(t_k)\|_{H^1(B_{L(t_k)})}+\|\eta(t_k)\|_{H^1(B_{L(t_k)})} \right)=0.
\]
The same holds as $t\to-\infty$. Finally, condition \eqref{eq:expanding condition} is satisfied, for instance, by the width $L(t)=t^{\alpha}$ with $0<\alpha<1$, and by $L(t)=t/(\log t)^{\beta}$ with $\beta>0$ and $T_0$ large.
\end{corollary}

The proof of Corollary~\ref{cor:expanding} is given in Subsection~\ref{sec:expanding}. The proof uses the virial identity obtained in Proposition~\ref{prop:Virial identity} with a time-dependent scale $L=L(t)$, and the fact that the dependence of a localized Hamiltonian at scale $L$ is of order $L^{-1}$.  The case $L(t)=t$ is excluded by \eqref{eq:expanding condition} and it is natural since small solitary waves may exist in this regime. In dimension one, decay was obtained in \cite{KMPP2019} on a proper subset of the light cone around $x=0$, and in \cite{KM2020} on the whole light cone, i.e.\ on intervals $|x-vt|\lesssim|t|/\log^2|t|$ for every $|v|<1$, as well as in the exterior region $|x|\gg|t|$. 

\begin{remark}
It is worth mentioning that the admissible set is non-empty; with $b=d$, \eqref{Conds} gives $a+c=\frac13-2b$. The three conditions $a,c<0$, $-a>b/9$, $-c>b/9$ are compatible if and only if $-(a+c)=2b-\frac13>\frac{2b}9$, i.e.\ if and only if $b=d>\frac3{16}$. In the symmetric case $a=c$, which implies $\theta^2=\frac23$ and $a=c=\frac16-b$, the condition reads $b=d>\frac3{16}$. 
\end{remark}

\begin{remark}
In one space dimension, decay of small solutions in the energy space was first obtained in \cite{KMPP2019} under the dispersive condition $3b(a+c)+2b^2<8ac$, and then extended in \cite{KM2020} to the weakly dispersive regime, under refined conditions, see \cite[Definition~2.3]{KM2020}. In terms of $\widetilde a=a/b$, $\widetilde c=c/b$, the latter result describes a region of the $(\widetilde a,\widetilde c)$ plane which is symmetric under the permutation $\widetilde a\leftrightarrow\widetilde c$. It is strictly larger than the (2D) set given by \eqref{eq:parameter condition}. The boundary of the 1D region touches the diagonal $\widetilde a=\widetilde c$ exactly at the point $\widetilde a=\widetilde c=-\frac19$, see \cite{KM2020}. In other words, every pair $(a,c)$ satisfying the two-dimensional radial condition $-a,-c>b/9$ also satisfies the one-dimensional conditions of \cite{KM2020}, while the converse fails outside the diagonal under the present results. On the diagonal $a=c$ the two thresholds coincide and both read $-a>b/9$. In \cite{KM2020} the threshold $-a=b/9$ is shown to be sharp in dimension one for $a=c$, because for $-a<b/9$ the group velocity of the linear $abcd$ waves vanishes at nonzero wave numbers satisfying $bk^2=3$. Here the same threshold appears through a suitable polynomial obtained in the Fourier variable, which has a double zero at frequencies $|\xi|^2=3$ in the Hamiltonian setting $b=d=1$. Thus the condition $-a,-c>b/9$ of Theorem~\ref{thm:MT} is the natural sharp two-dimensional counterpart of the sharp one-dimensional symmetric threshold. The nonsymmetric part of the region of \cite{KM2020} in our setting will be considered elsewhere.
\end{remark}

\begin{remark}
Notice that \eqref{2D ABCD} follows from \eqref{2D ABCD_e} through the scaling given by $\eta'(t',x')=\varepsilon\eta(\sqrt\varepsilon t',\sqrt\varepsilon x')$, $u'=\varepsilon u(\sqrt\varepsilon\,\cdot)$, and that under this scaling we have the identities $\|\eta'\|_{L^2}^2=\varepsilon\|\eta\|_{L^2}^2$, $\|\nabla\eta'\|_{L^2}^2=\varepsilon^2\|\nabla\eta\|_{L^2}^2$. The smallness assumption of Theorems~\ref{thm:GWP}--\ref{thm:MT} for \eqref{2D ABCD} is therefore equivalent to $O(1)$ data for the physical system \eqref{2D ABCD_e} in the Boussinesq regime $\varepsilon\ll1$.  The precise description of Theorem \ref{thm:MT} is given in Theorem \ref{thm:local energy decay}. 
\end{remark}

The proof of Theorem~\ref{thm:MT} is the heart of this paper. The strategy carried out in \cite{KMPP2019} is to consider the localized momentum 
\[
I_\phi(t)=\int\phi\,(\eta u+ (Du)\nabla\eta)\,dx,
\]
with $\phi$ a bounded weight supported at physical scales $\lesssim L$, and to show that its time derivative is, up to small errors controlled using the smallness of the solution, a positive definite quadratic form in $(u,\eta)$ controlling a suitable localized $H^1$ energy. If this condition is satisfied, we have that time integrability of the local energy is obtained. Moreover, its complete time decay will follow by improved localized energy methods. In dimension one, this strategy succeeds \cite{KMPP2019} under a suitable choice in the $abcd$ parameters. In dimension two this strategy encounters a genuine obstruction. Indeed, under radial symmetry $r=|x|$, we have the decomposition 
\[
u(t,x)=q(t,r)e_r, \quad \eta(t,x)=h(t,r),
\]
with $e_r$ the radial unit vector. Then the quadratic form of the derivative of the virial functional involves not only the radial weight $\phi_r$, as in one dimension, but also the tangential weight $\phi/r$, coming from the angular part of the divergence and of the vector Laplacian. We have then that the coefficients of the gradient terms in the derivative of the virial functional $q_r^2$ and $h_r^2$ are proportional to to the weighted terms $\phi_r-\phi/r$ and $\phi_r+\phi/r$, respectively, and they should be at least nonnegative to get the minimally possible coercivity. Since $(\phi/r)_r=(\phi_r-\phi/r)/r$, the condition $\phi_r\ge\phi/r$ would force $\phi/r$ to be nondecreasing, which is against the boundedness of $\phi$. Hence no bounded radial weight makes the quadratic form coercive and the direct two-dimensional analogue of the one-dimensional virial argument fails.

The main new idea is to modify $I_\phi$ by two additional functionals that exactly cancel the tangential contributions. We introduce two additional functionals
\[
\mathcal J_\phi(t)=2\pi\int_0^\infty\frac{\phi}{r}\,h\,D_rq\,rdr,\qquad \mathcal K_\phi(t)=2\pi\int_0^\infty\frac{\phi}{r}\,q\,D_r^{-1}h\,rdr,
\]
with $D_r=\partial_r+\frac1r$ the radial divergence and $D_r^{-1}z(r)=\frac1r\int_0^r sz(s)\,ds$ a sort of Hardy inverse, both of which are bounded by the energy. The second one is bounded thanks to the Hardy inequality \eqref{eq:radial inverse Hardy}. The first modification removes the tangential weights in front of the ``derivative'' terms and the second removes those in front of the ``mass'' terms. The functional 
\[
\mathcal H_\phi=I_\phi-\frac12\mathcal J_\phi-\frac12\mathcal K_\phi
\]
 then satisfies a virial identity whose principal part is represented by the radial weight $\phi_r$, together with other two nonnegative terms. We emphasize that these cancellations are exact and rely on the algebraic structure of the system, which are represented by algebraic-differential relations between the scalar and vector radial operators and some new nonlocal resolvent operators $K_s$, $K_v$. They do not rely on the choice of the weight.  The resulting quadratic form still contains nonlocal terms, which are indefinite in general. Following \cite{KMPP2019}, we pass to the canonical variables $f=K_vq$, $g=K_sh$, in which the quadratic form becomes a local expression in $f,g$ and their derivatives up to order three. Its positivity is then studied in detail. After a suitable completion of squares, it is proved that this form is coercive if and only if $-a>\frac19$ and $-c>\frac19$, in the normalized variables $b=d=1$. The threshold $\frac19$ is sharp for the functional, in the sense that for constant weights the form reduces, in Fourier variables, to a symbol which is nonnegative precisely in the considered region of coercivity. 

Now we provide a choice for the weight $\phi$. The weights used in dimensions 2 and 3 in \cite{HMM1,HMM2} have inspired our choice. The weight must satisfy several necessary requirements, such as 
\[
\phi_r>0, \quad \phi/r\ge\phi_r, \quad (\phi/r)_{rr}\ge\frac1r(\phi/r)_r.
\]
Additionally, the ratio $\phi/(r\phi_r)$, which measures the strength of the tangential weight relative to the radial one, must be controlled when estimating the cubic terms (recall that no $L^\infty$ control is allowed in 2D). We have chosen 
\[
\phi_r^L(r)= \frac1{L(1+ \frac{r^2}{L^2})^{3/4}}, 
\]
so that $\phi/(r\phi_r)\lesssim(1+r/L)^{1/2}$. The exponent $\frac34$ is the largest one for which this growth is compensated by the pointwise decay $|q(r)|+|h(r)|\lesssim\varepsilon r^{-1/2}$ of radial $H^1$ profiles (a.k.a. the Strauss lemma \cite{Strauss1977}). The resulting local energy can be described as 
\[
\mathcal E_L(t) \sim \int\langle x/L\rangle^{-3/2}(|u|^2+|Du|^2+\eta^2+|\nabla\eta|^2)\,dx, \quad \langle x\rangle =\sqrt{1+x^2},
\]
whose weight is not integrable in $\mathbb R^2$. Finally, the cubic terms are shown to be of size $\varepsilon$ times the quadratic form, which closes the argument and yields the classical time integrated boundedness 
\[
\int_{\mathbb R}\mathcal E_L(t)\,dt\lesssim\varepsilon^2.
\]
The desired decay of $\mathcal E_L(t)$ follows from a second virial-type identity for a suitably modified weighted energy, whose time derivative is again controlled by $\mathcal E_L$. The mechanism introduced here, namely the compensation of tangential terms by functionals built on the radial divergence and its inverse, seems of independent interest and should apply to other two-dimensional dispersive systems under radial symmetry. Removing the radial symmetry, which in our argument is needed both to reduce the problem to profiles and to obtain the pointwise decay, and extending the range of parameters beyond \eqref{eq:parameter condition}, will be considered elsewhere.

\subsection{Notations}\label{sec:notations}

For any positive numbers $\alpha$ and $\beta$, $\alpha \lesssim \beta$ means that there exists $C>0$ such that $\alpha \le C\beta$, and $\alpha \sim \beta$ means $\alpha \lesssim \beta$ and $\beta \lesssim \alpha$. We regard vectors in $\R^n$ as column vectors. We denote the horizontal variable by $x=(x_1,\ldots,x_n)$, and the Euclidean scalar product (inner product) by $\cdot$, that is to say that $x\cdot y=\sum_{j=1}^n x_jy_j$ for $x=(x_1,\ldots,x_n)$ and $y=(y_1,\ldots,y_n)$ in $\R^n$. $|\cdot|$ denotes the Euclidean norm, i.e., $|x|=\left(\sum_{j=1}^n x_j^2\right)^{\frac12}$.

We conveniently use the sub-index notation for partial derivatives as follows: for $v=(v_1,\ldots,v_n)$ and $f$ defined on $\R^n$, we write
\[\partial_{x_j}v_k=v_{k,x_j}\quad\mbox{and}\quad\partial_{x_j}f=f_{x_j},\qquad j,k=1,\ldots,n.\]

For $f\in\mathcal S'(\R^n)$ we denote by $\mathcal F(f)$ or $\widehat{f}$ the Fourier transform of $f$ with respect to the spatial variable,
\[\hat{f}(\xi)=\int_{\R^n}e^{-i(x\cdot\xi)}f(x)\;dx,\qquad \xi=(\xi_1,\ldots,\xi_n)\in\R^n.\]
Under suitable condition, $f$ is determined by its Fourier coefficient $\widehat{f}$ via the inverse Fourier transform (denoted by $\mathcal F^{-1}$):
\[f(x)=\mathcal F^{-1}(\widehat{f})(x):=\frac{1}{(2\pi)^n}\int_{\R^n}e^{i(x\cdot\xi)}\widehat{f}(\xi)\;d\xi.\]
We define Bessel and Riesz potentials ($J^s$ and $D^s$, respectively) of order $-s$, $s\in\R$, as Fourier multipliers by
\[J^sf(x):=\mathcal F^{-1}\left((1+|\xi|^2)^{\frac{s}{2}}\widehat{f}\right)\quad\mbox{and}\quad D^sf(x):=\mathcal F^{-1}\left(|\xi|^s\widehat{f}\right).\]
In particular, $\sqrt{-\Delta}=D^1=D$ is the Fourier multiplier of the symbol $|\xi|$.

We denote by $\norm{\cdot}_{L^p}$ ($1\le p\le\infty$) the standard norm of Lebesgue space $L^p(\R^n)$, and the standard notation $H^s(\R^n)$, simply $H^s$, denotes the $L^2$-based Sobolev space equipped with the norm
\[\norm{f}_{H^s}^2=\int_{\R^n}(1+|\xi|^2)^s|\wh{f}(\xi)|^2\;d\xi.\]
For $v=(v_1,\ldots,v_n)\in H^s(\R^n)^n$ and $f\in H^s(\R^n)$, we denote a $H^s$-level space of $(v,f)$ by $\vH^s(\R^n)$ (as an abbreviation for the product space), simply $\vH^s$, equipped with the norm
\[\norm{(v,f)}_{\vH^s}^2=\sum_{j=1}^n\norm{v_j}_{H^s}^2+\norm{f}_{H^s}^2.\]
From the energy \eqref{ENERGY}, $\vH^1(\R^n)$ is the natural energy space, and Theorem \ref{thm:GWP} establishes the global well-posedness of \eqref{nD ABCD} with sufficiently small initial data in the energy space $\vH^1(\R^n)$ for $2\le n\le4$.

Finally, let $a_1,a_2, a_3 \in\R$. The quantities $a_{max}\ge a_{med}\ge a_{min}$ can be conveniently defined to be the maximum, median and minimum values of $a_1,a_2,a_3$ respectively.

\subsection{Acknowledgments}
C. M. would like to thank support and hospitality during April and May 2026 at Ewha Womans U. in South Korea, where part of this work was done.  

\subsection{Organization of this work} This paper is organized as follows. Section \ref{sec:pre} presents preliminary results needed for the rest of the proof. Section \ref{sec:main} is devoted to the proof of local and global well-posedness (Theorems \ref{thm:LWP} and \ref{thm:GWP}). Section \ref{sec:4} contains the proof of the virial identities related to the \emph{abcd} dynamics. Section \ref{sec:quad} contains the analysis of quadratic terms appearing in virial identities. Finally, Section \ref{sec:MT} contains the end of the proof of Theorem \ref{thm:MT}, and in particular Subsection~\ref{sec:expanding} contains the proof of Corollary~\ref{cor:expanding}.

\section{Preliminaries}\label{sec:pre}

\subsection{Structural equivalence}\label{sec:reformualtion}
Analogously to \cite{DMS2007}, we write \eqref{nD ABCD} in the equivalent form (in the Fourier space, for fixed $\xi\in\R^n\setminus\{0\}$)
\[
\partial_t
\begin{pmatrix}
\wh{\eta}\\
\wh{u}
\end{pmatrix}
+i|\xi|\mathcal A(\xi)
\begin{pmatrix}
\wh{\eta}\\
\wh{u}
\end{pmatrix}
+i
\begin{pmatrix}
\displaystyle\frac{\sum_{j=1}^n\xi_j\wh{u_j\eta}}{1+b|\xi|^2}\\[2mm]
\displaystyle\frac{\xi\wh{|u|^2}}{2(1+d|\xi|^2)}
\end{pmatrix}
=0,
\]
where
\[\mathcal A(\xi)=
\begin{pmatrix}
0 & \displaystyle\frac{1-a|\xi|^2}{1+b|\xi|^2}\frac{\xi^T}{|\xi|}\\[2mm]
\displaystyle\frac{1-c|\xi|^2}{1+d|\xi|^2}\frac{\xi}{|\xi|} & 0_{n\times n}
\end{pmatrix}.
\]
Define the Fourier symbol $\sigma(\xi)$ by
\begin{equation}\label{eq:eigenvalues}
\sigma(\xi)=\left(\frac{(1-a|\xi|^2)(1-c|\xi|^2)}{(1+b|\xi|^2)(1+d|\xi|^2)}\right)^{\frac12}.
\end{equation}
Note that $|\sigma(\xi)|\sim1$ since $a,c<0$ and $b,d>0$. A direct computation shows that $\mathcal A(\xi)$ has the eigenvalues $0$, with multiplicity $n-1$, and $\sigma(\xi)$ and $-\sigma(\xi)$. More precisely, regarding $\mathcal A(\xi)$ as acting on $\C^{n+1}$, let
\[\mathcal A(\xi)\begin{pmatrix}z_0\\z\end{pmatrix}=0,\qquad z_0\in\C,\quad z\in\C^n.\]
Since $a,c<0$, $b,d>0$, and $\xi\neq0$, the vector
\[\frac{1-c|\xi|^2}{1+d|\xi|^2}\frac{\xi}{|\xi|}\]
is nonzero, and hence the second component implies $z_0=0$. The first component then implies that both $\operatorname{Re}z$ and $\operatorname{Im}z$ are orthogonal to $\xi$. Thus, the eigenspace corresponding to the zero eigenvalue consists precisely of vectors of the form
\[\begin{pmatrix}0\\z\end{pmatrix},\qquad z\in\C^n,\]
with $\operatorname{Re}z\perp\xi$ and $\operatorname{Im}z\perp\xi$.

We now use the irrotational condition \eqref{eq:irrotational condition}. In the Fourier variables, it gives
\[\xi_j\wh{u_k}(\xi)=\xi_k\wh{u_j}(\xi),\qquad 1\le j,k\le n.\]
For each fixed $\xi\in\R^n\setminus\{0\}$, choose $k$ such that $\xi_k\neq0$. Then
\[\wh{u_j}(\xi)=\frac{\xi_j}{\xi_k}\wh{u_k}(\xi),\qquad j=1,\ldots,n,\]
and hence
\[\wh{u}(\xi)=\frac{\wh{u_k}(\xi)}{\xi_k}\,\xi.\]
Thus, $\wh{u}(\xi)$ is parallel to $\xi$ for every $\xi\in\R^n\setminus\{0\}$. Since the zero eigenspace consists of the velocity directions orthogonal to $\xi$, $\wh{u}(\xi)$ has no component in the zero eigenspace.

We define the longitudinal component of $\wh{u}(\xi)$ by
\[\ell(\xi):=\frac{\xi\cdot\wh{u}(\xi)}{|\xi|}.\]
Since $\wh{u}(\xi)$ is parallel to $\xi$, we have
\[\wh{u}(\xi)=\left(\frac{\xi\cdot\wh{u}(\xi)}{|\xi|}\right)\frac{\xi}{|\xi|}=\frac{\xi}{|\xi|}\ell(\xi).\]
Substituting this representation into the action of $\mathcal A(\xi)$, we obtain
\[\mathcal A(\xi)\begin{pmatrix}\wh{\eta}\\ \wh{u}\end{pmatrix}
=
\begin{pmatrix}
\displaystyle\frac{1-a|\xi|^2}{1+b|\xi|^2}\ell\\[2mm]
\displaystyle\frac{\xi}{|\xi|}\frac{1-c|\xi|^2}{1+d|\xi|^2}\wh{\eta}
\end{pmatrix}.
\]
Thus, under the irrotational condition, the action of $\mathcal A(\xi)$ on $(\wh{\eta},\wh{u})$ is completely determined by the scalar pair $(\wh{\eta},\ell)$, which is transformed according to
\[
\begin{pmatrix}
\wh{\eta}\\
\ell
\end{pmatrix}
\longmapsto
\begin{pmatrix}
0 & \displaystyle\frac{1-a|\xi|^2}{1+b|\xi|^2}\\[2mm]
\displaystyle\frac{1-c|\xi|^2}{1+d|\xi|^2} & 0
\end{pmatrix}
\begin{pmatrix}
\wh{\eta}\\
\ell
\end{pmatrix}.
\]
This $2\times2$ matrix is diagonalized as
\[
\mathcal P^{-1}(\xi)
\begin{pmatrix}
0 & \displaystyle\frac{1-a|\xi|^2}{1+b|\xi|^2}\\[2mm]
\displaystyle\frac{1-c|\xi|^2}{1+d|\xi|^2} & 0
\end{pmatrix}
\mathcal P(\xi)
=
\begin{pmatrix}
\sigma(\xi) & 0\\
0 & -\sigma(\xi)
\end{pmatrix},
\]
where
\[
\mathcal P(\xi)=
\begin{pmatrix}
\varsigma(\xi) & -\varsigma(\xi)\\
1 & 1
\end{pmatrix},
\qquad
\mathcal P^{-1}(\xi)=\frac12
\begin{pmatrix}
\varsigma(\xi)^{-1} & 1\\
-\varsigma(\xi)^{-1} & 1
\end{pmatrix},
\]
with
\[\varsigma(\xi)=\left(\frac{(1-a|\xi|^2)(1+d|\xi|^2)}{(1+b|\xi|^2)(1-c|\xi|^2)}\right)^{\frac12}.\]
Note that $|\varsigma(\xi)|\sim1$ since $a,c<0$ and $b,d>0$.

Accordingly, we define the diagonal variables by
\[\begin{pmatrix}
\wh{\nu}_1\\
\wh{\nu}_2
\end{pmatrix}
=
\mathcal P^{-1}(\xi)
\begin{pmatrix}
\wh{\eta}\\
\ell
\end{pmatrix},\]
that is,
\[\wh{\nu}_1=\frac{\wh{\eta}}{2\varsigma(\xi)}+\frac12\ell,\qquad \wh{\nu}_2=-\frac{\wh{\eta}}{2\varsigma(\xi)}+\frac12\ell.\]
Taking the sum and the difference, we obtain
\[\wh{\eta}=\varsigma(\xi)\left(\wh{\nu}_1-\wh{\nu}_2\right),\qquad \ell=\wh{\nu}_1+\wh{\nu}_2.\]
Since
\[\wh{u}(\xi)=\frac{\xi}{|\xi|}\ell(\xi),\]
the inverse change of variables is
\[\wh{\eta}=\varsigma(\xi)\left(\wh{\nu}_1-\wh{\nu}_2\right),\qquad \wh{u_j}=\frac{\xi_j}{|\xi|}\left(\wh{\nu}_1+\wh{\nu}_2\right),\qquad j=1,\ldots,n.\]

We introduce the $j$-th Riesz transform $\mathcal R_j$ by
\[\widehat{\mathcal R_jf}(\xi):=-i\frac{\xi_j}{|\xi|}\wh f(\xi),\qquad j=1,\ldots,n,\]
and the Fourier multiplier $\varsigma(D)$ by
\[\widehat{\varsigma(D)f}(\xi):=\varsigma(\xi)\wh f(\xi).\]
Then
\begin{equation}\label{eq:eta u inverse}
\eta=\varsigma(D)(\nu_1-\nu_2),\qquad u_j=i\mathcal R_j(\nu_1+\nu_2),\qquad j=1,\ldots,n.
\end{equation}
Since $\varsigma(\xi)\sim1$ and the Riesz transforms are bounded on $H^s(\R^n)$, under the irrotational condition, $(\eta,u)\in H^s(\R^n)\times H^s(\R^n)^n$ if and only if $(\nu_1,\nu_2)\in H^s(\R^n)^2$, with equivalent norms.

Under the irrotational condition, the inverse transformation \eqref{eq:eta u inverse} recovers $(u,\eta)$ uniquely from $(\nu_1,\nu_2)$. Therefore, the dynamics of irrotational solutions to \eqref{nD ABCD} can be reformulated entirely in terms of $(\nu_1,\nu_2)$ as
\begin{equation}\label{eq:nd abcd}
\partial_t \begin{pmatrix}\nu_1 \\ \nu_2 \end{pmatrix}+\mathcal B(-i\nabla)\begin{pmatrix}\nu_1 \\ \nu_2 \end{pmatrix}+\mathcal N\begin{pmatrix}\nu_1 \\ \nu_2 \end{pmatrix}=0,\qquad \begin{pmatrix}\nu_1(0)\\ \nu_2(0)\end{pmatrix}=\begin{pmatrix}\nu_{1,0}\\ \nu_{2,0}\end{pmatrix}.
\end{equation}
where $\mathcal B(-i\nabla)$ is the $2\times2$ matrix Fourier multiplier with symbol
\[i|\xi|
\begin{pmatrix}
\sigma(\xi)&0\\
0&-\sigma(\xi)
\end{pmatrix},
\]
which is skew-Hermitian.

Using the inverse transformation \eqref{eq:eta u inverse},
\[u_j\eta=i\mathcal R_j(\nu_1+\nu_2)\varsigma(D)(\nu_1-\nu_2),\qquad |u|^2=-\sum_{j=1}^n\left(\mathcal R_j(\nu_1+\nu_2)\right)^2.\]
Hence
\begin{equation}\label{eq:nonlinear}
\mathcal N
\begin{pmatrix}
\nu_1\\
\nu_2
\end{pmatrix}
=\frac{i}{2}
\begin{pmatrix}
\mathcal N_+(\nu_1,\nu_2)\\
\mathcal N_-(\nu_1,\nu_2)
\end{pmatrix},
\end{equation}
where
\[\begin{aligned}
\mathcal N_\pm(\nu_1,\nu_2):=&~{}\mp\varsigma(D)^{-1}D(1-b\Delta)^{-1}\sum_{j=1}^n\mathcal R_j\left(\mathcal R_j(\nu_1+\nu_2)\varsigma(D)(\nu_1-\nu_2)\right)\\
&~{}-\frac12D(1-d\Delta)^{-1}\sum_{j=1}^n\left(\mathcal R_j(\nu_1+\nu_2)\right)^2.
\end{aligned}\]

\subsection{Littlewood-Paley Decomposition}
This section is devoted to explaining the Littlewood-Paley decomposition, which is the main tool in the proof of the bilinear estimate for the local well-posedness theory. As is well known, the Littlewood-Paley decomposition is a particular way to write a single function as a superposition of a countably infinite family of functions of varying frequencies.

\medskip

Let $\varphi(\xi)$ be a real-valued radially symmetric bump function on $\R^n$ supported in $\{\xi \in \R^n : |\xi| \le 2\}$, which is identical to $1$ on the set $\{\xi \in \R^n : |\xi| \le 1\}$ and is decreasing on $\{\xi \in \R^n : 1 \le |\xi| \le 2\}$. Define a dyadic number $N \in 2^{\Z_{\ge 0}}$ of the form $N = 2^k$, $k \in \Z_{\ge 0}$. Let $\varphi_1 = \varphi$ and define
\[\varphi_N(\xi) = \varphi\left(\frac{\xi}{N}\right) - \varphi\left(\frac{2\xi}{N}\right), \quad N \ge 2.\]
By construction, the sequence of $\varphi_N$ satisfies
\[\sum_{N \ge 1 \; : \; dyadic} \varphi_N(\xi) \equiv 1.\]
We simply write $\displaystyle \sum_{N \ge 1}$ by dropping ``dyadic". This provides a typical \emph{partition of unity} which allows to define the projection operator (one of the so-called \emph{Littlewood-Paley} projection operators) on $L^2(\R^n)$ by $P_Nf(x) = \mathcal F^{-1} \left(\varphi_N(\xi) \widehat{f} \right) (x)$. Using the projection operators, one decomposes any function $f$ in $L^2(\R^n)$ as
\[f = \sum_{N \ge 1} P_N f.\]
We sometimes denote $P_Nf$ by simply $f_N$. Note that for every $N$, $f_N\in H^\infty$, whenever $f \in L^2(\R^n)$.

\medskip

We next record a basic support property of dyadic frequency interactions. Let $f_i\in L^2(\R^n)$ be compactly supported functions with $\operatorname{supp}f_i\subset\{\xi\in\R^n:|\xi|\sim N_i\}$, $i=1,2,3$, with the convention that $|\xi|\sim 1$ means $|\xi|\lesssim1$ for the lowest dyadic block. We define
\[J(f_1,f_2,f_3) = \int_{\R^n} f_1(\xi_1)f_2(\xi_2)f_3(\xi_1+\xi_2) \; d\xi_1d\xi_2.\]
A direct change of variables in the integral yields
\[J(f_1,f_2,f_3) = J(f_2,f_1,f_3) = J(\wt{f}_1,f_3,f_2),\]
where $\wt{f}(\xi)=f(-\xi)$. On the support of the integrand in $J(f_1,f_2,f_3)$, we have
\[\xi_1+\xi_2=\xi_3.\]
Hence $J(f_1,f_2,f_3)$ vanishes unless
\begin{equation}\label{eq:support property}
N_{\max}\sim N_{\med} \gtrsim 1.
\end{equation}

The following lemma is the well-known \emph{Bernstein inequality}, which upgrades low to high Lebesgue integrability with the price of some powers of $N$:
\begin{lemma}[Bernstein's inequalities]
Let $f \in L^2(\R^n)$, $1 \le p \le q \le \infty$, and $s \ge 0$. Then, we have
\begin{equation}\label{eq:Bernstein_1}
\norm{D^{\pm s}f_N}_{L^p} \sim N^{\pm s} \norm{f_N}_{L^p}, \quad N\geq 2, \quad \norm{D^{s}f_1}_{L^p} \lesssim  \norm{f_1}_{L^p},  
\end{equation}
and
\begin{equation}\label{eq:Bernstein_2}
\norm{f_N}_{L^q} \lesssim N^{\frac{n}{p} - \frac{n}{q}} \norm{f_N}_{L^p}.
\end{equation}
The implicit constants in both \eqref{eq:Bernstein_1} and \eqref{eq:Bernstein_2} depend only on $s, p, q$ and $n$.
\end{lemma}

\section{Global well-posedness}\label{sec:main}

\subsection{Local well-posedness}
As mentioned in Section \ref{sec:WP theory}, for the proof of Theorem \ref{thm:LWP}, we prove Proposition \ref{prop:refined bilinear} when $n\ge3$. The case $n=2$ follows from Lemma \ref{lem:bi0}.

\begin{proof}[Proof of Proposition \ref{prop:refined bilinear}]
It suffices to prove the following separately:
\begin{equation}\label{eq:low}
\norm{(1-\Delta)^{-1}DP_1(fg)}_{H^s}\lesssim\norm{f}_{H^s}\norm{g}_{H^s}
\end{equation}
and
\begin{equation}\label{eq:high}
\norm{(1-\Delta)^{-1}DP_{> 1}(fg)}_{H^s}\lesssim\norm{f}_{H^s}\norm{g}_{H^s},
\end{equation}
where $N >1$.

\medskip

The proof of \eqref{eq:low} immediately follows from the boundedness of the multiplier $(1+|\xi|^2)^{\frac{s-2}{2}}|\xi|$ on $|\xi|\le2$, together with
\[|(\wh{f}\ast\wh{g})(\xi)|\lesssim\norm{f}_{L^2}\norm{g}_{L^2} \quad \text{and} \quad \int_{|\xi|\le2}\;d\xi\lesssim1.\]

\medskip

For the proof of \eqref{eq:high}, we use the duality argument. Then it suffices to estimate
\[\int_{\R^n} \left(J^{s-2}D(fg) \right) P_{>1}w \; dx,\]
where $w \in L^2$ with $\norm{w}_{L^2} =1$. We make the Littlewood-Paley decomposition of $f,g$ and $w$ as
\[f = \sum_{N_1 \ge 1} f_{N_1}, \quad g = \sum_{N_2 \ge 1}g_{N_2} \quad \mbox{and} \quad P_{>1}w = \sum_{N > 1} w_N,\]
respectively. Without loss of generality, we may assume $N_1 \le N_2$. Then we are now reduced to establishing
\begin{equation}\label{eq:bi2}
\sum_{N > 1}\sum_{N_1, N_2 \ge 1}N^{-1+s}\int f_{N_1}g_{N_2} w_N \; dx \lesssim \norm{f}_{H^s}\norm{g}_{H^s}.
\end{equation}
From \eqref{eq:support property}, we separate \eqref{eq:bi2} into two cases: $N \sim N_2 \gtrsim N_1$ and $N_2 \sim N_1 \gg N$.

\medskip

\textbf{(Case I.)} $N \sim N_2 \gtrsim N_1$. Using H\"older's and Bernstein's \eqref{eq:Bernstein_2} inequalities, one has
\[\int f_{N_1}g_{N_2} w_N \; dx \le \norm{f_{N_1}}_{L^{\infty}}\norm{g_{N_2}}_{L^2}\norm{w_N}_{L^2} \lesssim N_1^{\frac n2}\norm{f_{N_1}}_{L^2}\norm{g_{N_2}}_{L^2}\norm{w_N}_{L^2}.\]
With this, we further reduce \eqref{eq:bi2} to
\begin{equation}\label{eq:bi3}
\sum_{N > 1}\sum_{N_2 \sim N}\sum_{N_1 \lesssim N_2} N^{-1+s}N_2^{-s}N_1^{\frac n2-s} \norm{f_{N_1}}_{H^s}\norm{g_{N_2}}_{H^s}\norm{w_N}_{L^2} \lesssim \norm{f}_{H^s}\norm{g}_{H^s},
\end{equation}
thanks to \eqref{eq:Bernstein_1}. We denote the multiplier in the left-hand side of \eqref{eq:bi3} by $m(N,N_1,N_2)$, namely, $m(N,N_1,N_2) = N^{-1+s}N_2^{-s}N_1^{\frac n2-s}$. Note that the number of $N_2$ satisfying $N_2\sim N$ is uniformly finite, and hence
\[\sum_{N}\sum_{N_2 \sim N} \norm{g_{N_2}}_{H^s}^2 \sim \sum_{N}\norm{g_N}_{H^s}^2.\]
With this observation, the Cauchy-Schwarz inequality yields
\[\mbox{LHS of } \eqref{eq:bi3} \lesssim \left(\sup_{N > 1} \sum_{\substack{1 \le N_1 \lesssim N_2 \\
N_2 \sim N}} m^2(N,N_1,N_2) \right)^{\frac12}\norm{f}_{H^s}\norm{g}_{H^s}\norm{w}_{L^2}.\]
Hence, it remains to prove
\begin{equation}\label{eq:bi4}
\sup_{N > 1} \sum_{N_2 \sim N}\sum_{1 \le N_1 \lesssim N_2} N^{-2+2s}N_2^{-2s}N_1^{n-2s} \lesssim 1.
\end{equation}
This is sometimes referred to as \emph{Schur's test}; see, for instance, \cite[Lemma 3.11]{Tao2001} for more details. If $s < \frac n2$, one gets
\[\sum_{N_2 \sim N}\sum_{1 \le N_1 \lesssim N_2} N^{-2+2s}N_2^{-2s}N_1^{n-2s} \lesssim \sum_{N_2 \sim N}N^{-2+2s}N_2^{n-4s} \lesssim N^{n-2-2s} \lesssim 1,\]
provided that $s \ge \frac{n-2}{2}$. If $s=\frac n2$, then
\[\sum_{N_2 \sim N}\sum_{1 \le N_1 \lesssim N_2} N^{-2+2s}N_2^{-2s}N_1^{n-2s} \lesssim N^{-2}\log N \lesssim 1.\]
Finally, if $s>\frac n2$, then
\[\sum_{N_2 \sim N}\sum_{1 \le N_1 \lesssim N_2} N^{-2+2s}N_2^{-2s}N_1^{n-2s} \lesssim \sum_{N_2 \sim N}N^{-2+2s}N_2^{-2s} \lesssim N^{-2} \lesssim 1.\]
Thus, \eqref{eq:bi4} holds true if $s \ge \frac{n-2}{2}$.

\medskip

\textbf{(Case II.)} $N_2 \sim N_1 \gg N$. Analogously, using H\"older's and Bernstein's \eqref{eq:Bernstein_2} inequalities, one obtains
\[\int f_{N_1}g_{N_2} w_N \; dx \le \norm{f_{N_1}}_{L^2}\norm{g_{N_2}}_{L^2}\norm{w_N}_{L^{\infty}} \lesssim N^{\frac n2}\norm{f_{N_1}}_{L^2}\norm{g_{N_2}}_{L^2}\norm{w_N}_{L^2},\]
which ensures \eqref{eq:high} provided that
\begin{equation}\label{eq:bi5}
\sup_{N_2 > 1} \sum_{N_1 \sim N_2}\sum_{1 < N \lesssim N_1} N_2^{-2s}N_1^{-2s}N^{n-2+2s} \lesssim 1.
\end{equation}
Note that for $s \ge \frac{n-2}{2}$ and $n\ge3$, we have $n-2+2s>0$. Thus,
\[\sum_{N_1 \sim N_2}\sum_{1 < N \lesssim N_1} N_2^{-2s}N_1^{-2s}N^{n-2+2s} \lesssim N_2^{n-2-2s} \lesssim 1.\]
Therefore, \eqref{eq:bi5} holds for $s\ge\frac{n-2}{2}$ if $n\ge3$. This completes the proof.
\end{proof}

\begin{remark}
In dimension two, Case II is precisely the obstruction to extending Proposition \ref{prop:refined bilinear} to the endpoint $s=0$, since the corresponding dyadic summation exhibits a logarithmic divergence. Indeed, when $n=2$ and $s=0$, one has $N^{n-2+2s}=1$ in the summand on the left-hand side of \eqref{eq:bi5}, and hence
\[\sum_{1<N\lesssim N_1}N^{n-2+2s}=\sum_{1<N\lesssim N_1}1\sim\log N_1.\] 
On the other hand, the high-high to low interactions arising in this case are strictly non-resonant for the considered dispersive system. This suggests that the present logarithmic divergence at $s=0$ might be overcome by exploiting the oscillatory structure of the Duhamel term. 
\end{remark}

\begin{remark}\label{rem:persistence}
For $2 \le n \le 4$ and $s \ge 1$, we indeed have
\[\norm{(1-\Delta)^{-1}D(fg)}_{H^s(\R^n)} \lesssim \norm{f}_{H^1(\R^n)}\norm{g}_{H^s(\R^n)} + \norm{f}_{H^s(\R^n)}\norm{g}_{H^1(\R^n)}.\]

The proof follows the same Littlewood--Paley decomposition as in the proof of Proposition \ref{prop:refined bilinear}, with one important difference. Since the regularities assigned to the two factors are different, we split the summation into the regions $N_1\le N_2$ and $N_2<N_1$. In the first region, the low-frequency factor is measured in $H^1$ and the high-frequency factor in $H^s$, while the second region gives the symmetric term.

For $N_1\le N_2$, in Case I, where $N\sim N_2\gtrsim N_1$, the multiplier in \eqref{eq:bi3} is replaced by
\[N^{-1+s}N_2^{-s}N_1^{\frac n2-1},\]
and the corresponding summation is bounded by
\[\sup_{N>1}N^{-2}\sum_{1\le N_1\lesssim N}N_1^{n-2}\lesssim1\]
for $2\le n\le4$. In Case II, where $N_1\sim N_2\gg N$, the corresponding summation is bounded by
\[\sup_{N_2>1}N_2^{-2s-2}\sum_{1<N\lesssim N_2}N^{n-2+2s}\lesssim\sup_{N_2>1}N_2^{n-4}\lesssim1.\]
The low-frequency estimate is treated exactly as in \eqref{eq:low}.
\end{remark}

Proposition \ref{prop:refined bilinear}, together with a standard fixed-point  argument, yields Theorem \ref{thm:LWP}.

\begin{proof}[Proof of Theorem \ref{thm:LWP}]
By \eqref{eq:eta u inverse} and the fact that $|\varsigma(\xi)|\sim1$, together with the boundedness of the Riesz transforms on $H^s(\R^n)$, one has
\[(u,\eta)\in\vH^s(\R^n)\quad\Longleftrightarrow\quad(\nu_1,\nu_2)\in H^s(\R^n)\times H^s(\R^n),\]
with equivalent norms. Thus, it suffices to prove the local well-posedness of \eqref{eq:nd abcd}. Throughout the proof, the implicit constants may depend on the fixed parameters $a,b,c,d$, as well as on $n$ and $s$, but are independent of the solution and the existence time $T$.

\medskip

We recall from \eqref{eq:nonlinear} that the nonlinear part is given by
\[\mathcal N\binom{\nu_1}{\nu_2}=\frac i2\binom{\mathcal N_+(\nu_1,\nu_2)}{\mathcal N_-(\nu_1,\nu_2)}.\]
The nonlinearity $\mathcal N(\nu_1,\nu_2)$ involves the multipliers $\varsigma(D)^{\pm1}$ and $\mathcal R_j$, while the multipliers $D(1-b\Delta)^{-1}$ and $D(1-d\Delta)^{-1}$ are of the same order as $J^{-2}D$. Thus, applying Proposition \ref{prop:refined bilinear}, one obtains
\begin{equation}\label{eq:nonlinear estimates}
\norm{\mathcal N(\nu_1,\nu_2)}_{H^s\times H^s}\le C\norm{(\nu_1,\nu_2)}_{H^s\times H^s}^2.
\end{equation}

We now complete the proof of Theorem \ref{thm:LWP}. Let
\[\norm{(\nu_{1,0},\nu_{2,0})}_{H^s\times H^s}\le R,\]
for $R>0$ arbitrary but fixed. By Duhamel's principle, \eqref{eq:nd abcd} is equivalent to the following integral equation:
\begin{equation}\label{Duhamel}
\begin{pmatrix}\nu_1\\ \nu_2\end{pmatrix}=S(t)\begin{pmatrix}\nu_{1,0}\\ \nu_{2,0}\end{pmatrix}-\int_0^tS(t-\tau)\mathcal N(\nu_1,\nu_2)(\tau)\;d\tau,
\end{equation}
where $S(t)$ denotes the linear propagator associated with the linear operator $\mathcal B(-i\nabla)$, that is,
\[S(t)\begin{pmatrix}f\\ g\end{pmatrix}=\begin{pmatrix}\mathcal F^{-1}\left(e^{-it|\xi|\sigma(\xi)}\wh f(\xi)\right)\\ \mathcal F^{-1}\left(e^{it|\xi|\sigma(\xi)}\wh g(\xi)\right)\end{pmatrix},\]
where $\sigma(\xi)$ is defined in \eqref{eq:eigenvalues}. Note that $e^{\mp it|D|\sigma(D)}$ are unitary operators on $L^2$-based Sobolev spaces. Let $X_T^s$ denote $C([0,T];H^s(\R^n))\times C([0,T];H^s(\R^n))$ equipped with the norm
\[\norm{(f,g)}_{X_T^s}^2=\norm{f}_{L_T^\infty H^s}^2+\norm{g}_{L_T^\infty H^s}^2,\]
for some $T>0$ to be chosen later. Taking the $X_T^s$-norm on both sides of \eqref{Duhamel} and applying \eqref{eq:nonlinear estimates}, one has
\begin{equation}\label{eq:contraction1}
\norm{(\nu_1,\nu_2)}_{X_T^s}\le\norm{(\nu_{1,0},\nu_{2,0})}_{H^s\times H^s}+CT\norm{(\nu_1,\nu_2)}_{X_T^s}^2.
\end{equation}
Moreover, for different functions $(\wt{\nu}_1,\wt{\nu}_2)$, one similarly has
\begin{equation}\label{eq:contraction2}
\begin{aligned}
&~{}\norm{(\nu_1,\nu_2)-(\wt{\nu}_1,\wt{\nu}_2)}_{X_T^s}\\
\le&~{}\norm{(\nu_{1,0}-\wt{\nu}_{1,0},\nu_{2,0}-\wt{\nu}_{2,0})}_{H^s\times H^s}\\
&~{}+CT\left(\norm{(\nu_1,\nu_2)}_{X_T^s}+\norm{(\wt{\nu}_1,\wt{\nu}_2)}_{X_T^s}\right)\norm{(\nu_1-\wt{\nu}_1,\nu_2-\wt{\nu}_2)}_{X_T^s}.
\end{aligned}
\end{equation}
We define the iteration map by the right-hand side of \eqref{Duhamel}. Choosing $T>0$ in both \eqref{eq:contraction1} and \eqref{eq:contraction2} satisfying
\[4CRT<\frac12,\]
one concludes that this map is a contraction from the ball
\[\left\{(\nu_1,\nu_2)\in X_T^s:\norm{(\nu_1,\nu_2)}_{X_T^s}\le 2R\right\}\]
to itself, which completes the proof of Theorem \ref{thm:LWP}.
\end{proof}

\subsection{Global well-posedness}\label{sec:GWP}
This section is devoted to the proof of Theorem \ref{thm:GWP}. We first establish a Sobolev estimate to control the cubic term in the energy \eqref{ENERGY}. For $1\le n\le6$, recall the classical Sobolev inequality $H^1(\R^n)\hookrightarrow L^3(\R^n)$, that is, there exists $C_0=C_0(n)>0$ such that, for all $f\in H^1(\R^n)$,
\be\label{Sob}
\norm{f}_{L^3(\R^n)}\le C_0\norm{f}_{H^1(\R^n)}.
\ee

\begin{lemma}\label{lem:L3}
Let $1\le n\le6$ and $\alpha>0$ be arbitrary. Then, we have
\[\left|\int_{\R^n}\eta|u|^2\;dx\right|\le\frac{1}{2\alpha}\norm{(u,\eta)}_{\vH^1}^2+\frac{\alpha C_0^6}{2}\norm{(u,\eta)}_{\vH^1}^4,\]
where $C_0$ is the constant appearing in the Sobolev inequality \eqref{Sob}.
\end{lemma}

\begin{proof}
The proof follows from H\"older's, Minkowski's, Young's, and Sobolev inequalities. Indeed, one has
\[\begin{aligned}
\left|\int_{\R^n}\eta|u|^2\;dx\right| \le&~{}\norm{\eta}_{L^3}\norm{\sum_{j=1}^n u_j^2}_{L^{\frac32}}\\
\le&~{}\frac{\delta}{2}\norm{\eta}_{L^3}^2+\frac{1}{2\delta}\left(\sum_{j=1}^n\norm{u_j}_{L^3}^2\right)^2\\
\le&~{}\frac{\delta C_0^2}{2}\norm{\eta}_{H^1}^2+\frac{C_0^4}{2\delta}\left(\sum_{j=1}^n\norm{u_j}_{H^1}^2\right)^2\\
\le&~{}\frac{\delta C_0^2}{2}\norm{(u,\eta)}_{\vH^1}^2+\frac{C_0^4}{2\delta}\norm{(u,\eta)}_{\vH^1}^4,
\end{aligned}\]
for all $\delta>0$. Taking $\delta=(\alpha C_0^2)^{-1}>0$, one obtains the desired result.
\end{proof}

Let us define parameters for given $a$ and $c$ by
\begin{equation}\label{eq:beta}
\beta=\max(|a|^{-1},|c|^{-1},1)\quad\mbox{and}\quad\widetilde{\beta}=\beta\max(|a|,|c|,1).
\end{equation}

\begin{proof}[Proof of Theorem \ref{thm:GWP}]
Let $\epsilon_0>0$ be given as in Remark \ref{rem:SLWP}. We choose the initial data $(u_0,\eta_0)\in\vH^1(\R^n)$ satisfying
\begin{equation}\label{eq:smallness}
\norm{(u_0,\eta_0)}_{\vH^1(\R^n)}^2+C_2\norm{(u_0,\eta_0)}_{\vH^1(\R^n)}^4<\epsilon^2,
\end{equation}
where
\begin{equation}\label{eq:C}
C_1=2\widetilde{\beta}+1>1, \quad C_2=\beta^2C_0^6\quad\mbox{and}\quad\epsilon^2=\min\left(\frac{1}{4C_1C_2},\frac{\epsilon_0^2}{2C_1}\right).
\end{equation}
The smallness condition on the initial data \eqref{eq:smallness} yields
\[\norm{(u_0,\eta_0)}_{\vH^1(\R^n)}^2<\epsilon^2 \le \frac{\epsilon_0^2}{2C_1} < \epsilon_0^2,\]
which ensures the local well-posedness in $\vH^1(\R^n)$ on $[0,1]$, thanks to Remark \ref{rem:SLWP}. By $\widetilde{\beta}$ as in \eqref{eq:beta} and Lemma \ref{lem:L3} applied to the initial data with $\alpha=\beta$, we have
\begin{equation}\label{eq:initial energy}
\begin{aligned}
2\beta E(u_0,\eta_0) \le&~{}\widetilde{\beta}\norm{(u_0,\eta_0)}_{\vH^1}^2+\beta\left|\int_{\R^n}\eta_0|u_0|^2\;dx\right|\\
\le&~{}\left(\widetilde{\beta}+\frac12\right)\norm{(u_0,\eta_0)}_{\vH^1}^2+\frac{\beta^2C_0^6}{2}\norm{(u_0,\eta_0)}_{\vH^1}^4.
\end{aligned}
\end{equation}
On the other hand, using the conservation of the energy, Lemma \ref{lem:L3} with $\alpha=\beta$, and \eqref{eq:initial energy}, we obtain
\[\begin{aligned}
\norm{(u,\eta)}_{\vH^1}^2 \le&~{}2\beta E(u,\eta)-\beta\int_{\R^n}\eta|u|^2\;dx\\
\le&~{}2\beta E(u_0,\eta_0)+\frac12\norm{(u,\eta)}_{\vH^1}^2+\frac{\beta^2C_0^6}{2}\norm{(u,\eta)}_{\vH^1}^4\\
\le&~{}\left(\widetilde{\beta}+\frac12\right)\norm{(u_0,\eta_0)}_{\vH^1}^2+\frac{\beta^2C_0^6}{2}\norm{(u_0,\eta_0)}_{\vH^1}^4\\
&~{}+\frac12\norm{(u,\eta)}_{\vH^1}^2+\frac{\beta^2C_0^6}{2}\norm{(u,\eta)}_{\vH^1}^4.
\end{aligned}\]
Thus, the following is valid for every solution $(u,\eta)\in C_t\vH^1$:
\begin{equation}\label{eq:gwp1}
\norm{(u,\eta)}_{\vH^1(\R^n)}^2\le C_1\left(\norm{(u_0,\eta_0)}_{\vH^1}^2+C_2\norm{(u_0,\eta_0)}_{\vH^1}^4\right)+C_2\norm{(u,\eta)}_{\vH^1}^4,
\end{equation}
where $C_1$ and $C_2$ are as in \eqref{eq:C}.

\medskip

On the other hand, let $f(x)=x-C_2x^2$, $x\ge0$. One immediately has $f'(x)=1-2C_2x$, and thus $f$ is increasing on $I=[0,(2C_2)^{-1}]$ from $0$ to $(4C_2)^{-1}$. Additionally, one has $\frac12x \le f(x)$ for all $x\in I$.

By \eqref{eq:gwp1}, \eqref{eq:smallness}, and the definition of $\epsilon$, one has
\[f\left(\norm{(u,\eta)(t)}_{\vH^1(\R^n)}^2\right)<C_1\epsilon^2\le\frac{1}{4C_2},\qquad t\in[0,1].\]
Moreover,
\[\norm{(u_0,\eta_0)}_{\vH^1(\R^n)}^2<\epsilon^2\le \frac{1}{2C_2}.\]
By the continuity of the map $t\mapsto(u,\eta)(t)$, it follows that
\[\norm{(u,\eta)(t)}_{\vH^1(\R^n)}^2\in I,\qquad t\in[0,1].\] 
Therefore,
\[\norm{(u,\eta)(t)}_{\vH^1(\R^n)}^2\le2f\left(\norm{(u,\eta)(t)}_{\vH^1(\R^n)}^2\right)<2C_1\epsilon^2\le\epsilon_0^2,\qquad t\in[0,1].\]
Since $\norm{(u,\eta)(1)}_{\vH^1(\R^n)}<\epsilon_0$, Remark \ref{rem:SLWP} provides a local solution on $[1,2]$, and hence on $[0,2]$. Since \eqref{eq:gwp1} remains valid on $[0,2]$, the same argument applies. Iterating this procedure yields a solution on $[0,T]$ for every $T>0$, which completes the proof of Theorem \ref{thm:GWP}.
\end{proof}

\begin{remark}
The continuity argument used to keep the solution in the coercive regime of the conserved energy is standard in nonlinear dispersive equations, see, for instance, Herr, Tataru and Tzvetkov \cite{HTT}.
\end{remark}

\begin{remark}
The system is time-reversible under
\[(u,\eta)(t,x)\mapsto(-u,\eta)(-t,x).\]
Thus the same continuation argument applies backward in time, and the global solution in Theorem \ref{thm:GWP} is defined for all $t\in\R$.
\end{remark}

\begin{remark}[Persistence of regularity]\label{rem:persistence regularity}
Let $2 \le n \le 4$ and $s \ge 1$, and let $\nu=(\nu_1,\nu_2)$ denote the solution to \eqref{eq:nd abcd} corresponding to $(u,\eta)$. By Remark \ref{rem:persistence}, the nonlinear term satisfies
\[\|\mathcal N(\nu)\|_{H^s(\R^n)\times H^s(\R^n)} \lesssim \|\nu\|_{H^1(\R^n)\times H^1(\R^n)}\|\nu\|_{H^s(\R^n)\times H^s(\R^n)}.\]
Therefore, the Duhamel formula and Gronwall's inequality yield
\[\|\nu(t)\|_{H^s(\R^n)\times H^s(\R^n)}\le\|\nu_0\|_{H^s(\R^n)\times H^s(\R^n)}\exp\left(C\int_0^{|t|}\|\nu(\tau)\|_{H^1(\R^n)\times H^1(\R^n)}\,d\tau\right).\]
In particular, as long as the $H^1(\R^n)\times H^1(\R^n)$ norm remains bounded, no higher $H^s(\R^n)\times H^s(\R^n)$ norm can blow up. Hence, if $(u_0,\eta_0)\in\vH^s(\R^n)$ and its $\vH^1(\R^n)$ norm is sufficiently small as in Theorem \ref{thm:GWP}, then the corresponding global solution satisfies
\[(u,\eta)\in C(\R;\vH^s(\R^n)).\]
\end{remark}

\section{Virial operators in radial symmetry}\label{sec:4}

In this section we rewrite the \emph{abcd} dynamics in radial variables. 

\subsection{Basic definitions}
As usual, we write
\[x=(x_1,x_2)=(r\cos\theta,r\sin\theta),\qquad r=|x|.\]
Let $e_r$ and $e_{\theta}$ be the orthonormal polar basis given by
\[e_r=(\cos\theta,\sin\theta),\qquad e_\theta=(-\sin\theta,\cos\theta).\]
Note by the chain rule that
\begin{equation}\label{eq:cartesian derivatives}
\partial_{x_1}=\cos\theta\partial_r-\frac{\sin\theta}{r}\partial_\theta,\qquad \partial_{x_2}=\sin\theta\partial_r+\frac{\cos\theta}{r}\partial_\theta.
\end{equation}

\medskip

In what follows, we consider radially symmetric solutions of the form
\[u(t,x)=q(t,r)e_r,\qquad \eta(t,x)=h(t,r).\]
Hence
\[\nabla h=(\partial_{x_1}h,\partial_{x_2}h) = h_re_r 
\]
and
\[\Delta h=\partial_{x_1}^2 h+\partial_{x_2}^2 h=h_{rr}+\frac1r h_r.\]
For the right-hand side, we define the scalar radial Laplacian
\begin{equation}\label{eq:scalar radial Laplacian}
L_sh:=h_{rr}+\frac1r h_r.
\end{equation}

We now consider the radial vector field
\[u=qe_r=(q\cos\theta,q\sin\theta)=:(u_1,u_2).\]
Using \eqref{eq:cartesian derivatives}, we have
\begin{equation}\label{eq:cartesian derivatives u}
\begin{aligned}
\partial_{x_1}u_1&=\cos^2\theta q_r+\frac{\sin^2\theta}{r}q, \qquad &\partial_{x_2}u_1&=\sin\theta\cos\theta\left(q_r-\frac qr\right), \\
\partial_{x_1}u_2&=\sin\theta\cos\theta\left(q_r-\frac qr\right), \qquad &\partial_{x_2}u_2&=\sin^2\theta q_r+\frac{\cos^2\theta}{r}q.
\end{aligned}
\end{equation}
Hence
\[\nabla\cdot u=\partial_{x_1}u_1+\partial_{x_2}u_2=q_r+\frac qr.\]
For the right-hand side, we define the radial divergence operator
\[
D_rq:=q_r+\frac qr.
\]
Note that $q_r$ represents the radial derivative, while $q/r$ is the tangential contribution coming from the angular direction.

Note that $u$ is irrotational, so its Laplacian can be written as
\[\Delta u=\nabla(\nabla\cdot u)-\nabla^\perp(\operatorname{curl}u)=\nabla(\nabla\cdot u).\]
Since $D_rq$ is radial,
\[\nabla(D_rq)=\partial_r(D_rq)e_r.\]
A direct computation gives
\[\partial_r(D_rq)=q_{rr}+\frac1r q_r-\frac1{r^2}q.\]
Therefore
\[\Delta(qe_r)=\left(q_{rr}+\frac1r q_r-\frac1{r^2}q\right)e_r.\]
For the vector profile in the right-hand side, we define the radial vector profile Laplacian
\begin{equation}\label{eq:vector radial Laplacian}
L_vq:=q_{rr}+\frac1r q_r-\frac1{r^2}q.
\end{equation}
Using the definitions above, we have
\[D_r(h_r)=h_{rr}+\frac1r h_r=L_sh,\qquad \partial_r(D_rq)=q_{rr}+\frac1r q_r-\frac1{r^2}q=L_vq.\]
Hence
\begin{equation}\label{eq:radial factorization}
L_s=D_r\partial_r,\qquad L_v=\partial_rD_r.
\end{equation}
Moreover,
\begin{equation}\label{eq:radial intertwining}
D_rL_v=D_r\partial_rD_r=L_sD_r,\qquad \partial_rL_s=\partial_rD_r\partial_r=L_v\partial_r.
\end{equation}

\subsection{Polar expression of the Hamiltonian}

We next compute the Jacobian matrix of $u$. By definition,
\[Du=
\begin{pmatrix}
\partial_{x_1}u_1&\partial_{x_2}u_1\\
\partial_{x_1}u_2&\partial_{x_2}u_2
\end{pmatrix}.\]
Using \eqref{eq:cartesian derivatives u}, we obtain
\[Du=
\begin{pmatrix}
\cos^2\theta q_r+\frac{\sin^2\theta}{r}q&
\sin\theta\cos\theta\left(q_r-\frac qr\right)\\
\sin\theta\cos\theta\left(q_r-\frac qr\right)&
\sin^2\theta q_r+\frac{\cos^2\theta}{r}q
\end{pmatrix}.\]
To simplify the computations now, we use tensor product notation. For vectors $a,b\in\R^2$, we define
\[(a\otimes b)_{ij}:=a_ib_j.\]
Since
\[e_r\otimes e_r=
\begin{pmatrix}
\cos^2\theta&\sin\theta\cos\theta\\
\sin\theta\cos\theta&\sin^2\theta
\end{pmatrix}
\quad \text{and} \quad
e_\theta\otimes e_\theta=
\begin{pmatrix}
\sin^2\theta&-\sin\theta\cos\theta\\
-\sin\theta\cos\theta&\cos^2\theta
\end{pmatrix},\]
we write
\[Du=q_re_r\otimes e_r+\frac qr e_\theta\otimes e_\theta.\]
For real $2 \times 2$ matrices $A=(A_{ij})$ and $B=(B_{ij})$, we define the Frobenius inner product by
\[\langle A,B \rangle_F :=\sum_{i,j=1}^2A_{ij}B_{ij},\]
and the corresponding Frobenius norm by
\[|A|^2:=\langle A,A \rangle_F.\]
For $a,b,c,d\in\R^2$, we have
\[\langle a\otimes b , c\otimes d \rangle_F =\sum_{i,j=1}^2a_ib_jc_id_j=\left(\sum_{i=1}^2a_ic_i\right)\left(\sum_{j=1}^2b_jd_j\right)=(a\cdot c)(b\cdot d).\]
By the orthonormality of $e_r$ and $e_\theta$,
\[\langle e_r\otimes e_r , e_r\otimes e_r \rangle_F = \langle e_\theta\otimes e_\theta , e_\theta\otimes e_\theta \rangle_F =1, \quad  \langle e_r\otimes e_r , e_\theta\otimes e_\theta \rangle_F=0.\]
Therefore
\[|Du|^2= \langle q_re_r\otimes e_r+\frac qr e_\theta\otimes e_\theta , q_re_r\otimes e_r+\frac qr e_\theta\otimes e_\theta \rangle_F = q_r^2+\frac{q^2}{r^2}.\]
Finally, note the simple identity
\[|\nabla h|^2=h_r^2.\]


\begin{remark}
For smooth functions $f$ and $g$ with vanishing boundary terms,
\[\int_0^\infty f_rg\;rdr=\int_0^\infty f_r(rg)\;dr=-\int_0^\infty f(rg)_r\;dr=-\int_0^\infty f\left(g_r+\frac gr\right)\;rdr.\]
Hence
\begin{equation}\label{eq:radial integration by parts}
\int_0^\infty f_rg\;rdr=-\int_0^\infty fD_rg\;rdr.
\end{equation}
\end{remark}

\begin{remark}\label{rem:behavior at origin}
Smooth radial scalar and vector fields behave at the origin as follows:

\begin{enumerate}
\item When $h(x)=h(|x|)$ is smooth near $x=0$, its radial profile extends as an even smooth function of $r$, i.e.,  locally,
\[h(r)=h(0)+O(r^2),\qquad h_r(r)=O(r)\quad \text{as} \;\; r\to0,\]
so that
\[h_r(0)=0.\]
Note that no vanishing condition on $h(0)$ is required.

\item For a smooth radial vector field $u(x)=q(r)e_r$, define
\[G(x):=\int_0^{|x|}q(s)\,ds.\]
Then, for $x\ne0$,
\[\nabla G(x)=q(r)e_r=u(x).\]
Thus $G$ is a radial potential for $u$. Since $u$ is smooth and irrotational, $G$ extends smoothly across the origin. By Item (1),
\[G(r)=G(0)+O(r^2),\qquad G_r(r)=O(r)\quad\text{as }\;r\to0.\]
Since $G_r(r)=q(r)$ and $G$ is smooth, we obtain
\[q(r)=O(r),\qquad q_r(r)=O(1)\quad\text{as }\;r\to0,\]
so that
\[q(0)=0.\]
\end{enumerate}
\end{remark}

\subsection{Reduction to the radial system}
We first normalize the parameter $b=d>0$. For a solution $(u,\eta)$ of the original Hamiltonian system \eqref{2D ABCD}, set
\[\widetilde u(t,x):=u(\sqrt b\,t,\sqrt b\,x),\qquad\widetilde\eta(t,x):=\eta(\sqrt b\,t,\sqrt b\,x).\]
The rescaled system has $b=d=1$ and coefficients
\[\widetilde a:=\frac ab,\qquad\widetilde c:=\frac cb.\]

We drop the tildes and consider the normalized two dimensional Hamiltonian \emph{abcd} system
\begin{equation}\label{eq:normalized abcd}
\begin{aligned}
(1-\Delta)\eta_t+\nabla\cdot(a\Delta u+u+\eta u)=&~{}0,\\
(1-\Delta)u_t+\nabla\left(c\Delta\eta+\eta+\frac12|u|^2\right)=&~{}0,
\end{aligned}
\end{equation}
where $a,c<0$.

We consider radially symmetric solutions
\[u(t,x)=q(t,r)e_r,\qquad\eta(t,x)=h(t,r),\qquad r=|x|.\]
Using \eqref{eq:scalar radial Laplacian} and \eqref{eq:vector radial Laplacian}, we write
\[a\Delta u+u+\eta u=(aL_vq+q+qh)e_r\]
and
\[c\Delta\eta+\eta+\frac12|u|^2=cL_sh+h+\frac12q^2.\]
Set
\begin{equation}\label{eq:radial fluxes}
\mathcal A:=aL_vq+q+qh,\qquad \mathcal B:=cL_sh+h+\frac12q^2.
\end{equation}
A direct computation reduces the system \eqref{eq:normalized abcd} to
\begin{equation}\label{eq:radial system}
(1-L_s)h_t+D_r\mathcal A=0,\qquad (1-L_v)q_t+\partial_r\mathcal B=0.
\end{equation}
From the previous computations, we introduce the radial resolvent operators for canonical variables:
\be\label{KsKv}
K_s:=(1-L_s)^{-1},\qquad K_v:=(1-L_v)^{-1}.
\ee
Since $(1-\Delta)^{-1}$ preserves radial irrotational vector fields, $K_s$ and $K_v$ are well-defined, because $K_sh$ is the radial profile of $(1-\Delta)^{-1}\eta$ for $\eta=h(|x|)$, and $K_vq$ the profile of $(1-\Delta)^{-1}u$ for $u=q(|x|)e_r$. This identification gives existence and uniqueness in $L^2(rdr)$ for \eqref{KsKv}.

Using \eqref{eq:radial intertwining}, we obtain
\begin{equation}\label{eq:resolvent intertwining}
K_sD_r=D_rK_v,\qquad K_v\partial_r=\partial_rK_s.
\end{equation}
Moreover, the resolvents commute with the corresponding radial Laplacians, i.e.,
\begin{equation}\label{eq:resolvent commutation}
L_sK_s=K_sL_s=K_s-1,\qquad L_vK_v=K_vL_v=K_v-1.
\end{equation}

\begin{lemma}
The radial system \eqref{eq:radial system} can be written as
\begin{equation}\label{eq:radial evolution}
\begin{aligned}
h_t=&~{}aD_rq-(a+1)K_sD_rq-K_sD_r(qh),\\
q_t=&~{}ch_r-(c+1)K_vh_r-\frac12K_v\partial_r(q^2).
\end{aligned}
\end{equation}
Consequently,
\begin{equation}\label{eq:radial derivative evolution}
\begin{aligned}
h_{tr}=&~{}aL_vq-(a+1)\partial_rK_sD_rq-\partial_rK_sD_r(qh),\\
q_{tr}=&~{}ch_{rr}-(c+1)\partial_rK_vh_r-\frac12\partial_rK_v\partial_r(q^2).
\end{aligned}
\end{equation}
\end{lemma}

\begin{proof}
Using \eqref{eq:radial fluxes} in the first equation of \eqref{eq:radial system}, we have
\[(1-L_s)h_t=-aD_rL_vq-D_rq-D_r(qh).\]
By \eqref{eq:radial intertwining},
\[(1-L_s)h_t=-aL_sD_rq-D_rq-D_r(qh).\]
Applying $K_s$ and using \eqref{eq:resolvent commutation}, we obtain
\[\begin{aligned}
h_t=&~{}-aK_sL_sD_rq-K_sD_rq-K_sD_r(qh)\\
=&~{}-a(K_s-1)D_rq-K_sD_rq-K_sD_r(qh)\\
=&~{}aD_rq-(a+1)K_sD_rq-K_sD_r(qh).
\end{aligned}\]

Similarly, the second equation of \eqref{eq:radial system} gives
\[(1-L_v)q_t=-c\partial_rL_sh-h_r-\frac12\partial_r(q^2).\]
By \eqref{eq:radial intertwining},
\[(1-L_v)q_t=-cL_vh_r-h_r-\frac12\partial_r(q^2).\]
Applying $K_v$ and using \eqref{eq:resolvent commutation}, we obtain
\[\begin{aligned}
q_t=&~{}-cK_vL_vh_r-K_vh_r-\frac12K_v\partial_r(q^2)\\
=&~{}-c(K_v-1)h_r-K_vh_r-\frac12K_v\partial_r(q^2)\\
=&~{}ch_r-(c+1)K_vh_r-\frac12K_v\partial_r(q^2).
\end{aligned}\]
This proves \eqref{eq:radial evolution}. Differentiating in $r$ gives \eqref{eq:radial derivative evolution}.
\end{proof}

\subsection{Radial virial identities}

For a smooth bounded radial weight $\phi=\phi(r)$ independent of time, define
\begin{equation}\label{eq:principal virial functional}
I_\phi(t):=2\pi\int_0^\infty\phi(r)(qh+q_rh_r)(t,r)\;rdr.
\end{equation}
Recall that
\[\nabla\eta=h_re_r,\qquad Du=q_re_r\otimes e_r+\frac qr e_\theta\otimes e_\theta.\]
Using $(a\otimes b)v=a(b\cdot v)$ and orthonormality of $e_r$ and $e_\theta$, we obtain
\[(Du)\nabla\eta = h_r\left(q_re_r\otimes e_r+\frac qr e_\theta\otimes e_\theta\right)e_r = q_rh_re_r.\]
Hence
\[\eta u+(Du)\nabla\eta=(qh+q_rh_r)e_r.\]
Therefore \eqref{eq:principal virial functional} can equivalently be written as
\[I_\phi(t)=\int_{\R^2}\phi(|x|)e_r\cdot\left(\eta u+(Du)\nabla\eta\right)(t,x)dx,\]
which is the corresponding localized momentum in spatial space. Moreover, for $(u,\eta)\in H^1(\R^2)\times H^1(\R^2)$,
\begin{equation}\label{eq:I_phi bound}
\begin{aligned}
|I_\phi(t)|\le&~{}\|\phi\|_{L^\infty}\left(\|\eta\|_{L^2}\|u\|_{L^2}+\|Du\|_{L^2}\|\nabla\eta\|_{L^2}\right)\\
\lesssim&~{}\|\phi\|_{L^\infty}\left(\|u\|_{H^1}^2+\|\eta\|_{H^1}^2\right).
\end{aligned}
\end{equation}
Thus $I_\phi(t)$ is well defined and bounded for energy solutions.

\begin{lemma}[Principal virial identity]\label{lem:principal virial identity}
Let $(q,h)$ be the profiles of a smooth radial solution. Suppose that $\phi(r)=O(r)$ as $r\to0$ and that $\phi$ and the solution have sufficient decay at infinity so that all integrals converge and all boundary terms vanish. Then
\begin{equation}\label{eq:principal virial decomposition}
\frac1{2\pi}\frac{d}{dt}I_\phi(t)=I_\phi^{(2)}(t)+I_\phi^{(3)}(t),
\end{equation}
where the quadratic contribution is
\begin{equation}\label{eq:principal virial quadratic}
\begin{aligned}
I_\phi^{(2)}(t):=&~{}\int_0^\infty\left(\phi_r-\frac{\phi}{r}\right)\left(\frac12q^2-\frac a2q_r^2+\frac a2\frac{q^2}{r^2}\right)\;rdr \\
&~{}+\int_0^\infty\left(\phi_r+\frac{\phi}{r}\right)\left(\frac12h^2-\frac c2h_r^2\right)\;rdr\\
&~{}+(a+1)\int_0^\infty\phi_rq(\partial_rK_sD_rq)\;rdr\\
&~{}+(c+1)\int_0^\infty\phi_rh(\partial_rK_vh_r) + \frac{\phi}{r^2}h(K_vh_r)\;rdr\\
\end{aligned}
\end{equation}
and the cubic contribution is
\begin{equation}\label{eq:principal virial cubic}
\begin{aligned}
I_\phi^{(3)}(t):=&~{}\int_0^\infty\phi_rhq^2\;rdr +\int_0^\infty\phi_rq(\partial_rK_sD_r(qh))\;rdr\\
&~{}+\frac12\int_0^\infty\phi_rh(\partial_rK_v\partial_r(q^2)) + \frac{\phi}{r^2}h(K_v\partial_r(q^2))\;rdr\\
\end{aligned}
\end{equation}
\end{lemma}

\begin{proof}
Differentiating \eqref{eq:principal virial functional} in time gives
\[\frac1{2\pi}\frac{d}{dt}I_\phi(t)=\int_0^\infty\phi(qh_t+hq_t+q_rh_{tr}+h_rq_{tr})\;rdr.\]
We first consider the terms containing $h_t$. Using \eqref{eq:scalar radial Laplacian}, we have
\[(1-L_s)h_t=h_t-h_{trr}-\frac1r h_{tr} = h_t - \frac1r(rh_{tr})_r.\]
Therefore
\[\begin{aligned}
\int_0^\infty\phi q(1-L_s)h_t\;rdr =&~{}\int_0^\infty\phi qh_t\;rdr-\int_0^\infty\phi q(rh_{tr})_r\;dr\\
=&~{}\int_0^\infty\phi(qh_t+q_rh_{tr})\;rdr+\int_0^\infty\phi_rqh_{tr}\;rdr.
\end{aligned}\]
Hence
\begin{equation}\label{eq:virial ht terms}
\int_0^\infty\phi(qh_t+q_rh_{tr})\;rdr =\int_0^\infty\phi q(1-L_s)h_t\;rdr-\int_0^\infty\phi_rqh_{tr}\;rdr.
\end{equation}
Similarly, using \eqref{eq:vector radial Laplacian}, we have
\[(1-L_v)q_t=q_t-q_{trr}-\frac1r q_{tr}+\frac1{r^2}q_t = q_t- \frac1r(rq_{tr})_r+\frac1{r^2}q_t.\]
Therefore
\[\begin{aligned}
\int_0^\infty\phi h(1-L_v)q_t\;rdr =&~{}\int_0^\infty\phi hq_t\;rdr-\int_0^\infty\phi h(rq_{tr})_r\;dr+\int_0^\infty\frac{\phi}{r}hq_t\;dr\\
=&~{}\int_0^\infty\phi(hq_t+h_rq_{tr})\;rdr+\int_0^\infty\phi_rhq_{tr}\;rdr+\int_0^\infty\frac{\phi}{r}hq_t\;dr.
\end{aligned}\]
Hence
\begin{equation}\label{eq:virial qt terms}
\int_0^\infty\phi(hq_t+h_rq_{tr})\;rdr =\int_0^\infty\phi h(1-L_v)q_t\;rdr -\int_0^\infty\phi_rhq_{tr}\;rdr-\int_0^\infty\frac{\phi}{r}hq_t\;dr.
\end{equation}

Combining \eqref{eq:virial ht terms} and \eqref{eq:virial qt terms}, and using \eqref{eq:radial system}, we obtain
\[\begin{aligned}
\frac1{2\pi}\frac{d}{dt}I_\phi(t) =&~{}-\int_0^\infty\phi qD_r\mathcal A\;rdr-\int_0^\infty\phi h\partial_r\mathcal B\;rdr\\
&~{}-\int_0^\infty\phi_rqh_{tr}\;rdr-\int_0^\infty\phi_rhq_{tr}\;rdr-\int_0^\infty\frac{\phi}{r}hq_t\;dr.
\end{aligned}\]
By \eqref{eq:radial integration by parts},
\[-\int_0^\infty\phi qD_r\mathcal A\;rdr =\int_0^\infty\phi q_r\mathcal A\;rdr+\int_0^\infty\phi_rq\mathcal A\;rdr.\]
Moreover, integration by parts yields
\[-\int_0^\infty\phi h\partial_r\mathcal B\;rdr =\int_0^\infty\phi h_r\mathcal B\;rdr+\int_0^\infty\phi_rh\mathcal B\;rdr+\int_0^\infty\phi h\mathcal B\;dr.\]
Therefore
\begin{equation}\label{eq:radial virial first expansion}
\begin{aligned}
\frac1{2\pi}\frac{d}{dt}I_\phi(t)=&~{}\int_0^\infty\phi q_r\mathcal A\;rdr+\int_0^\infty\phi_rq\mathcal A\;rdr\\
&~{}+\int_0^\infty\phi h_r\mathcal B\;rdr+\int_0^\infty\phi_rh\mathcal B\;rdr+\int_0^\infty\phi h\mathcal B\;dr\\
&~{}-\int_0^\infty\phi_rqh_{tr}\;rdr-\int_0^\infty\phi_rhq_{tr}\;rdr-\int_0^\infty\frac{\phi}{r}hq_t\;dr.
\end{aligned}
\end{equation}

Substituting \eqref{eq:radial fluxes}, \eqref{eq:radial evolution}, and \eqref{eq:radial derivative evolution} into \eqref{eq:radial virial first expansion}, we obtain
\begin{equation}\label{eq:virial A first}
\int_0^\infty\phi q_r\mathcal A\;rdr =a\int_0^\infty\phi q_rL_vq\;rdr+\int_0^\infty\phi qq_r\;rdr+\int_0^\infty\phi hqq_r\;rdr,
\end{equation}
\begin{equation}\label{eq:virial A second}
\int_0^\infty\phi_rq\mathcal A\;rdr =a\int_0^\infty\phi_rqL_vq\;rdr+\int_0^\infty\phi_rq^2\;rdr+\int_0^\infty\phi_rhq^2\;rdr,
\end{equation}
\begin{equation}\label{eq:virial B first}
\int_0^\infty\phi h_r\mathcal B\;rdr =c\int_0^\infty\phi h_rL_sh\;rdr+\int_0^\infty\phi hh_r\;rdr+\frac12\int_0^\infty\phi q^2h_r\;rdr,
\end{equation}
\begin{equation}\label{eq:virial B second}
\int_0^\infty\phi_rh\mathcal B\;rdr =c\int_0^\infty\phi_rhL_sh\;rdr+\int_0^\infty\phi_rh^2\;rdr+\frac12\int_0^\infty\phi_rhq^2\;rdr,
\end{equation}
\begin{equation}\label{eq:virial B third}
\int_0^\infty\phi h\mathcal B\;dr =c\int_0^\infty\phi hL_sh\;dr+\int_0^\infty\phi h^2\;dr+\frac12\int_0^\infty\phi hq^2\;dr.
\end{equation}
Moreover,
\begin{equation}\label{eq:virial htr}
\begin{aligned}
-\int_0^\infty\phi_rqh_{tr}\;rdr =&~{}-a\int_0^\infty\phi_rqL_vq\;rdr +(a+1)\int_0^\infty\phi_rq(\partial_rK_sD_rq)\;rdr \\
&~{}+\int_0^\infty\phi_rq(\partial_rK_sD_r(qh))\;rdr,
\end{aligned}
\end{equation}
\begin{equation}\label{eq:virial qtr}
\begin{aligned}
-\int_0^\infty\phi_rhq_{tr}\;rdr =&~{}-c\int_0^\infty\phi_rhh_{rr}\;rdr+(c+1)\int_0^\infty\phi_rh(\partial_rK_vh_r)\;rdr\\
&~{}+\frac12\int_0^\infty\phi_rh(\partial_rK_v\partial_r(q^2))\;rdr,
\end{aligned}
\end{equation}
and
\begin{equation}\label{eq:virial qt}
\begin{aligned}
-\int_0^\infty\frac{\phi}{r}hq_t\;dr =&~{}-c\int_0^\infty\frac{\phi}{r}hh_r\;dr +(c+1)\int_0^\infty\frac{\phi}{r}h(K_vh_r)\;dr\\
&~{}+\frac12\int_0^\infty\frac{\phi}{r}h(K_v\partial_r(q^2))\;dr.
\end{aligned}
\end{equation}

We first collect the local quadratic terms involving $q$. Combining \eqref{eq:virial A second} and \eqref{eq:virial htr}, the terms involving $a\phi_rqL_vq$ cancel exactly. Using \eqref{eq:vector radial Laplacian} and integrating by parts,
\[\begin{aligned}
\int_0^\infty\phi q_rL_vq\;rdr =&~{}\int_0^\infty\phi q_r \left(q_{rr}+\frac1rq_r-\frac1{r^2}q\right)\;rdr\\ 
=&~{}-\frac12\int_0^\infty(r\phi_r+\phi)q_r^2\;dr+\int_0^\infty\phi q_r^2\;dr+\frac12\int_0^\infty\left(\frac{\phi_r}{r}-\frac{\phi}{r^2}\right)q^2\;dr\\
=&~{}\frac12\int_0^\infty\left(\phi_r-\frac{\phi}{r}\right)\left(-q_r^2+\frac{q^2}{r^2}\right)\;rdr.
\end{aligned}\]
Moreover,
\[\int_0^\infty\phi qq_r\;rdr+\int_0^\infty\phi_rq^2\;rdr=\frac12\int_0^\infty\left(\phi_r-\frac{\phi}{r}\right)q^2\;rdr.\]
Hence the local quadratic $q$ terms in \eqref{eq:virial A first}, \eqref{eq:virial A second}, and \eqref{eq:virial htr} reduce to
\begin{equation}\label{eq:local q quadratic}
\int_0^\infty\left(\phi_r-\frac{\phi}{r}\right)\left(\frac12q^2-\frac a2q_r^2+\frac a2\frac{q^2}{r^2}\right)\;rdr.
\end{equation}

We next collect the local quadratic terms involving $h$. Combining \eqref{eq:virial B second} and \eqref{eq:virial qtr}, and using \eqref{eq:scalar radial Laplacian}, we obtain
\[c\int_0^\infty\phi_rhL_sh\;rdr-c\int_0^\infty\phi_rhh_{rr}\;rdr=c\int_0^\infty\phi_rhh_r\;dr.\]
Similarly, combining \eqref{eq:virial B third} and \eqref{eq:virial qt}, we obtain
\[c\int_0^\infty\phi hL_sh\;dr-c\int_0^\infty\frac{\phi}{r}hh_r\;dr=c\int_0^\infty\phi hh_{rr}\;dr.\]
Adding the two right hand sides, an integration by parts gives
\begin{equation}\label{eq:h_r^2}
c\int_0^\infty\phi_rhh_r\;dr+c\int_0^\infty\phi hh_{rr}\;dr =-c\int_0^\infty\phi h_r^2\;dr.
\end{equation}
Using \eqref{eq:scalar radial Laplacian} and integrating by parts,
\[\int_0^\infty\phi h_rL_sh\;rdr=\int_0^\infty\phi h_r\left(h_{rr}+\frac1r h_r\right)\;rdr =\frac12\int_0^\infty(\phi-r\phi_r)h_r^2\;dr.\]
Together with \eqref{eq:h_r^2}, we have
\[c\int_0^\infty\phi h_rL_sh\;rdr+c\int_0^\infty\phi_rhh_r\;dr+c\int_0^\infty\phi hh_{rr}\;dr=-\frac c2\int_0^\infty\left(\phi_r+\frac{\phi}{r}\right)h_r^2\;rdr.\]
Similarly,
\[\int_0^\infty\phi hh_r\;rdr+\int_0^\infty\phi_rh^2\;rdr+\int_0^\infty\phi h^2\;dr=\frac12\int_0^\infty\left(\phi_r+\frac{\phi}{r}\right)h^2\;rdr.\]
Hence the local quadratic $h$ terms in \eqref{eq:virial B first}--\eqref{eq:virial B third}, \eqref{eq:virial qtr}, and \eqref{eq:virial qt} reduce to
\begin{equation}\label{eq:local h quadratic}
\int_0^\infty\left(\phi_r+\frac{\phi}{r}\right)\left(\frac12h^2-\frac c2h_r^2\right)\;rdr.
\end{equation}

\medskip

Finally, we collect the local cubic terms. For cubic terms in \eqref{eq:virial A first} and \eqref{eq:virial B first}, integration by parts gives
\[\int_0^\infty\phi hqq_r\;rdr+\frac12\int_0^\infty\phi q^2h_r\;rdr=-\frac12\int_0^\infty(r\phi_r+\phi)hq^2\;dr.\]
Combining this with the remaining local cubic contributions in \eqref{eq:virial A second}, \eqref{eq:virial B second}, and \eqref{eq:virial B third}, we obtain
\begin{equation}\label{eq:local cubic}
\int_0^\infty\phi_rhq^2\;rdr.
\end{equation}

Combining \eqref{eq:local q quadratic} and \eqref{eq:local h quadratic} with the remaining quadratic resolvent contributions in \eqref{eq:virial htr}, \eqref{eq:virial qtr}, and \eqref{eq:virial qt} gives \eqref{eq:principal virial quadratic}. Adding the remaining cubic resolvent terms from \eqref{eq:virial htr}, \eqref{eq:virial qtr}, and \eqref{eq:virial qt} to \eqref{eq:local cubic} yields \eqref{eq:principal virial cubic}. This proves \eqref{eq:principal virial decomposition}.
\end{proof}

\begin{remark}
Assume that $\phi$ is a nontrivial nonnegative bounded weight. The quadratic identity \eqref{eq:principal virial quadratic} contains different combinations of the radial $\phi_r$ and tangential $\phi/r$ weights. Since $a,c<0$, nonnegativity of the coefficients of $q_r^2$ and $h_r^2$ would require
\[\phi_r-\frac{\phi}{r}\ge0,\qquad \phi_r+\frac{\phi}{r}\ge0.\]
The obstruction comes from the first condition, since
\[\left(\frac{\phi}{r}\right)_r=\frac1r\left(\phi_r-\frac{\phi}{r}\right).\]
Thus $\phi_r-\phi/r\ge0$ would imply that $\phi(r)/r$ is nondecreasing.  On the other hand, requiring $\phi$ bounded implies $\frac{\phi(r)}r\to0$  as $r\to\infty$. Therefore, no nontrivial nonnegative bounded weight can satisfy the required condition globally in space. It is then necessary to introduce a modification on the virial functionals in order to compensate the tangential contribution $\phi/r$ in the (radial) derivative terms.
\end{remark}

\medskip

Let us define
\begin{equation}\label{eq:tangential derivative correction}
\mathcal J_\phi(t):=2\pi\int_0^\infty\frac{\phi}{r}hD_rq\;rdr.
\end{equation}
Notice that $\phi(r)=O(r)$ as $r\to0$, and $\phi$ is bounded. From this, we have
\[\sup_{r>0}\left|\frac{\phi(r)}r\right|<\infty.\]
Using the Cauchy--Schwarz inequality and the fact that $D_rq=\nabla\cdot u$,
\begin{equation}\label{eq:J_phi bound}
\begin{aligned}
|\mathcal J_\phi(t)| \le&~{}2\pi\sup_{r>0}\left|\frac{\phi(r)}r\right|\|h\|_{L^2(rdr)}\|D_rq\|_{L^2(rdr)}\\
\lesssim&~{}\sup_{r>0}\left|\frac{\phi(r)}r\right|\|\eta\|_{L^2(\R^2)}\|\nabla u\|_{L^2(\R^2)}\\
\lesssim&~{}\sup_{r>0}\left|\frac{\phi(r)}r\right|\left(\|u\|_{H^1(\R^2)}^2+\|\eta\|_{H^1(\R^2)}^2\right).
\end{aligned}
\end{equation}
Thus $\mathcal J_\phi(t)$ is well defined and bounded for radial energy solutions.

\begin{lemma}[Tangential derivative correction]\label{lem:tangential derivative virial}
Under the assumptions of Lemma \ref{lem:principal virial identity},
\begin{equation}\label{eq:tangential derivative decomposition}
\frac1{2\pi}\frac{d}{dt}\mathcal J_\phi(t)=\mathcal J_\phi^{(2)}(t)+\mathcal J_\phi^{(3)}(t),
\end{equation}
where
\begin{equation}\label{eq:tangential derivative quadratic}
\begin{aligned}
\mathcal J_\phi^{(2)}(t):=&~{}a\int_0^\infty\frac{\phi}{r}q_r^2\;rdr-c\int_0^\infty\frac{\phi}{r}h_r^2\;rdr\\
&~{}+a\int_0^\infty\left(\frac{\phi}{r^3}-\frac1r\left(\frac{\phi}{r}\right)_r\right)q^2\;rdr+\frac c2\int_0^\infty L_s\left(\frac{\phi}{r}\right)h^2\;rdr\\
&~{}-(a+1)\int_0^\infty\frac{\phi}{r}q^2\;rdr+(c+1)\int_0^\infty\frac{\phi}{r}h^2\;rdr\\
&~{}+(a+1)\int_0^\infty\frac{\phi}{r}q(K_vq)\;rdr-(c+1)\int_0^\infty\frac{\phi}{r}h(K_sh)\;rdr\\
&~{}+(a+1)\int_0^\infty\left(\frac{\phi}{r}\right)_rq(D_rK_vq)\;rdr,
\end{aligned}
\end{equation}
and
\begin{equation}\label{eq:tangential derivative cubic}
\begin{aligned}
\mathcal J_\phi^{(3)}(t):=&~{}\int_0^\infty\frac{\phi}{r}qK_v(qh)\;rdr-\frac12\int_0^\infty\frac{\phi}{r}hq^2\;rdr\\
&~{}-\frac12\int_0^\infty\frac{\phi}{r}hK_s(q^2)\;rdr+\int_0^\infty\left(\frac{\phi}{r}\right)_rqD_rK_v(qh)\;rdr.
\end{aligned}
\end{equation}
\end{lemma}

\begin{proof}
Differentiating \eqref{eq:tangential derivative correction} in time gives
\[\frac1{2\pi}\frac{d}{dt}\mathcal J_\phi(t)=\int_0^\infty\frac{\phi}{r}h_tD_rq\;rdr+\int_0^\infty\frac{\phi}{r}hD_rq_t\;rdr.\]
Using \eqref{eq:radial evolution}, together with \eqref{eq:radial factorization}, \eqref{eq:resolvent intertwining}, and \eqref{eq:resolvent commutation}, we obtain
\begin{equation}\label{eq:tangential derivative ht expansion}
\begin{aligned}
\int_0^\infty\frac{\phi}{r}h_tD_rq\;rdr=&~{}\int_0^\infty\frac{\phi}{r}\left(aD_rq-(a+1)K_sD_rq-K_sD_r(qh)\right)D_rq\;rdr\\
=&~{}a\int_0^\infty\frac{\phi}{r}(D_rq)^2\;rdr-(a+1)\int_0^\infty\frac{\phi}{r}(D_rK_vq)(D_rq)\;rdr\\
&~{}-\int_0^\infty\frac{\phi}{r}(D_rK_v(qh))(D_rq)\;rdr,
\end{aligned}
\end{equation}
and
\begin{equation}\label{eq:tangential derivative qt expansion}
\begin{aligned}
\int_0^\infty\frac{\phi}{r}hD_rq_t\;rdr=&~{}\int_0^\infty\frac{\phi}{r}hD_r\left(ch_r-(c+1)K_vh_r-\frac12K_v\partial_r(q^2)\right)\;rdr\\
=&~{}c\int_0^\infty\frac{\phi}{r}hL_sh\;rdr-(c+1)\int_0^\infty\frac{\phi}{r}hL_sK_sh\;rdr\\
&~{}-\frac12\int_0^\infty\frac{\phi}{r}hL_sK_s(q^2)\;rdr.
\end{aligned}
\end{equation}

We first consider \eqref{eq:tangential derivative ht expansion}. Expanding $D_rq=q_r+q/r$, we have
\[\begin{aligned}
\int_0^\infty\frac{\phi}{r}(D_rq)^2\;rdr
=&~{}\int_0^\infty\frac{\phi}{r}q_r^2\;rdr+2\int_0^\infty\frac{\phi}{r}qq_r\;dr+\int_0^\infty\frac{\phi}{r^3}q^2\;rdr\\
=&~{}\int_0^\infty\frac{\phi}{r}q_r^2\;rdr-\int_0^\infty\left(\frac{\phi}{r}\right)_rq^2\;dr+\int_0^\infty\frac{\phi}{r^3}q^2\;rdr.
\end{aligned}\]
Hence
\begin{equation}\label{eq:tangential derivative q identity}
\int_0^\infty\frac{\phi}{r}(D_rq)^2\;rdr=\int_0^\infty\frac{\phi}{r}q_r^2\;rdr+\int_0^\infty\left(\frac{\phi}{r^3}-\frac1r\left(\frac{\phi}{r}\right)_r\right)q^2\;rdr.
\end{equation}
For the second term, using \eqref{eq:radial integration by parts}, we obtain
\[-\int_0^\infty\frac{\phi}{r}(D_rK_vq)(D_rq)\;rdr =\int_0^\infty q\partial_r\left(\frac{\phi}{r}D_rK_vq\right)\;rdr.\]
By the product rule,
\[\partial_r\left(\frac{\phi}{r}D_rK_vq\right) =\left(\frac{\phi}{r}\right)_rD_rK_vq+\frac{\phi}{r}\partial_rD_rK_vq.\]
Using \eqref{eq:radial factorization} and \eqref{eq:resolvent commutation}, we obtain
\[
\begin{aligned}
-\int_0^\infty\frac{\phi}{r}(D_rK_vq)(D_rq)\;rdr =&~{}\int_0^\infty\frac{\phi}{r}qL_vK_vq\;rdr +\int_0^\infty\left(\frac{\phi}{r}\right)_rqD_rK_vq\;rdr\\
=&~{}\int_0^\infty\frac{\phi}{r}q(K_vq-q)\;rdr+\int_0^\infty\left(\frac{\phi}{r}\right)_rqD_rK_vq\;rdr.
\end{aligned}
\]
Replacing $q$ by $qh$, we have
\begin{equation}\label{eq:tangential derivative nonlinear resolvent}
\ba
-\int_0^\infty\frac{\phi}{r}(D_rK_v(qh))(D_rq)\;rdr= &~{} \int_0^\infty\frac{\phi}{r}q\left(K_v(qh)-qh\right)\;rdr \\
&~{} +\int_0^\infty\left(\frac{\phi}{r}\right)_rqD_rK_v(qh)\;rdr.
\ea
\end{equation}

We next consider \eqref{eq:tangential derivative qt expansion}. For the first term, using \eqref{eq:radial factorization} and \eqref{eq:radial integration by parts},
\[\begin{aligned}
\int_0^\infty\frac{\phi}{r}hL_sh\;rdr=&~{}\int_0^\infty\frac{\phi}{r}hD_rh_r\;rdr\\
=&~{}-\int_0^\infty\partial_r\left(\frac{\phi}{r}h\right)h_r\;rdr\\
=&~{}-\int_0^\infty\frac{\phi}{r}h_r^2\;rdr-\int_0^\infty r\left(\frac{\phi}{r}\right)_rhh_r\;dr\\
=&~{}-\int_0^\infty\frac{\phi}{r}h_r^2\;rdr+\frac12\int_0^\infty\left(r\left(\frac{\phi}{r}\right)_r\right)_rh^2\;dr.
\end{aligned}\]
Note by \eqref{eq:scalar radial Laplacian} that
\[rL_s\left(\frac{\phi}{r}\right)=r\left(\frac{\phi}{r}\right)_{rr}+ \left(\frac{\phi}{r}\right)_r = \left(r\left(\frac{\phi}{r}\right)_r\right)_r.\]
Thus
\[
\int_0^\infty\frac{\phi}{r}hL_sh\;rdr=-\int_0^\infty\frac{\phi}{r}h_r^2\;rdr+\frac12\int_0^\infty L_s\left(\frac{\phi}{r}\right)h^2\;rdr.
\]
For the second term, using \eqref{eq:resolvent commutation},
\[
\ba
-\int_0^\infty\frac{\phi}{r}hL_sK_sh\;rdr = &~{} -\int_0^\infty\frac{\phi}{r}h(K_s-1)h\;rdr \\
=&~{} \int_0^\infty\frac{\phi}{r}h^2\;rdr-\int_0^\infty\frac{\phi}{r}hK_sh\;rdr.
\ea
\]
Similarly, we have
\begin{equation}\label{eq:tangential derivative scalar nonlinear resolvent}
-\frac12\int_0^\infty\frac{\phi}{r}hL_sK_s(q^2)\;rdr =\frac12\int_0^\infty\frac{\phi}{r}hq^2\;rdr-\frac12\int_0^\infty\frac{\phi}{r}hK_s(q^2)\;rdr.
\end{equation}

Substituting \eqref{eq:tangential derivative q identity}--\eqref{eq:tangential derivative scalar nonlinear resolvent} into \eqref{eq:tangential derivative ht expansion} and \eqref{eq:tangential derivative qt expansion}, and collecting the quadratic terms, we obtain \eqref{eq:tangential derivative quadratic}. For the cubic terms, the local contributions from \eqref{eq:tangential derivative nonlinear resolvent} and \eqref{eq:tangential derivative scalar nonlinear resolvent} combine as
\[-\int_0^\infty\frac{\phi}{r}hq^2\;rdr+\frac12\int_0^\infty\frac{\phi}{r}hq^2\;rdr=-\frac12\int_0^\infty\frac{\phi}{r}hq^2\;rdr,\]
while the remaining nonlinear resolvent terms give \eqref{eq:tangential derivative cubic}. This proves \eqref{eq:tangential derivative decomposition}.
\end{proof}

\begin{remark}
The purpose of the correction is seen by combining \eqref{eq:principal virial quadratic} with $-\mathcal J_\phi^{(2)}/2$. For the radial derivative terms,
\[-\frac a2\left(\phi_r-\frac{\phi}{r}\right)-\frac a2\frac{\phi}{r}=-\frac a2\phi_r,\qquad -\frac c2\left(\phi_r+\frac{\phi}{r}\right)+\frac c2\frac{\phi}{r}=-\frac c2\phi_r.\]
For the term containing $q^2/r^2$, using $r(\phi/r)_r=\phi_r-\phi/r$, the resulting coefficient is
\[\frac a2\left(\phi_r-\frac{\phi}{r}\right)-\frac a2\left(\frac{\phi}{r}-r\left(\frac{\phi}{r}\right)_r\right)=-\frac a2\left(3\frac{\phi}{r}-2\phi_r\right).\]
Hence the highest order local terms become
\begin{equation}\label{eq:tangential corrected derivative terms}
-\frac a2\int_0^\infty\phi_rq_r^2\;rdr-\frac c2\int_0^\infty\phi_rh_r^2\;rdr-\frac a2\int_0^\infty\left(3\frac{\phi}{r}-2\phi_r\right)\frac{q^2}{r^2}\;rdr.
\end{equation}
They have nonnegative coefficients when $\phi_r\ge0$ and $3\phi/r\ge2\phi_r$. In particular, both conditions hold for our choice of the weight in Section \ref{sec:local weight}.
\end{remark}

\medskip

We now examine the lower order quadratic terms in $I_\phi^{(2)}-\mathcal J_\phi^{(2)}/2$. By \eqref{eq:principal virial quadratic} and \eqref{eq:tangential derivative quadratic}, the quadratic mass terms carrying the weights $\phi_r$ or $\phi/r$ combine as
\begin{equation}\label{eq:tangential correction remaining mass terms}
\frac12\int_0^\infty\phi_r(q^2+h^2)\;rdr
+\frac12\int_0^\infty\frac{\phi}{r}\left(aq^2-(a+1)q(K_vq)-ch^2+(c+1)h(K_sh)\right)\;rdr.
\end{equation}
The second integral still carries the tangential weight $\phi/r$. By \eqref{eq:radial evolution} and \eqref{eq:resolvent intertwining}, $h_t$ can be written as a radial divergence, while $q_t$ is a radial derivative. This suggests pairing $q$ with an inverse radial divergence of $h$ to compensate for this remaining tangential contribution.

We first introduce the inverse radial divergence operator. Since
\[D_rF=F_r+\frac Fr=\frac1r(rF)_r,\]
the regular solution of $D_rF=z$ is naturally given by the weighted Hardy type operator
\begin{equation}\label{eq:regular radial inverse}
(D_r^{-1}z)(r):=\frac1r\int_0^r sz(s) \;ds.
\end{equation}

\begin{lemma}\label{lem:regular radial inverse}
For $z\in L^2((0,\infty),rdr)$, the operator defined in \eqref{eq:regular radial inverse} satisfies
\[D_rD_r^{-1}z=z\]
almost everywhere on $(0,\infty)$.

Moreover, for every smooth radial vector profile $F$,
\[D_r^{-1}D_rF=F.\]

The following Hardy type inequality also holds:
\begin{equation}\label{eq:radial inverse Hardy}
\int_0^\infty\frac{|D_r^{-1}z|^2}{r^2}\;rdr\le\int_0^\infty|z|^2\;rdr.
\end{equation}
If, in addition, $z$ is a smooth radial scalar profile, then
\[D_r^{-1}z=O(r)\]
as $r\to0$.
\end{lemma}

\begin{proof}
Since $sz(s)\in L^1_{\mathrm{loc}}(0,\infty)$, the function
\[rD_r^{-1}z(r)=\int_0^r sz(s) \;ds\]
is locally absolutely continuous. Hence
\[(rD_r^{-1}z)_r=rz\]
almost everywhere, and therefore
\begin{equation}\label{eq:inverse identity}
D_rD_r^{-1}z=\frac1r(rD_r^{-1}z)_r=z
\end{equation}
almost everywhere on $(0,\infty)$.

For a smooth radial vector profile $F$, Remark \ref{rem:behavior at origin} gives $F(r)=O(r)$ as $r\to0$. Hence
\[
D_r^{-1}D_rF = \frac1r\int_0^r s\left(F_s+\frac Fs\right) \; ds = \frac1r\int_0^r(sF(s))_s \; ds = F.
\]
For $\rho>0$, by \eqref{eq:inverse identity},
\[\begin{aligned}\int_0^\rho zD_r^{-1}z\;dr=&~{}\int_0^\rho\left(\partial_rD_r^{-1}z+\frac1rD_r^{-1}z\right)D_r^{-1}z\;dr\\
=&~{}\frac12|D_r^{-1}z(\rho)|^2+\int_0^\rho\frac{|D_r^{-1}z|^2}{r}\;dr.\end{aligned}
\]
Note that $D_r^{-1}z(\rho)\to0$ as $\rho\to0$ for $z\in L^2(rdr)$ (by Cauchy--Schwarz $|D_r^{-1}z(r)|\le\|z\|_{L^2((0,r),rdr)}/\sqrt2$), which is what makes the boundary term at $0$ vanish. On the other hand, Cauchy--Schwarz gives
\[\int_0^\rho zD_r^{-1}z\;dr\le\left(\int_0^\rho z^2\;rdr\right)^{1/2}\left(\int_0^\rho \frac{|D_r^{-1}z|^2}{r}\;dr\right)^{1/2}.\]
Dropping the nonnegative boundary term, we obtain
\[\int_0^\rho\frac{|D_r^{-1}z|^2}{r}\;dr\le\left(\int_0^\rho z^2\;rdr\right)^{1/2}\left(\int_0^\rho\frac{|D_r^{-1}z|^2}{r}\;dr\right)^{1/2}.\]
Hence
\[\int_0^\rho\frac{|D_r^{-1}z|^2}{r}\;dr\le\int_0^\rho z^2\;rdr\le\int_0^\infty z^2\;rdr.\]
Letting $\rho\to\infty$ and using monotone convergence yields \eqref{eq:radial inverse Hardy}.

Finally, if $z$ is a smooth radial scalar profile, then by Remark \ref{rem:behavior at origin},
\[z(r)=z(0)+O(r^2)\]
as $r\to0$. Hence
\[\int_0^r sz(s) \; ds=\frac12z(0)r^2+O(r^4),\]
and therefore, by \eqref{eq:regular radial inverse},
\[D_r^{-1}z(r)=\frac1r\int_0^r sz(s) \;ds=\frac12z(0)r+O(r^3).\]
\end{proof}

For the computation of the virial identity, we set
\begin{equation}\label{eq:inverse divergence}
v:=D_r^{-1}h.
\end{equation}
We now introduce the inverse tangential correction
\begin{equation}\label{eq:inverse tangential correction}
\mathcal K_\phi(t):=2\pi\int_0^\infty\frac{\phi}{r}qD_r^{-1}h\;rdr.
\end{equation}
By Cauchy--Schwarz, \eqref{eq:radial inverse Hardy}, and polar integration,
\begin{equation}\label{eq:K_phi bound}
\begin{aligned}
|\mathcal K_\phi(t)| \le&~{} 2\pi\sup_{r>0}|\phi(r)|\|q\|_{L^2(rdr)} \left\|\frac{D_r^{-1}h}{r}\right\|_{L^2(rdr)}\\
\le&~{}2\pi\sup_{r>0}|\phi(r)|\|q\|_{L^2(rdr)}\|h\|_{L^2(rdr)}\\
=&~{}\sup_{r>0}|\phi(r)|\|u\|_{L^2(\R^2)}\|\eta\|_{L^2(\R^2)}\\
\le&~{}\frac12\sup_{r>0}|\phi(r)|\left(\|u\|_{H^1(\R^2)}^2+\|\eta\|_{H^1(\R^2)}^2\right).
\end{aligned}
\end{equation}
Thus $\mathcal K_\phi(t)$ is well defined and bounded for radial energy solutions, even though $D_r^{-1}h$ itself need not belong to $L^2((0,\infty),rdr)$.

\begin{lemma}[Inverse divergence correction]\label{lem:inverse divergencel virial}
Under the assumptions of Lemma \ref{lem:principal virial identity},
\begin{equation}\label{eq:inverse divergence decomposition}
\frac1{2\pi}\frac{d}{dt}\mathcal K_\phi(t)=\mathcal K_\phi^{(2)}(t)+\mathcal K_\phi^{(3)}(t),
\end{equation}
where
\begin{equation}\label{eq:inverse divergence quadratic}
\begin{aligned}\mathcal K_\phi^{(2)}(t):=&~{}a\int_0^\infty\frac{\phi}{r}q^2\;rdr-(a+1)\int_0^\infty\frac{\phi}{r}q(K_vq)\;rdr\\
&~{}-c\int_0^\infty\frac{\phi}{r}h^2\;rdr+(c+1)\int_0^\infty\frac{\phi}{r}h(K_sh)\;rdr\\
&~{}+\int_0^\infty\left(\frac{\phi}{r}\right)_r\left((c+1)K_sh-ch\right)D_r^{-1}h\;rdr,\end{aligned}
\end{equation}
and
\begin{equation}\label{eq:inverse divergence cubic}
\begin{aligned}\mathcal K_\phi^{(3)}(t):=&~{}-\int_0^\infty\frac{\phi}{r}qK_v(qh)\;rdr+\frac12\int_0^\infty\frac{\phi}{r}hK_s(q^2)\;rdr\\
&~{}+\frac12\int_0^\infty\left(\frac{\phi}{r}\right)_r(D_r^{-1}h)K_s(q^2)\;rdr.\end{aligned}
\end{equation}
\end{lemma}

\begin{proof}
We first compute the time derivative of $D_r^{-1}h$. Using \eqref{eq:radial evolution} and \eqref{eq:resolvent intertwining}, we have
\[h_t=D_r\left(aq-(a+1)K_vq-K_v(qh)\right).\]
Since $aq-(a+1)K_vq-K_v(qh)$ is a smooth radial vector profile, and $\partial_t$ commutes with $D_r^{-1}$, Lemma \ref{lem:regular radial inverse} gives
\begin{equation}\label{eq:inverse profile evolution}
\partial_tD_r^{-1}h=aq-(a+1)K_vq-K_v(qh).
\end{equation}
Differentiating \eqref{eq:inverse tangential correction} in time gives
\[\frac1{2\pi}\frac{d}{dt}\mathcal K_\phi(t)=\int_0^\infty\frac{\phi}{r}q_tD_r^{-1}h\;rdr+\int_0^\infty\frac{\phi}{r}q\partial_tD_r^{-1}h\;rdr.\]
Substituting \eqref{eq:radial evolution} and \eqref{eq:inverse profile evolution}, and using \eqref{eq:resolvent intertwining}, we obtain
\begin{equation}\label{eq:inverse tangential qt expansion}
\begin{aligned}
\int_0^\infty\frac{\phi}{r}q_tD_r^{-1}h\;rdr=&~{}c\int_0^\infty\frac{\phi}{r}h_rD_r^{-1}h\;rdr-(c+1)\int_0^\infty\frac{\phi}{r}(\partial_rK_sh)D_r^{-1}h\;rdr\\
&~{}-\frac12\int_0^\infty\frac{\phi}{r}\partial_rK_s(q^2)D_r^{-1}h\;rdr,\end{aligned}
\end{equation}
and
\begin{equation}\label{eq:inverse tangential inverse evolution term}
\begin{aligned}
\int_0^\infty\frac{\phi}{r}q\partial_tD_r^{-1}h\;rdr=&~{}a\int_0^\infty\frac{\phi}{r}q^2\;rdr-(a+1)\int_0^\infty\frac{\phi}{r}q(K_vq)\;rdr\\
&~{}-\int_0^\infty\frac{\phi}{r}qK_v(qh)\;rdr.\end{aligned}
\end{equation}

We now evaluate the three terms in \eqref{eq:inverse tangential qt expansion}. By the product rule and Lemma \ref{lem:regular radial inverse},
\begin{equation}\label{eq:inverse tangential product identity}
\begin{aligned}
D_r\left(\frac{\phi}{r}D_r^{-1}h\right) =&~{}\partial_r\left(\frac{\phi}{r}D_r^{-1}h\right)+\frac{\phi}{r^2}D_r^{-1}h\\
=&~{}\left(\frac{\phi}{r}\right)_rD_r^{-1}h+\frac{\phi}{r}\left(\partial_rD_r^{-1}h+\frac1rD_r^{-1}h\right)\\
=&~{}\left(\frac{\phi}{r}\right)_rD_r^{-1}h+\frac{\phi}{r}h.
\end{aligned}
\end{equation}
For the first term, \eqref{eq:radial integration by parts} and \eqref{eq:inverse tangential product identity} give
\begin{equation}\label{eq:inverse tangential h term}
\begin{aligned}
c\int_0^\infty\frac{\phi}{r}h_rD_r^{-1}h\;rdr =&~{}-c\int_0^\infty hD_r\left(\frac{\phi}{r}D_r^{-1}h\right)\;rdr\\
=&~{}-c\int_0^\infty\frac{\phi}{r}h^2\;rdr-c\int_0^\infty\left(\frac{\phi}{r}\right)_rhD_r^{-1}h\;rdr.
\end{aligned}
\end{equation}
For the second term, \eqref{eq:radial integration by parts} and \eqref{eq:inverse tangential product identity} give
\[
\begin{aligned}
&~{} -(c+1)\int_0^\infty\frac{\phi}{r}(\partial_rK_sh)D_r^{-1}h\;rdr\\ 
&~{} \quad = (c+1)\int_0^\infty(K_sh)D_r\left(\frac{\phi}{r}D_r^{-1}h\right)\;rdr\\
&~{} \quad = (c+1)\int_0^\infty\frac{\phi}{r}h(K_sh)\;rdr+(c+1)\int_0^\infty\left(\frac{\phi}{r}\right)_r(K_sh)D_r^{-1}h\;rdr.
\end{aligned}
\]
For the third term, the same identities give
\begin{equation}\label{eq:inverse tangential nonlinear term}
\begin{aligned}
&~{}-\frac12\int_0^\infty\frac{\phi}{r}\partial_rK_s(q^2)D_r^{-1}h\;rdr \\
=&~{}\frac12\int_0^\infty K_s(q^2)D_r\left(\frac{\phi}{r}D_r^{-1}h\right)\;rdr\\
=&~{}\frac12\int_0^\infty\frac{\phi}{r}hK_s(q^2)\;rdr+\frac12\int_0^\infty\left(\frac{\phi}{r}\right)_r(D_r^{-1}h)K_s(q^2)\;rdr.
\end{aligned}
\end{equation}
Combining \eqref{eq:inverse tangential h term}--\eqref{eq:inverse tangential nonlinear term} with \eqref{eq:inverse tangential inverse evolution term}, the quadratic terms give \eqref{eq:inverse divergence quadratic}, while the cubic terms give \eqref{eq:inverse divergence cubic}. This proves \eqref{eq:inverse divergence decomposition}.
\end{proof}

The first four terms of \eqref{eq:inverse divergence quadratic} are precisely twice the second integral in \eqref{eq:tangential correction remaining mass terms}. We therefore combine the two corrections with the principal virial functional as
\[\mathcal H_\phi(t):=I_\phi(t)-\frac12\mathcal J_\phi(t)-\frac12\mathcal K_\phi(t).\]
By \eqref{eq:I_phi bound}, \eqref{eq:J_phi bound}, and \eqref{eq:K_phi bound},
\begin{equation}\label{eq:H_phi bound}
|\mathcal H_\phi(t)|\le C\left(\|\phi\|_{L^\infty}+\left\|\frac{\phi}{r}\right\|_{L^\infty}\right)\left(\|u(t)\|_{H^1(\R^2)}^2+\|\eta(t)\|_{H^1(\R^2)}^2\right).
\end{equation}
Recall $v=D_r^{-1}h$ defined in \eqref{eq:inverse divergence}.

\begin{proposition}[Virial identity]\label{prop:Virial identity}
Under the assumptions of Lemma \ref{lem:principal virial identity},
\begin{equation}\label{eq:Virial decomposition}
\frac1{2\pi}\frac{d}{dt}\mathcal H_\phi(t)=\mathcal Q_\phi^{\mathrm{main}}(t)+\mathcal Q_\phi^{\mathrm{err}}(t)+\mathcal N_\phi(t),
\end{equation}
where
\begin{equation}\label{eq:Virial quadratic main}
\begin{aligned}
\mathcal Q_\phi^{\mathrm{main}}:=&~{}\frac12\int_0^\infty\phi_r\left(q^2+h^2-a\left(q_r^2+\frac{q^2}{r^2}\right)-ch_r^2\right)\;rdr\\
&~{}+(a+1)\int_0^\infty\phi_rqL_vK_vq\;rdr
+(c+1)\int_0^\infty\phi_rhL_sK_sh\;rdr\\
&~{}-\frac{3a}{2}\int_0^\infty\left(\frac{\phi}{r}-\phi_r\right)\frac{q^2}{r^2}\;rdr\\
&~{}+\frac14\int_0^\infty\left(\left(\frac{\phi}{r}\right)_{rr}-\frac1r\left(\frac{\phi}{r}\right)_r\right)v^2\;rdr,
\end{aligned}
\end{equation}
and
\begin{equation}\label{eq:Virial quadratic error}
\begin{aligned}
\mathcal Q_\phi^{\mathrm{err}}:=&~{}-\frac{a+1}{2}\int_0^\infty\left(\frac{\phi}{r}\right)_rqD_rK_vq\;rdr -\frac{c+1}{2}\int_0^\infty\left(\frac{\phi}{r}\right)_rh\partial_rK_sh\;rdr\\
&~{}-\frac c4\int_0^\infty L_s\left(\frac{\phi}{r}\right)h^2\;rdr\\
&~{}+\frac{c+1}{2}\int_0^\infty\left(\left(\frac{\phi}{r}\right)_{rr}-\frac1r\left(\frac{\phi}{r}\right)_r\right)v\partial_rK_sh\;rdr,
\end{aligned}
\end{equation}
and
\begin{equation}\label{eq:Virial cubic}
\begin{aligned}
\mathcal N_\phi:=&~{}\int_0^\infty\left(-\frac{\phi_r}{2}+\frac{\phi}{4r}\right)hq^2\;rdr
+\int_0^\infty\phi_rqK_v(qh)\;rdr
+\frac12\int_0^\infty\phi_rhK_s(q^2)\;rdr\\
&~{}-\frac12\int_0^\infty\left(\frac{\phi}{r}\right)_rh\partial_rK_s(q^2)\;rdr
-\frac12\int_0^\infty\left(\frac{\phi}{r}\right)_rqD_rK_v(qh)\;rdr\\
&~{}-\frac14\int_0^\infty\left(\frac{\phi}{r}\right)_rvK_s(q^2)\;rdr.
\end{aligned}
\end{equation}
\end{proposition}

\begin{remark}
The last two terms in \eqref{eq:Virial quadratic main},
\[-\frac{3a}{2}\int_0^\infty\left(\frac{\phi}{r}-\phi_r\right)\frac{q^2}{r^2}\;rdr+\frac14\int_0^\infty\left(\left(\frac{\phi}{r}\right)_{rr}-\frac1r\left(\frac{\phi}{r}\right)_r\right)v^2\;rdr,\]
are nonnegative for the localized weight chosen below (later verified in \eqref{eq:size positivity}). In particular, the second term provides the positive contribution needed to absorb the quadratic error in \eqref{eq:Virial quadratic error},
\[\frac{c+1}{2}\int_0^\infty\left(\left(\frac{\phi}{r}\right)_{rr}-\frac1r\left(\frac{\phi}{r}\right)_r\right)v\partial_rK_sh\;rdr,\]
and the cubic term in \eqref{eq:Virial cubic},
\[-\frac14\int_0^\infty\left(\frac{\phi}{r}\right)_rvK_s(q^2)\;rdr.\]
The precise sign properties of these weights and the corresponding error estimates are established in Sections \ref{sec:local weight} and \ref{sec:coercivity}, respectively.
\end{remark}

\begin{proof}[Proof of Proposition \ref{prop:Virial identity}]
Lemmas \ref{lem:principal virial identity}, \ref{lem:tangential derivative virial}, and \ref{lem:inverse divergencel virial} give
\[\mathcal Q_\phi^{\mathrm{main}}+\mathcal Q_\phi^{\mathrm{err}}=I_\phi^{(2)}-\frac12\mathcal J_\phi^{(2)}-\frac12\mathcal K_\phi^{(2)},\qquad \mathcal N_\phi=I_\phi^{(3)}-\frac12\mathcal J_\phi^{(3)}-\frac12\mathcal K_\phi^{(3)}.\]

We first collect the quadratic terms. In \eqref{eq:tangential derivative quadratic} and \eqref{eq:inverse divergence quadratic}, the terms containing $\frac{\phi}{r}q(K_vq)$ and $\frac{\phi}{r}h(K_sh)$ cancel. Their local mass terms combine as
\[-(a+1)\frac{\phi}{r}q^2+a\frac{\phi}{r}q^2=-\frac{\phi}{r}q^2,\qquad (c+1)\frac{\phi}{r}h^2-c\frac{\phi}{r}h^2=\frac{\phi}{r}h^2.\]
Consequently, the mass terms in $\mathcal Q_\phi^{\mathrm{main}}+\mathcal Q_\phi^{\mathrm{err}}$ become
\[\frac12\left(\phi_r-\frac{\phi}{r}\right)q^2+\frac12\left(\phi_r+\frac{\phi}{r}\right)h^2-\frac12\left(-\frac{\phi}{r}q^2+\frac{\phi}{r}h^2\right)=\frac12\phi_r(q^2+h^2).\]
By \eqref{eq:tangential corrected derivative terms}, the derivative contribution is
\[\begin{aligned}
&~{}-\frac a2\int_0^\infty \phi_rq_r^2\;rdr-\frac c2\int_0^\infty \phi_rh_r^2\;rdr-\frac a2\int_0^\infty\left(3\frac{\phi}{r}-2\phi_r\right)\frac{q^2}{r^2}\;rdr\\
=&~{}-\frac a2\int_0^\infty \phi_r\left(q_r^2+\frac{q^2}{r^2}\right)\;rdr-\frac c2\int_0^\infty \phi_rh_r^2\;rdr-\frac{3a}{2}\int_0^\infty\left(\frac{\phi}{r}-\phi_r\right)\frac{q^2}{r^2}\;rdr.
\end{aligned}\]
The vector resolvent term in \eqref{eq:principal virial quadratic} is
\[(a+1)\int_0^\infty\phi_rqL_vK_vq\;rdr.\]

We next consider the scalar resolvent terms in \eqref{eq:principal virial quadratic}. By \eqref{eq:resolvent intertwining},
\[K_vh_r=\partial_rK_sh,\qquad \partial_rK_vh_r=\partial_{rr}K_sh,\]
which implies
\[\begin{aligned}
&~{}(c+1)\int_0^\infty\phi_rh(\partial_rK_vh_r) + \frac{\phi}{r^2}h(K_vh_r)\;rdr \\
=&~{} (c+1)\int_0^\infty\phi_rh\partial_{rr}K_sh\;rdr+(c+1)\int_0^\infty\frac{\phi}{r}h\partial_rK_sh\;dr.
\end{aligned}\]
By \eqref{eq:scalar radial Laplacian},
\[\partial_{rr}K_sh=L_sK_sh-\frac1r\partial_rK_sh.\]
Hence
\[\begin{aligned}
&~{}(c+1)\int_0^\infty\phi_rh\partial_{rr}K_sh\;rdr+(c+1)\int_0^\infty\frac{\phi}{r}h\partial_rK_sh\;dr\\
=&~{}(c+1)\int_0^\infty\phi_rhL_sK_sh\;rdr-(c+1)\int_0^\infty\phi_rh\partial_rK_sh\;dr+(c+1)\int_0^\infty\frac{\phi}{r}h\partial_rK_sh\;dr\\
=&~{}(c+1)\int_0^\infty\phi_rhL_sK_sh\;rdr+(c+1)\int_0^\infty\left(\frac{\phi}{r}-\phi_r\right)h\partial_rK_sh\;dr.
\end{aligned}\]
Using
\[\left(\frac{\phi}{r}\right)_r=\frac1r\left(\phi_r-\frac{\phi}{r}\right),\]
we conclude
\[\begin{aligned}
&~{}(c+1)\int_0^\infty\phi_rh(\partial_rK_vh_r) + \frac{\phi}{r^2}h(K_vh_r)\;rdr\\
=&~{}(c+1)\int_0^\infty\phi_rhL_sK_sh\;rdr-(c+1)\int_0^\infty\left(\frac{\phi}{r}\right)_rh\partial_rK_sh\;rdr.
\end{aligned}\]

After these combinations, the quadratic contribution takes the intermediate form
\begin{equation}\label{eq:combined quadratic intermediate}
\begin{aligned}
\mathcal Q_\phi^{\mathrm{main}}+\mathcal Q_\phi^{\mathrm{err}}=&~{}\frac12\int_0^\infty \phi_r\left(q^2+h^2-a\left(q_r^2+\frac{q^2}{r^2}\right)-ch_r^2\right)\;rdr\\
&~{}-\frac{3a}{2}\int_0^\infty\left(\frac{\phi}{r}-\phi_r\right)\frac{q^2}{r^2}\;rdr\\
&~{}+(a+1)\int_0^\infty \phi_rqL_vK_vq\;rdr+(c+1)\int_0^\infty \phi_rhL_sK_sh\;rdr\\
&~{}-\frac{a+1}{2}\int_0^\infty \left(\frac{\phi}{r}\right)_rqD_rK_vq\;rdr-(c+1)\int_0^\infty \left(\frac{\phi}{r}\right)_rh\partial_rK_sh\;rdr\\
&~{}-\frac c4\int_0^\infty L_s\left(\frac{\phi}{r}\right)h^2\;rdr-\frac12\int_0^\infty \left(\frac{\phi}{r}\right)_r\left((c+1)K_sh-ch\right)v\;rdr.
\end{aligned}
\end{equation}

By \eqref{eq:resolvent commutation},
\[K_sh=h+L_sK_sh.\]
Hence the last integral splits as
\[-\frac12\int_0^\infty \left(\frac{\phi}{r}\right)_r\left((c+1)K_sh-ch\right)v\;rdr=-\frac12\int_0^\infty \left(\frac{\phi}{r}\right)_rhv\;rdr-\frac{c+1}{2}\int_0^\infty \left(\frac{\phi}{r}\right)_rvL_sK_sh\;rdr.\]
We compute these two terms separately. Using $h=D_rv=v_r+v/r$, integration by parts gives
\begin{equation}\label{eq:inverse tangential positive potential}
\begin{aligned}
-\frac12\int_0^\infty \left(\frac{\phi}{r}\right)_rhv\;rdr
=&~{}-\frac14\int_0^\infty r\left(\frac{\phi}{r}\right)_r(v^2)_r\;dr-\frac12\int_0^\infty \left(\frac{\phi}{r}\right)_rv^2\;dr\\
=&~{}\frac14\int_0^\infty\left(r\left(\frac{\phi}{r}\right)_{rr}-\left(\frac{\phi}{r}\right)_r\right)v^2\;dr\\
=&~{}\frac14\int_0^\infty\left(\left(\frac{\phi}{r}\right)_{rr}-\frac1r\left(\frac{\phi}{r}\right)_r\right)v^2\;rdr.
\end{aligned}
\end{equation}
For the other term, using \eqref{eq:radial factorization} and \eqref{eq:radial integration by parts}, we obtain
\[\int_0^\infty \left(\frac{\phi}{r}\right)_rvL_sK_sh\;rdr=-\int_0^\infty\left(\left(\frac{\phi}{r}\right)_rv\right)_r\partial_rK_sh\;rdr.\]
Since $v_r=h-\frac vr$, we have
\[\left(\left(\frac{\phi}{r}\right)_rv\right)_r=\left(\frac{\phi}{r}\right)_{rr}v+\left(\frac{\phi}{r}\right)_rv_r=\left(\frac{\phi}{r}\right)_rh+\left(\left(\frac{\phi}{r}\right)_{rr}-\frac1r\left(\frac{\phi}{r}\right)_r\right)v.\]
Therefore
\begin{equation}\label{eq:inverse tangential mixed term}
\begin{aligned}
&~{}\int_0^\infty \left(\frac{\phi}{r}\right)_rvL_sK_sh\;rdr \\
=&~{}-\int_0^\infty \left(\frac{\phi}{r}\right)_rh\partial_rK_sh\;rdr -\int_0^\infty\left(\left(\frac{\phi}{r}\right)_{rr}-\frac1r\left(\frac{\phi}{r}\right)_r\right)v\partial_rK_sh\;rdr.
\end{aligned}
\end{equation}
Substituting \eqref{eq:inverse tangential positive potential} and \eqref{eq:inverse tangential mixed term} into \eqref{eq:combined quadratic intermediate} gives \eqref{eq:Virial quadratic main} and \eqref{eq:Virial quadratic error}.

\medskip

We next collect the cubic terms. We first consider the two resolvent terms in \eqref{eq:principal virial cubic}. For the first one, using \eqref{eq:resolvent intertwining}, \eqref{eq:radial factorization}, and \eqref{eq:resolvent commutation},
\begin{equation}\label{eq:combined cubic vector resolvent}
\partial_rK_sD_r(qh)=\partial_rD_rK_v(qh)=L_vK_v(qh)=K_v(qh)-qh.
\end{equation}
For the second resolvent term in \eqref{eq:principal virial cubic}, using \eqref{eq:resolvent intertwining}, \eqref{eq:scalar radial Laplacian}, and \eqref{eq:resolvent commutation},
\begin{equation}\label{eq:combined cubic scalar resolvent}
\partial_rK_v\partial_r(q^2)=\partial_{rr}K_s(q^2)=L_sK_s(q^2)-\frac1r\partial_rK_s(q^2)=K_s(q^2)-q^2-\frac1r\partial_rK_s(q^2).
\end{equation}
Substituting \eqref{eq:combined cubic vector resolvent} and \eqref{eq:combined cubic scalar resolvent} into \eqref{eq:principal virial cubic}, and using $\left(\frac{\phi}{r}\right)_r=\left(\phi_r-\frac{\phi}{r}\right)/r$, we obtain
\begin{equation}\label{eq:main cubic resolvent form}
\begin{aligned}
I_\phi^{(3)}=&~{}-\frac12\int_0^\infty\phi_rhq^2\;rdr+\int_0^\infty\phi_rqK_v(qh)\;rdr+\frac12\int_0^\infty\phi_rhK_s(q^2)\;rdr\\
&~{}-\frac12\int_0^\infty\left(\frac{\phi}{r}\right)_rh\partial_rK_s(q^2)\;rdr.
\end{aligned}
\end{equation}

We next combine the cubic terms from \eqref{eq:tangential derivative cubic} and \eqref{eq:inverse divergence cubic}. The terms containing $\frac{\phi}{r}qK_v(qh)$ and $\frac{\phi}{r}hK_s(q^2)$ cancel, and the remaining terms give
\begin{equation}\label{eq:combined correction cubic terms}
\begin{aligned}
\mathcal J_\phi^{(3)}+\mathcal K_\phi^{(3)}
=&~{}-\frac12\int_0^\infty\frac{\phi}{r}hq^2\;rdr+\int_0^\infty\left(\frac{\phi}{r}\right)_rqD_rK_v(qh)\;rdr\\
&~{}+\frac12\int_0^\infty\left(\frac{\phi}{r}\right)_rvK_s(q^2)\;rdr.
\end{aligned}
\end{equation}
Combining \eqref{eq:main cubic resolvent form} and \eqref{eq:combined correction cubic terms} gives \eqref{eq:Virial cubic}. This proves \eqref{eq:Virial decomposition}.
\end{proof}

\section{Analysis of quadratic terms}\label{sec:quad}

Now we proceed to analyze the terms appearing in the virial identity of Proposition \ref{prop:Virial identity}.

\subsection{Canonical variables}

We first analyze the $\phi_r$ weighted scalar and vector terms in $\mathcal Q_\phi^{\mathrm{main}}$ from \eqref{eq:Virial quadratic main}, namely
\begin{equation}\label{eq:phi_r quadratic}
\begin{aligned}
&~{}\frac12\int_0^\infty\phi_r\left(q^2-a\left(q_r^2+\frac{q^2}{r^2}\right)\right)\;rdr+(a+1)\int_0^\infty\phi_rqL_vK_vq\;rdr\\
&~{}+\frac12\int_0^\infty\phi_r\left(h^2-ch_r^2\right)\;rdr+(c+1)\int_0^\infty\phi_rhL_sK_sh\;rdr.
\end{aligned}
\end{equation}
Introduce the canonical variables
\begin{equation}\label{eq:canonical radial variables}
f:=K_vq,\qquad g:=K_sh.
\end{equation}
Equivalently,
\[q=(1-L_v)f,\qquad h=(1-L_s)g.\]
Hence the resolvent terms in \eqref{eq:phi_r quadratic} become
\[(a+1)\int_0^\infty\phi_rqL_vf\;rdr,\qquad (c+1)\int_0^\infty\phi_rhL_sg\;rdr.\]

\begin{lemma}\label{lem:phi_r canonical quadratic}
Let $\gamma>0$.
\begin{enumerate}
\item[(i)] For the scalar case, with $U=(1-L_s)z$,
\begin{equation}\label{eq:phi_r canonical scalar quadratic}
\begin{aligned}
&~{}\frac12\int_0^\infty\phi_r\left(U^2+\gamma U_r^2\right)\;rdr+(1-\gamma)\int_0^\infty\phi_rUL_sz\;rdr\\
&~{} \quad = \frac12\int_0^\infty\phi_r\left(z^2+3\gamma z_r^2+(4\gamma-1)(L_sz)^2+\gamma(\partial_rL_sz)^2\right)\;rdr\\
&~{}\quad \quad -\frac\gamma2\int_0^\infty\left((L_s\phi_r)z^2+\left(\phi_{rrr}-\frac{\phi_{rr}}{r}\right)z_r^2\right)\;rdr.
\end{aligned}
\end{equation}

\item[(ii)] For the vector case, with $U=(1-L_v)z$,
\begin{equation}\label{eq:phi_r canonical vector quadratic}
\begin{aligned}
&~{}\frac12\int_0^\infty\phi_r\left(U^2+\gamma\left(U_r^2+\frac{U^2}{r^2}\right)\right)\;rdr+(1-\gamma)\int_0^\infty\phi_rUL_vz\;rdr\\
&~{}= \frac12\int_0^\infty\phi_r\left(z^2+3\gamma\left(z_r^2+\frac{z^2}{r^2}\right)+(4\gamma-1)(L_vz)^2+\gamma\left((\partial_rL_vz)^2+\frac{(L_vz)^2}{r^2}\right)\right)\;rdr\\
&~{}\quad -\frac\gamma2\int_0^\infty\left((L_s\phi_r)z^2+\left(\phi_{rrr}-\frac{\phi_{rr}}{r}\right)\left(z_r^2-\frac{z^2}{r^2}\right)\right)\;rdr.
\end{aligned}
\end{equation}
\end{enumerate}
\end{lemma}

\begin{proof}
For the scalar operator, using \eqref{eq:radial factorization}, \eqref{eq:radial integration by parts}, and \eqref{eq:scalar radial Laplacian},
\[\begin{aligned}
\int_0^\infty\phi zL_sz\;rdr =&~{}\int_0^\infty\phi zD_rz_r\;rdr\\
=&~{}-\int_0^\infty\phi z_r^2\;rdr-\int_0^\infty\phi_rzz_r\;rdr\\
=&~{}-\int_0^\infty\phi z_r^2\;rdr+\frac12\int_0^\infty(r\phi_r)_rz^2\;dr\\
=&~{}-\int_0^\infty\phi z_r^2\;rdr+\frac12\int_0^\infty(L_s\phi)z^2\;rdr.
\end{aligned}\]
For the vector operator, using \eqref{eq:vector radial Laplacian}, \eqref{eq:radial integration by parts}, and \eqref{eq:scalar radial Laplacian},
\[\begin{aligned}
\int_0^\infty\phi zL_vz\;rdr =&~{}\int_0^\infty\phi z\left(z_{rr}+\frac1rz_r-\frac1{r^2}z\right)\;rdr\\
=&~{}-\int_0^\infty\phi\left(z_r^2+\frac{z^2}{r^2}\right)\;rdr-\int_0^\infty\phi_rzz_r\;rdr\\
=&~{}-\int_0^\infty\phi\left(z_r^2+\frac{z^2}{r^2}\right)\;rdr+\frac12\int_0^\infty(L_s\phi)z^2\;rdr.
\end{aligned}\]
Hence
\begin{equation}\label{eq:canonical L integration by parts}
\begin{aligned}
\int_0^\infty\phi zL_sz\;rdr=&~{}-\int_0^\infty\phi z_r^2\;rdr+\frac12\int_0^\infty(L_s\phi)z^2\;rdr,\\
\int_0^\infty\phi zL_vz\;rdr=&~{}-\int_0^\infty\phi\left(z_r^2+\frac{z^2}{r^2}\right)\;rdr+\frac12\int_0^\infty(L_s\phi)z^2\;rdr.
\end{aligned}
\end{equation}

We first expand the mass contribution. Using \eqref{eq:canonical L integration by parts}, a direct computation gives, for the scalar case,
\begin{equation}\label{eq:canonical scalar mass expansion}
\int_0^\infty\phi_rU^2\;rdr=\int_0^\infty\phi_r\left(z^2+2z_r^2+(L_sz)^2\right)\;rdr-\int_0^\infty(L_s\phi_r)z^2\;rdr,
\end{equation}
and, for the vector case,
\begin{equation}\label{eq:canonical vector mass expansion}
\int_0^\infty\phi_rU^2\;rdr=\int_0^\infty\phi_r\left(z^2+2\left(z_r^2+\frac{z^2}{r^2}\right)+(L_vz)^2\right)\;rdr-\int_0^\infty(L_s\phi_r)z^2\;rdr.
\end{equation}

We next expand the derivative contribution. In the scalar case,
\[\int_0^\infty\phi_rU_r^2\;rdr=\int_0^\infty\phi_r\left(z_r^2+(\partial_rL_sz)^2\right)\;rdr-2\int_0^\infty\phi_rz_r(L_sz)_r\;rdr.\]
Using \eqref{eq:radial integration by parts}, \eqref{eq:radial factorization}, and
\begin{equation}\label{eq:D_rmz_r}
D_r(\phi_rz_r)=\phi_{rr}z_r+\phi_r\left(z_{rr}+\frac1rz_r\right)=\phi_rL_sz+\phi_{rr}z_r,
\end{equation}
we have
\begin{equation}\label{eq:mz_rL_szr}
\begin{aligned}
-2\int_0^\infty\phi_rz_r(L_sz)_r\;rdr =&~{}2\int_0^\infty D_r(\phi_rz_r)L_sz\;rdr\\
=&~{}2\int_0^\infty\phi_r(L_sz)^2\;rdr+2\int_0^\infty\phi_{rr}z_rL_sz\;rdr,
\end{aligned}
\end{equation}
Hence
\begin{equation}\label{eq:canonical scalar derivative expansion_1}
\int_0^\infty\phi_rU_r^2\;rdr=\int_0^\infty\phi_r\left(z_r^2+2(L_sz)^2+(\partial_rL_sz)^2\right)\;rdr+2\int_0^\infty\phi_{rr}z_rL_sz\;rdr.
\end{equation}
For the vector case,
\[\begin{aligned}
\int_0^\infty\phi_r\left(U_r^2+\frac{U^2}{r^2}\right)\;rdr
=&~{}\int_0^\infty\phi_r\left(z_r^2+\frac{z^2}{r^2}+(\partial_rL_vz)^2+\frac{(L_vz)^2}{r^2}\right)\;rdr\\
&~{}-2\int_0^\infty\phi_rz_r(L_vz)_r\;rdr-2\int_0^\infty\frac{\phi_r}{r^2}zL_vz\;rdr.
\end{aligned}\]
Using \eqref{eq:radial integration by parts} and \eqref{eq:D_rmz_r}, we obtain
\[-2\int_0^\infty\phi_rz_r(L_vz)_r\;rdr=2\int_0^\infty\phi_rL_szL_vz\;rdr+2\int_0^\infty\phi_{rr}z_rL_vz\;rdr.\]
Moreover, by \eqref{eq:scalar radial Laplacian} and \eqref{eq:vector radial Laplacian},
\begin{equation}\label{eq:L_sL_v}
L_sz=L_vz+\frac{z}{r^2}.
\end{equation}
Therefore
\[\begin{aligned}
2\int_0^\infty\phi_rL_szL_vz\;rdr =&~{}2\int_0^\infty\phi_r\left(L_vz+\frac{z}{r^2}\right)L_vz\;rdr\\
=&~{}2\int_0^\infty\phi_r(L_vz)^2\;rdr+2\int_0^\infty\frac{\phi_r}{r^2}zL_vz\;rdr.
\end{aligned}\]
Consequently,
\[-2\int_0^\infty\phi_rz_r(L_vz)_r\;rdr=2\int_0^\infty\phi_r(L_vz)^2\;rdr+2\int_0^\infty\frac{\phi_r}{r^2}zL_vz\;rdr+2\int_0^\infty\phi_{rr}z_rL_vz\;rdr.\]
The second term on the right cancels the remaining mixed term
\[-2\int_0^\infty\frac{\phi_r}{r^2}zL_vz\;rdr.\]
Hence
\begin{equation}\label{eq:canonical vector derivative expansion_1}
\begin{aligned}
\int_0^\infty\phi_r\left(U_r^2+\frac{U^2}{r^2}\right)\;rdr=&~{}\int_0^\infty\phi_r\left(z_r^2+\frac{z^2}{r^2}+2(L_vz)^2+(\partial_rL_vz)^2+\frac{(L_vz)^2}{r^2}\right)\;rdr\\
&~{}+2\int_0^\infty\phi_{rr}z_rL_vz\;rdr.
\end{aligned}
\end{equation}

It remains to compute the last term in \eqref{eq:canonical scalar derivative expansion_1} and \eqref{eq:canonical vector derivative expansion_1}. For the scalar operator, using \eqref{eq:scalar radial Laplacian} and integration by parts,
\begin{equation}\label{eq:m_rz_rL_sz}
\begin{aligned}
2\int_0^\infty\phi_{rr}z_rL_sz\;rdr=&~{}2\int_0^\infty r\phi_{rr}z_rz_{rr}\;dr+2\int_0^\infty\phi_{rr}z_r^2\;dr\\
=&~{}-\int_0^\infty(r\phi_{rr})_rz_r^2\;dr+2\int_0^\infty\phi_{rr}z_r^2\;dr\\
=&~{}-\int_0^\infty\left(\phi_{rrr}-\frac{\phi_{rr}}{r}\right)z_r^2\;rdr.
\end{aligned}
\end{equation}
Substituting this identity into \eqref{eq:canonical scalar derivative expansion_1}, we obtain
\begin{equation}\label{eq:canonical scalar derivative expansion_2}
\ba
\int_0^\infty\phi_rU_r^2\;rdr= &~{} \int_0^\infty\phi_r\left(z_r^2+2(L_sz)^2+(\partial_rL_sz)^2\right)\;rdr \\
&~{} -\int_0^\infty\left(\phi_{rrr}-\frac{\phi_{rr}}{r}\right)z_r^2\;rdr.
\ea
\end{equation}

For the vector operator, using \eqref{eq:L_sL_v} and \eqref{eq:m_rz_rL_sz}, we have
\[\begin{aligned}
2\int_0^\infty\phi_{rr}z_rL_vz\;rdr=&~{}2\int_0^\infty\phi_{rr}z_rL_sz\;rdr-2\int_0^\infty\frac{\phi_{rr}}{r}zz_r\;dr\\
=&~{}-\int_0^\infty\left(\phi_{rrr}-\frac{\phi_{rr}}{r}\right)z_r^2\;rdr-2\int_0^\infty\frac{\phi_{rr}}{r}zz_r\;dr.
\end{aligned}\]
For the last term, integration by parts gives
\[\begin{aligned}
-2\int_0^\infty\frac{\phi_{rr}}{r}zz_r\;dr=&~{}-\int_0^\infty\frac{\phi_{rr}}{r}(z^2)_r\;dr\\
=&~{}\int_0^\infty\left(\frac{\phi_{rr}}{r}\right)_rz^2\;dr\\
=&~{}\int_0^\infty\left(\phi_{rrr}-\frac{\phi_{rr}}{r}\right)\frac{z^2}{r^2}\;rdr.
\end{aligned}\]
Hence
\[2\int_0^\infty\phi_{rr}z_rL_vz\;rdr=-\int_0^\infty\left(\phi_{rrr}-\frac{\phi_{rr}}{r}\right)\left(z_r^2-\frac{z^2}{r^2}\right)\;rdr.\]
Substituting this identity into \eqref{eq:canonical vector derivative expansion_1}, we obtain
\begin{equation}\label{eq:canonical vector derivative expansion_2}
\begin{aligned}
\int_0^\infty\phi_r\left(U_r^2+\frac{U^2}{r^2}\right)\;rdr=&~{}\int_0^\infty\phi_r\left(z_r^2+\frac{z^2}{r^2}+2(L_vz)^2+(\partial_rL_vz)^2+\frac{(L_vz)^2}{r^2}\right)\;rdr\\
&~{}-\int_0^\infty\left(\phi_{rrr}-\frac{\phi_{rr}}{r}\right)\left(z_r^2-\frac{z^2}{r^2}\right)\;rdr.
\end{aligned}
\end{equation}

Finally, since $U=z-L_sz$ in the scalar case, \eqref{eq:canonical L integration by parts} gives
\begin{equation}\label{eq:canonical scalar derivative expansion_3}
\int_0^\infty\phi_rUL_sz\;rdr=-\int_0^\infty\phi_r\left(z_r^2+(L_sz)^2\right)\;rdr+\frac12\int_0^\infty(L_s\phi_r)z^2\;rdr.
\end{equation}
Similarly, since $U=z-L_vz$ in the vector case,
\begin{equation}\label{eq:canonical vector derivative expansion_3}
\int_0^\infty\phi_rUL_vz\;rdr=-\int_0^\infty\phi_r\left(z_r^2+\frac{z^2}{r^2}+(L_vz)^2\right)\;rdr+\frac12\int_0^\infty(L_s\phi_r)z^2\;rdr.
\end{equation}
Substituting \eqref{eq:canonical scalar mass expansion}, \eqref{eq:canonical scalar derivative expansion_2}, and \eqref{eq:canonical scalar derivative expansion_3} into the scalar expression gives \eqref{eq:phi_r canonical scalar quadratic}. Similarly, using \eqref{eq:canonical vector mass expansion}, \eqref{eq:canonical vector derivative expansion_2}, and \eqref{eq:canonical vector derivative expansion_3} gives \eqref{eq:phi_r canonical vector quadratic}.
\end{proof}

\begin{proposition}[Canonical representation of the quadratic terms]
Let
\[f:=K_vq,\qquad g:=K_sh,\qquad v:=D_r^{-1}h.\]
Then the quadratic terms in \eqref{eq:Virial decomposition} admit the following canonical representation:
\begin{equation}\label{eq:canonical representation quadratic}
\begin{aligned}
\mathcal Q_\phi^{\mathrm{main}}=&~{}\frac12\int_0^\infty\phi_r\left(f^2-3a\left(f_r^2+\frac{f^2}{r^2}\right)-(4a+1)(L_vf)^2-a\left((\partial_rL_vf)^2+\frac{(L_vf)^2}{r^2}\right)\right)\;rdr\\
&~{}+\frac12\int_0^\infty\phi_r\left(g^2-3cg_r^2-(4c+1)(L_sg)^2-c(\partial_rL_sg)^2\right)\;rdr\\
&~{}+\frac a2\int_0^\infty\left((L_s\phi_r)f^2+\left(\phi_{rrr}-\frac{\phi_{rr}}r\right)\left(f_r^2-\frac{f^2}{r^2}\right)\right)\;rdr\\
&~{}+\frac c2\int_0^\infty\left((L_s\phi_r)g^2+\left(\phi_{rrr}-\frac{\phi_{rr}}r\right)g_r^2\right)\;rdr\\
&~{}-\frac{3a}{2}\int_0^\infty\left(\frac{\phi}{r}-\phi_r\right)\frac{\left((1-L_v)f\right)^2}{r^2}\;rdr\\
&~{}+\frac14\int_0^\infty\left(\left(\frac{\phi}{r}\right)_{rr}-\frac1r\left(\frac{\phi}{r}\right)_r\right)v^2\;rdr,
\end{aligned}
\end{equation}
and
\begin{equation}\label{eq:canonical quadratic error}
\begin{aligned}
\mathcal Q_\phi^{\mathrm{err}}=&~{}-\frac{a+1}{2}\int_0^\infty\left(\frac{\phi}{r}\right)_r(1-L_v)fD_rf\;rdr\\
&~{}-\frac{c+1}{2}\int_0^\infty\left(\frac{\phi}{r}\right)_r(1-L_s)gg_r\;rdr\\
&~{}-\frac c4\int_0^\infty L_s\left(\frac{\phi}{r}\right)\left((1-L_s)g\right)^2\;rdr\\
&~{}+\frac{c+1}{2}\int_0^\infty\left(\left(\frac{\phi}{r}\right)_{rr}-\frac1r\left(\frac{\phi}{r}\right)_r\right)vg_r\;rdr.
\end{aligned}
\end{equation}
\end{proposition}

\begin{proof}
By \eqref{eq:canonical radial variables},
\[q=(1-L_v)f,\qquad h=(1-L_s)g,\qquad D_rK_vq=D_rf,\qquad\partial_rK_sh=g_r.\]
Apply Lemma \ref{lem:phi_r canonical quadratic} with $(\gamma,z,U)=(-a,f,q)$ in the vector case and $(\gamma,z,U)=(-c,g,h)$ in the scalar case to the first three integrals in \eqref{eq:Virial quadratic main}. Substituting $q=(1-L_v)f$ in its remaining vector term gives \eqref{eq:canonical representation quadratic}. Substituting the identities above into \eqref{eq:Virial quadratic error} gives \eqref{eq:canonical quadratic error}.
\end{proof}

\subsection{Choice of localized weight}\label{sec:local weight}

The quadratic structure in \eqref{eq:canonical representation quadratic} indicates the conditions to be imposed on the weight. The principal local terms require
\[\phi_r>0,\qquad \frac{\phi}{r}-\phi_r\ge0,\]
while the quadratic term involving $v^2$ requires
\[\left(\frac{\phi}{r}\right)_{rr}-\frac1r\left(\frac{\phi}{r}\right)_r\ge0.\]
In addition, the terms carrying derivatives of the weight in \eqref{eq:canonical representation quadratic} and \eqref{eq:canonical quadratic error} require quantitative control of
\begin{equation}\label{eq:error weight}
L_s\phi_r,\qquad \phi_{rrr}-\frac{\phi_{rr}}r,\qquad \left(\frac{\phi}{r}\right)_r,\qquad L_s\left(\frac{\phi}{r}\right),\qquad \left(\frac{\phi}{r}\right)_{rr}-\frac1r\left(\frac{\phi}{r}\right)_r.
\end{equation}
We therefore seek a bounded weight satisfying these sign conditions and introduce a large spatial scale $L$ so that the quantities in \eqref{eq:error weight} are small relative to $\phi_r$. More precisely, under the rescaling $r\mapsto r/L$, each additional radial derivative produces a factor of $L^{-1}$, which allows the quantities in \eqref{eq:error weight} to be treated perturbatively for sufficiently large $L$. At the same time, the scale $L$ determines the spatial region measured by the resulting localized energy.

For this purpose, we choose
\begin{equation}\label{eq:weight choice}
\phi(r):=\int_0^r(1+s^2)^{-\frac34} \;ds,\qquad \phi^L(r):=\phi\left(\frac rL\right).
\end{equation}
Notice that $\|\phi^L\|_{L^\infty}=C_*$ independent of $L$.

\begin{remark}[Choice of the weight]
The form of \eqref{eq:weight choice} is motivated by the relation between the radial and tangential weights. More generally, suppose that $\phi(0)=0$ and $\phi_r$ is positive and decreasing. Then
\[\frac{\phi(r)}r=\frac1r\int_0^r\phi_r(s)ds\ge\phi_r(r)>0.\]
Hence
\[\frac{\phi}{r}-\phi_r\ge0.
\]
We consider the polynomial family prescribed by
\[\phi_r(s)=(1+s^2)^{-\alpha},\qquad\phi(0)=0.\]
The corresponding primitive is bounded provided $\alpha>1/2$. In this range,
\[\phi_r(s)\sim s^{-2\alpha},\qquad\frac{\phi(s)}s=\frac1s\int_0^s(1+\tau^2)^{-\alpha}d\tau\sim s^{-1}\]
as $s\to\infty$, and therefore
\[\frac{\phi(s)}{s\phi_r(s)}\sim s^{2\alpha-1}.\]
The cubic estimates below use the radial pointwise decay from Lemma \ref{l6p4}, which is of order $r^{-1/2}$. Thus, in order that the growth of the tangential weight relative to the radial weight be compatible with this decay, it is natural to require
\[2\alpha-1\le\frac12.\]
Together with $\alpha>1/2$, this gives
\[\frac12<\alpha\le\frac34.\]
We choose the endpoint
\[\alpha=\frac34,\]
which gives the fastest decay of $\phi_r$ in this range and the corresponding critical growth
\[\frac{\phi(s)}{s\phi_r(s)}\sim s^{1/2}.\]
For the weight \eqref{eq:weight choice}, the growth of the tangential weight relative to the radial weight is quantified by \eqref{eq:weight comparison}.
\end{remark}

\begin{lemma}[Properties of the weight]\label{lem:coercive weight bounds}
There is a constant $C$, independent of $L\ge1$, such that the following estimates hold:
\begin{enumerate}
\item[(i)] Size and positivity:
\begin{equation}\label{eq:size positivity}
\begin{aligned}
0<\phi_r^L\le\frac{\phi^L}{r}\le\frac1L,\qquad \|\phi^L\|_{L^\infty}\le C,\\
\left(\frac{\phi^L}{r}\right)_{rr}-\frac1r\left(\frac{\phi^L}{r}\right)_r\ge0.
\end{aligned}
\end{equation}

\item[(ii)] Derivative bounds:
\begin{equation}\label{eq:derivative bounds_1}
|\phi_{rr}^L|+\left|\left(\frac{\phi^L}{r}\right)_r\right|\le\frac CL\phi_r^L,
\end{equation}
and (by \eqref{eq:size positivity})
\begin{equation}\label{eq:derivative bounds_2}
|L_s\phi_r^L|+\left|\phi_{rrr}^L-\frac{\phi_{rr}^L}{r}\right|+\left|L_s\left(\frac{\phi^L}{r}\right)\right| +\left(\frac{\phi^L}{r}\right)_{rr}-\frac1r\left(\frac{\phi^L}{r}\right)_r\le\frac C{L^2}\phi_r^L.
\end{equation}

\item[(iii)] Weight comparison:
\begin{equation}\label{eq:weight comparison}
\begin{aligned}
\left(\left(\frac{\phi^L}{r}\right)_r\right)^2\le&~{}\frac37\frac{\phi^L}{r}\left(\left(\frac{\phi^L}{r}\right)_{rr}-\frac1r\left(\frac{\phi^L}{r}\right)_r\right),\\
\frac{\phi^L}{r\phi_r^L}\le&~{}C\left(1+\frac rL\right)^{1/2}.
\end{aligned}
\end{equation}
\end{enumerate}
The constant $C>0$ may change from line to line and is independent of $L$ and $r$.
\end{lemma}

\begin{proof}
By \eqref{eq:weight choice} and the change of variables $s=tr/L$,
\begin{equation}\label{eq:weight scaling}
\begin{aligned}
&~{}\phi_r^L(r)=\frac1L\left(1+\frac{r^2}{L^2}\right)^{-3/4},\\
&~{}\frac{\phi^L(r)}r=\frac1r\int_0^{r/L}(1+s^2)^{-3/4} \; ds =\frac1L\int_0^1\left(1+t^2\frac{r^2}{L^2}\right)^{-3/4} \; dt.
\end{aligned}
\end{equation}
We first prove \eqref{eq:size positivity}. For $0\le t\le1$,
\[0<\left(1+\frac{r^2}{L^2}\right)^{-3/4}\le\left(1+t^2\frac{r^2}{L^2}\right)^{-3/4}\le1.\]
Integrating with respect to $t$ over $[0,1]$ and using \eqref{eq:weight scaling}, we obtain
\begin{equation}\label{eq:weights bound}
0<\phi_r^L(r)\le\frac{\phi^L(r)}r\le\frac1L.
\end{equation}
We next verify the uniform boundedness of $\phi^L$. By \eqref{eq:weight choice},
\[\phi^L(r)=\int_0^{r/L}(1+s^2)^{-3/4} \; ds.\]
Since $(1+s^2)^{-3/4}\in L^1(0,\infty)$,
\[\|\phi^L\|_{L^\infty}\le\int_0^\infty(1+s^2)^{-3/4} \; ds\le C.\]
For $t\in[0,1]$, direct differentiation gives
\begin{equation}\label{eq:weight derivatives_1}
\begin{aligned}
\partial_r\left(1+t^2\frac{r^2}{L^2}\right)^{-3/4}
=&~{}-\frac{3t^2r}{2L^2}\left(1+t^2\frac{r^2}{L^2}\right)^{-7/4},\\
\left(\partial_r^2-\frac1r\partial_r\right)\left(1+t^2\frac{r^2}{L^2}\right)^{-3/4}
=&~{}\frac{21t^4r^2}{4L^4}\left(1+t^2\frac{r^2}{L^2}\right)^{-11/4}\ge0.
\end{aligned}
\end{equation}
Differentiating the second identity in \eqref{eq:weight scaling} under the integral sign, we obtain
\[\left(\frac{\phi^L(r)}r\right)_r=\frac1L\int_0^1\partial_r\left(1+t^2\frac{r^2}{L^2}\right)^{-3/4} \; dt,\]
and
\begin{equation}\label{eq:weight derivatives_2}
\left(\frac{\phi^L(r)}r\right)_{rr}-\frac1r\left(\frac{\phi^L(r)}r\right)_r=\frac1L\int_0^1\left(\partial_r^2-\frac1r\partial_r\right)\left(1+t^2\frac{r^2}{L^2}\right)^{-3/4} \; dt.
\end{equation}
By \eqref{eq:weight derivatives_1}, the integrand in \eqref{eq:weight derivatives_2} is nonnegative, and therefore
\[\left(\frac{\phi^L(r)}r\right)_{rr}-\frac1r\left(\frac{\phi^L(r)}r\right)_r\ge0.\]
This completes the proof of \eqref{eq:size positivity}.

\medskip

We next prove \eqref{eq:derivative bounds_1} and \eqref{eq:derivative bounds_2}. Using \eqref{eq:weight derivatives_1} with $t=1$, we obtain
\[\phi_{rr}^L(r)=-\frac{3r}{2L^3}\left(1+\frac{r^2}{L^2}\right)^{-7/4}\]
and
\[\phi_{rrr}^L(r)-\frac{\phi_{rr}^L(r)}r=\frac{21r^2}{4L^5}\left(1+\frac{r^2}{L^2}\right)^{-11/4}.\]
Moreover, combining the two identities gives
\[\phi_{rrr}^L(r)+\frac{\phi_{rr}^L(r)}r=-\frac3{L^3}\left(1+\frac{r^2}{L^2}\right)^{-7/4}+\frac{21r^2}{4L^5}\left(1+\frac{r^2}{L^2}\right)^{-11/4}.\]
Dividing these expressions by $\phi_r^L(r)$, we obtain
\[\frac{|\phi_{rr}^L(r)|}{\phi_r^L(r)}=\frac{3r}{2L^2\left(1+r^2/L^2\right)},\qquad \frac{\left|\phi_{rrr}^L(r)-\phi_{rr}^L(r)/r\right|}{\phi_r^L(r)}=\frac{21r^2}{4L^4\left(1+r^2/L^2\right)^2},\]
and
\[\frac{\left|\phi_{rrr}^L(r)+\phi_{rr}^L(r)/r\right|}{\phi_r^L(r)}\le\frac3{L^2\left(1+r^2/L^2\right)}+\frac{21r^2}{4L^4\left(1+r^2/L^2\right)^2}.\]
Since
\[\frac{r/L}{1+r^2/L^2}\le C,\qquad \frac{r^2/L^2}{\left(1+r^2/L^2\right)^2}\le C,\]
we obtain
\begin{equation}\label{eq:weight derivatives_3}
|\phi_{rr}^L(r)|\le\frac CL\phi_r^L(r),\qquad \left|\phi_{rrr}^L(r)-\frac{\phi_{rr}^L(r)}r\right|+\left|\phi_{rrr}^L(r)+\frac{\phi_{rr}^L(r)}r\right|\le\frac C{L^2}\phi_r^L(r).
\end{equation}
Note that
\begin{equation}\label{eq:sign of phi_rr}
\phi_{rr}^L\le0.
\end{equation}
Using \eqref{eq:scalar radial Laplacian},
\[L_s\phi_r^L=\phi_{rrr}^L+\frac{\phi_{rr}^L}{r}.\]
Thus \eqref{eq:weight derivatives_3} gives
\[|L_s\phi_r^L| \le\frac C{L^2}\phi_r^L.\]

We now estimate the derivatives of $\phi^L/r$. Differentiating the second identity in \eqref{eq:weight scaling}, we obtain
\begin{equation}\label{eq:weight derivatives_4}
\left(\frac{\phi^L(r)}r\right)_r=-\frac{3r}{2L^3}\int_0^1t^2\left(1+t^2\frac{r^2}{L^2}\right)^{-7/4} \; dt.
\end{equation}
Differentiating once more,
\begin{equation}\label{eq:weight derivatives_5}
\begin{aligned}
\left(\frac{\phi^L(r)}r\right)_{rr}-\frac1r\left(\frac{\phi^L(r)}r\right)_r =  &~{} \frac{21r^2}{4L^5}\int_0^1t^4\left(1+t^2\frac{r^2}{L^2}\right)^{-11/4} \; dt,\\
L_s\left(\frac{\phi^L(r)}r\right) =&~{} -\frac3{L^3}\int_0^1t^2\left(1+t^2\frac{r^2}{L^2}\right)^{-7/4} \; dt \\
&~{} +\frac{21r^2}{4L^5}\int_0^1t^4\left(1+t^2\frac{r^2}{L^2}\right)^{-11/4}\;dt.
\end{aligned}
\end{equation}
We estimate these expressions separately on $0\le r\le L$ and $r\ge L$. For $0\le r\le L$,
\begin{equation}\label{eq:phi_r lower_1}
\phi_r^L(r)\ge\frac{2^{-3/4}}L.
\end{equation}
Using \eqref{eq:weight derivatives_4} and \eqref{eq:weight derivatives_5},
\[\left|\left(\frac{\phi^L}{r}\right)_r\right|\le \frac{3r}{2L^3}\int_0^1t^2\; dt \le \frac C{L^2},\]
\[\left(\frac{\phi^L}{r}\right)_{rr}-\frac1r\left(\frac{\phi^L}{r}\right)_r\le \frac{21r^2}{4L^5}\int_0^1t^4 \; dt \le \frac C{L^3},\]
and 
\[\left|L_s\left(\frac{\phi^L}{r}\right)\right| \le \frac3{L^3}\int_0^1t^2 \; dt+\frac{21r^2}{4L^5}\int_0^1t^4\;dt \le\frac C{L^3}.\]
Hence by \eqref{eq:phi_r lower_1},
\[\left|\left(\frac{\phi^L}{r}\right)_r\right|\le\frac CL\phi_r^L,\]
and
\[\left(\frac{\phi^L}{r}\right)_{rr}-\frac1r\left(\frac{\phi^L}{r}\right)_r+\left|L_s\left(\frac{\phi^L}{r}\right)\right|\le\frac C{L^2}\phi_r^L\]
on $0\le r\le L$.

\medskip

For $r\ge L$, by \eqref{eq:size positivity},
\begin{equation}\label{eq:phi/r_r bound}
\frac{\phi^L(r)}r\le\frac Cr.
\end{equation}
Moreover, by \eqref{eq:weight scaling},
\[\phi_r^L(r)\ge\frac{2^{-3/4}L^{1/2}}{r^{3/2}},\]
which implies
\begin{equation}\label{eq:phi_r lower_2}
\frac{\phi_r^L(r)}{L^j} \ge \frac{2^{-3/4}}{r^{j+1}}\left(\frac{r}{L}\right)^{j-\frac12} \quad j \ge 1.
\end{equation}
Using
\[\left(\frac{\phi^L}{r}\right)_r=\frac1r\left(\phi_r^L-\frac{\phi^L}{r}\right)\]
and \eqref{eq:phi_r lower_2}, we obtain
\[\left|\left(\frac{\phi^L}{r}\right)_r\right|\le\frac{\phi_r^L}{r}+\frac C{r^2}\le\frac CL\phi_r^L.\]
Also, by \eqref{eq:weight derivatives_3},
\[\frac{|\phi_{rr}^L|}{r}\le\frac C{L^2}\phi_r^L.\]
Using
\[\left(\frac{\phi^L}{r}\right)_{rr}-\frac1r\left(\frac{\phi^L}{r}\right)_r=\frac{\phi_{rr}^L}{r}+\frac3{r^2}\left(\frac{\phi^L}{r}-\phi_r^L\right),\]
together with \eqref{eq:phi/r_r bound} and \eqref{eq:phi_r lower_2}, we obtain
\[\left(\frac{\phi^L}{r}\right)_{rr}-\frac1r\left(\frac{\phi^L}{r}\right)_r\le\frac C{L^2}\phi_r^L.\]
Similarly,
\[L_s\left(\frac{\phi^L}{r}\right)=\frac{\phi_{rr}^L}{r}-\frac1{r^2}\left(\phi_r^L-\frac{\phi^L}{r}\right),\]
and therefore
\[\left|L_s\left(\frac{\phi^L}{r}\right)\right|\le\frac C{L^2}\phi_r^L.\]
Combining the two regions, we conclude that
\[\left|\left(\frac{\phi^L}{r}\right)_r\right|\le\frac CL\phi_r^L,\qquad \left(\frac{\phi^L}{r}\right)_{rr}-\frac1r\left(\frac{\phi^L}{r}\right)_r+\left|L_s\left(\frac{\phi^L}{r}\right)\right|\le\frac C{L^2}\phi_r^L.\]
Together with \eqref{eq:weight derivatives_3}, this proves \eqref{eq:derivative bounds_1} and \eqref{eq:derivative bounds_2}.

\medskip

We finally prove the two estimates in \eqref{eq:weight comparison}. For the first one, \eqref{eq:weight derivatives_1} gives
\[\begin{aligned}
\left|\partial_r\left(1+t^2\frac{r^2}{L^2}\right)^{-3/4}\right|^2=&~{}\frac{9t^4r^2}{4L^4}\left(1+t^2\frac{r^2}{L^2}\right)^{-7/2}\\
=&~{}\frac37\left(1+t^2\frac{r^2}{L^2}\right)^{-3/4}\frac{21t^4r^2}{4L^4}\left(1+t^2\frac{r^2}{L^2}\right)^{-11/4}\\
=&~{}\frac37\left(1+t^2\frac{r^2}{L^2}\right)^{-3/4}\left(\partial_r^2-\frac1r\partial_r\right)\left(1+t^2\frac{r^2}{L^2}\right)^{-3/4}.
\end{aligned}\]
Since $\left(1+t^2r^2/L^2\right)^{-3/4}>0$, by \eqref{eq:weight scaling}, Cauchy--Schwarz, and \eqref{eq:weight derivatives_2},
\[\begin{aligned}
\left|\left(\frac{\phi^L(r)}r\right)_r\right|^2=&~{}\frac1{L^2}\left|\int_0^1\partial_r\left(1+t^2\frac{r^2}{L^2}\right)^{-3/4} \; dt\right|^2\\
\le&~{}\frac1{L^2}\left(\int_0^1\left(1+t^2\frac{r^2}{L^2}\right)^{-3/4} \; dt\right)\left(\int_0^1\frac{\left|\partial_r\left(1+t^2\frac{r^2}{L^2}\right)^{-3/4}\right|^2}{\left(1+t^2\frac{r^2}{L^2}\right)^{-3/4}} \; dt\right)\\
=&~{}\frac37\frac{\phi^L(r)}r\left(\left(\frac{\phi^L(r)}r\right)_{rr}-\frac1r\left(\frac{\phi^L(r)}r\right)_r\right).
\end{aligned}\]

For the second estimate, if $0\le r\le L$, then by \eqref{eq:weights bound} and \eqref{eq:phi_r lower_1},
\[\frac{\phi^L(r)}{r\phi_r^L(r)}\le2^{3/4}.\]
If $r\ge L$, by \eqref{eq:phi/r_r bound} and \eqref{eq:phi_r lower_2} with $j=1$, we obtain
\[\frac{\phi^L(r)}{r\phi_r^L(r)}\le C\left(\frac rL\right)^{1/2}.\]
Therefore
\[\frac{\phi^L(r)}{r\phi_r^L(r)}\le C\left(1+\frac rL\right)^{1/2},\qquad r\ge0.\]
This proves \eqref{eq:weight comparison} and completes the proof.
\end{proof}

We conclude this section with the weighted norm equivalences for the radial resolvents that will be used in the coercivity argument.

\begin{lemma}\label{lem:weighted radial resolvent}
There exist $L_1\ge1$ and $C>0$ such that, for $L\ge L_1$, the following estimates hold. If $U=(1-L_s)z$, then
\begin{equation}\label{eq:norm equi scalar_1}
\begin{aligned}
C^{-1}\int_0^\infty\phi_r^L\left(z^2+z_r^2+(L_sz)^2\right)\;rdr\le&~{}\int_0^\infty\phi_r^LU^2\;rdr\\
\le&~{}C\int_0^\infty\phi_r^L\left(z^2+z_r^2+(L_sz)^2\right)\;rdr,
\end{aligned}
\end{equation}
and
\begin{equation}\label{eq:norm equi scalar_2}
\begin{aligned}
C^{-1}\int_0^\infty\frac{\phi^L}{r}\left(z^2+z_r^2+(L_sz)^2\right)\;rdr\le&~{}\int_0^\infty\frac{\phi^L}{r}U^2\;rdr\\
\le&~{}C\int_0^\infty\frac{\phi^L}{r}\left(z^2+z_r^2+(L_sz)^2\right)\;rdr.
\end{aligned}
\end{equation}
If $U=(1-L_v)z$, then
\begin{equation}\label{eq:norm equi vector_1}
\begin{aligned}
& C^{-1}\int_0^\infty\phi_r^L\left(z^2+z_r^2+\frac{z^2}{r^2}+(L_vz)^2\right)\;rdr \\
&\quad \le \int_0^\infty\phi_r^LU^2\;rdr
\le C\int_0^\infty\phi_r^L\left(z^2+z_r^2+\frac{z^2}{r^2}+(L_vz)^2\right)\;rdr,
\end{aligned}
\end{equation}
and
\begin{equation}\label{eq:norm equi vector_2}
\begin{aligned}
& C^{-1}\int_0^\infty\frac{\phi^L}{r}\left(z^2+z_r^2+\frac{z^2}{r^2}+(L_vz)^2\right)\;rdr \\ 
&\quad \le \int_0^\infty\frac{\phi^L}{r}U^2\;rdr \le C\int_0^\infty\frac{\phi^L}{r}\left(z^2+z_r^2+\frac{z^2}{r^2}+(L_vz)^2\right)\;rdr.
\end{aligned}
\end{equation}
\end{lemma}

\begin{proof}
The upper bounds follow immediately from $U=z-L_jz$ in the corresponding scalar and vector cases. It remains to show the lower bounds.

For the scalar case, replacing the weight in the first identity of \eqref{eq:canonical L integration by parts} by $\phi_r^L$ gives
\[\int_0^\infty\phi_r^LUz\;rdr=\int_0^\infty\phi_r^L\left(z^2+z_r^2\right)\;rdr-\frac12\int_0^\infty(L_s\phi_r^L)z^2\;rdr.\]
Let $C>0$ be the constant in \eqref{eq:derivative bounds_2}, and choose $L_1\ge1$ sufficiently large so that
\[\frac{C}{2L_1^2}\le\frac12.\]
Then, for every $L\ge L_1$,
\[\begin{aligned}
\int_0^\infty\phi_r^LUz\;rdr
\ge&~{}\left(1-\frac{C}{2L^2}\right)\int_0^\infty\phi_r^Lz^2\;rdr+\int_0^\infty\phi_r^Lz_r^2\;rdr\\
\ge&~{}\frac12\int_0^\infty\phi_r^L\left(z^2+z_r^2\right)\;rdr.
\end{aligned}\]
Cauchy--Schwarz then yields
\[\int_0^\infty\phi_r^LUz\;rdr \le \left(\int_0^\infty\phi_r^LU^2\;rdr\right)^{\frac12}\left(\int_0^\infty\phi_r^L(z^2 + z_r^2)\;rdr\right)^{\frac12},\]
which implies
\begin{equation}\label{eq:scalar lower}
\int_0^\infty\phi_r^L\left(z^2+z_r^2\right)\;rdr\le4\int_0^\infty\phi_r^LU^2\;rdr.
\end{equation}
Since $L_sz=z-U$, using \eqref{eq:scalar lower},
\[\int_0^\infty\phi_r^L(L_sz)^2\;rdr\le 2\int_0^\infty\phi_r^LU^2\;rdr + 2\int_0^\infty\phi_r^Lz^2\;rdr \le C\int_0^\infty\phi_r^LU^2\;rdr.\]
This proves the lower bound in \eqref{eq:norm equi scalar_1}.

For the vector case, replacing the weight in the second identity of \eqref{eq:canonical L integration by parts} by $\phi_r^L$ gives 
\[\int_0^\infty\phi_r^LUz\;rdr=\int_0^\infty\phi_r^L\left(z^2+z_r^2+\frac{z^2}{r^2}\right)\;rdr-\frac12\int_0^\infty(L_s\phi_r^L)z^2\;rdr.\]
For the same choice of $L_1$, the same argument gives
\[\int_0^\infty\phi_r^L\left(z^2+z_r^2+\frac{z^2}{r^2}+(L_vz)^2\right)\;rdr\le C\int_0^\infty\phi_r^LU^2\;rdr,\]
which proves \eqref{eq:norm equi vector_1}.

For the weight $\phi^L/r$, \eqref{eq:derivative bounds_2} and \eqref{eq:size positivity} give
\[\left|L_s\left(\frac{\phi^L}{r}\right)\right|\le\frac{C}{L^2}\phi_r^L\le\frac{C}{L^2}\frac{\phi^L}{r}.\]
Hence the same choice of $L_1$ and the same arguments with $\phi_r^L$ replaced by $\phi^L/r$ prove \eqref{eq:norm equi scalar_2} and \eqref{eq:norm equi vector_2}. This completes the proof.
\end{proof}

\section{Local energy decay}\label{sec:MT}

In this section we prove Theorem \ref{thm:MT}. 

\subsection{Coercivity}\label{sec:coercivity}

Recall from \eqref{eq:weight choice} that
\[
\phi_r^L(r)=\frac1L\left(1+\frac{r^2}{L^2}\right)^{-3/4}.
\]
By Lemma \ref{lem:phi_r canonical quadratic}, for $U=(1-L_s)z$ and $\phi=\phi^L$,
\begin{equation}\label{eq:Q main canonical scalar}
\begin{aligned}
&~{}\frac12\int_0^\infty\phi_r^L\left(U^2+\gamma U_r^2\right)\;rdr+(1-\gamma)\int_0^\infty\phi_r^LUL_sz\;rdr\\
=&~{}\frac12\int_0^\infty\phi_r^L\left(z^2+3\gamma z_r^2+(4\gamma-1)(L_sz)^2+\gamma(\partial_rL_sz)^2\right)\;rdr\\
&~{}-\frac\gamma2\int_0^\infty\left((L_s\phi_r^L)z^2+\left(\phi_{rrr}^L-\frac{\phi_{rr}^L}{r}\right)z_r^2\right)\;rdr.
\end{aligned}
\end{equation}
For $U=(1-L_v)z$, the corresponding vector identity is
\begin{equation}\label{eq:Q main canonical vector}
\begin{aligned}
&~{}\frac12\int_0^\infty\phi_r^L\left(U^2+\gamma\left(U_r^2+\frac{U^2}{r^2}\right)\right)\;rdr+(1-\gamma)\int_0^\infty\phi_r^LUL_vz\;rdr\\
=&~{}\frac12\int_0^\infty\phi_r^L\left(z^2+3\gamma\left(z_r^2+\frac{z^2}{r^2}\right)+(4\gamma-1)(L_vz)^2+\gamma\left((\partial_rL_vz)^2+\frac{(L_vz)^2}{r^2}\right)\right)\;rdr\\
&~{}-\frac\gamma2\int_0^\infty\left((L_s\phi_r^L)z^2+\left(\phi_{rrr}^L-\frac{\phi_{rr}^L}{r}\right)\left(z_r^2-\frac{z^2}{r^2}\right)\right)\;rdr.
\end{aligned}
\end{equation}
With $f=K_vq$ and $g=K_sh$ as in \eqref{eq:canonical radial variables}, the first three integrals in $\mathcal Q_{\phi^L}^{\mathrm{main}}$ correspond to $(\gamma,z,U)=(-a,f,q)$ in the vector case and $(-c,g,h)$ in the scalar case. We first establish their lower bounds.

\begin{lemma}\label{lem:canonical coercivity}
Let $\gamma>1/9$. There exist $L_2 =L_2(\gamma)\ge1$ and $\kappa_\gamma>0$ such that, for $L\ge L_2$, the following estimates hold. If $U=(1-L_s)z$, then
\[\frac12\int_0^\infty\phi_r^L\left(U^2+\gamma U_r^2\right)\;rdr+(1-\gamma)\int_0^\infty\phi_r^LUL_sz\;rdr\ge\kappa_\gamma\int_0^\infty\phi_r^L\left(U^2+U_r^2\right)\;rdr.\]
If $U=(1-L_v)z$, then
\[
\ba
& \frac12\int_0^\infty\phi_r^L\left(U^2+\gamma\left(U_r^2+\frac{U^2}{r^2}\right)\right)\;rdr+(1-\gamma)\int_0^\infty\phi_r^LUL_vz\;rdr \\
&\quad \ge\kappa_\gamma\int_0^\infty\phi_r^L\left(U^2+U_r^2+\frac{U^2}{r^2}\right)\;rdr.
\ea
\]
Here $z$ is a smooth real profile of the corresponding scalar or vector type, with sufficient decay to justify the canonical identities above. The constants are independent of $L$ and $z$.
\end{lemma}

\begin{remark}
The condition $\gamma>1/9$ above is sharp. See Lemma~\ref{prop:sharp}. Also, $L_2(\gamma)\ge L_1$ and $c_\gamma=\min\{\frac14,\frac9{64}(\gamma-\frac19)\}$ (and the same in the vector case). This shows $\kappa_\gamma\to0$ as $\gamma\downarrow1/9$, consistent with Lemma~\ref{prop:sharp}.
\end{remark}

\begin{proof}
Observe that, for any $\rho\in\R$,
\[U^2\le2\left(z+\rho L_jz\right)^2+2(1+\rho)^2(L_jz)^2.\]
In the scalar case,
\[U_r^2\le2\left(\partial_r(z+\rho L_sz)\right)^2+2(1+\rho)^2(\partial_rL_sz)^2,\]
while in the vector case,
\[U_r^2+\frac{U^2}{r^2}\le2\left(\left(\partial_r(z+\rho L_vz)\right)^2+\frac{(z+\rho L_vz)^2}{r^2}\right)+2(1+\rho)^2\left((\partial_rL_vz)^2+\frac{(L_vz)^2}{r^2}\right).\]
Thus, to control the weighted energy of $U$, it suffices to obtain positive terms involving $z+\rho L_jz$, $L_jz$, and the corresponding derivatives.

We now determine an appropriate $\rho$. The $\phi_r^L$ weighted principal terms have coefficients
\[1,\qquad 3\gamma,\qquad 4\gamma-1,\qquad \gamma\]
in front of the corresponding terms involving $z$, its derivative, $L_jz$, and the derivative of $L_jz$. On the other hand, the square contribution has coefficients
\[1,\qquad 1-2\rho,\qquad \rho^2-2\rho,\qquad \rho^2,\]
as can be verified from the computations below for the specific choice $\rho=1/3$. Subtracting these coefficients, the remaining ones are
\[3\gamma-1+2\rho,\qquad 4\gamma-(1-\rho)^2,\qquad \gamma-\rho^2.\]
We therefore require
\[3\gamma-1+2\rho\ge0,\qquad 4\gamma-(1-\rho)^2>0,\qquad \gamma-\rho^2>0.\]
Accordingly, we impose
\[\gamma>\max\left\{\frac{1-2\rho}{3},\frac{(1-\rho)^2}{4},\rho^2\right\}.\]
The right hand side is minimized at $\rho=1/3$, where
\[\frac{1-2\rho}{3}=\frac{(1-\rho)^2}{4}=\rho^2=\frac19.\]
We therefore choose $\rho=1/3$, which leads to the condition
\[\gamma>\frac19.\]

With this choice, the estimates at the beginning of the proof become
\begin{equation}\label{eq:scalar U control}
\begin{aligned}
U^2+U_r^2\le&~{}2\left(\left(z+\frac13L_sz\right)^2+\left(\partial_r\left(z+\frac13L_sz\right)\right)^2\right)\\
&~{}+\frac{32}{9}\left((L_sz)^2+(\partial_rL_sz)^2\right)
\end{aligned}
\end{equation}
in the scalar case, and
\begin{equation}\label{eq:vector U control}
\begin{aligned}
U^2+U_r^2+\frac{U^2}{r^2}\le&~{}2\left(\left(z+\frac13L_vz\right)^2+\left(\partial_r\left(z+\frac13L_vz\right)\right)^2+\frac{\left(z+\frac13L_vz\right)^2}{r^2}\right)\\
&~{}+\frac{32}{9}\left((L_vz)^2+(\partial_rL_vz)^2+\frac{(L_vz)^2}{r^2}\right)
\end{aligned}
\end{equation}
in the vector case.

We now compute the weighted decomposition. Recall from \eqref{eq:canonical L integration by parts} that
\[\int_0^\infty\phi_r^LzL_sz\;rdr=-\int_0^\infty\phi_r^Lz_r^2\;rdr+\frac12\int_0^\infty(L_s\phi_r^L)z^2\;rdr,\]
and
\[\int_0^\infty\phi_r^LzL_vz\;rdr=-\int_0^\infty\phi_r^L\left(z_r^2+\frac{z^2}{r^2}\right)\;rdr+\frac12\int_0^\infty(L_s\phi_r^L)z^2\;rdr.\]

For the scalar case, using \eqref{eq:canonical L integration by parts},
\begin{equation}\label{eq:common cal scalar}
\begin{aligned}
\int_0^\infty\phi_r^L\left(z+\frac13L_sz\right)^2\;rdr
=&~{}\int_0^\infty\phi_r^L\left(z^2-\frac23z_r^2+\frac19(L_sz)^2\right)\;rdr\\
&~{}+\frac13\int_0^\infty(L_s\phi_r^L)z^2\;rdr.
\end{aligned}
\end{equation}
On the other hand, by \eqref{eq:mz_rL_szr} and \eqref{eq:m_rz_rL_sz}, with $\phi=\phi^L$,
\[\int_0^\infty\phi_r^Lz_r(\partial_rL_sz)\;rdr=-\int_0^\infty\phi_r^L(L_sz)^2\;rdr+\frac12\int_0^\infty\left(\phi_{rrr}^L-\frac{\phi_{rr}^L}{r}\right)z_r^2\;rdr.\]
Therefore,
\begin{equation}\label{eq:scalar cal}
\begin{aligned}
\int_0^\infty\phi_r^L\left(\partial_r\left(z+\frac13L_sz\right)\right)^2\;rdr
=&~{}\int_0^\infty\phi_r^L\left(z_r^2-\frac23(L_sz)^2+\frac19(\partial_rL_sz)^2\right)\;rdr\\
&~{}+\frac13\int_0^\infty\left(\phi_{rrr}^L-\frac{\phi_{rr}^L}{r}\right)z_r^2\;rdr.
\end{aligned}
\end{equation}
Adding \eqref{eq:common cal scalar} and \eqref{eq:scalar cal}, we obtain
\[\begin{aligned}
&~{}\int_0^\infty\phi_r^L\left(\left(z+\frac13L_sz\right)^2+\left(\partial_r\left(z+\frac13L_sz\right)\right)^2\right)\;rdr\\
=&~{}\int_0^\infty\phi_r^L\left(z^2+\frac13z_r^2-\frac59(L_sz)^2+\frac19(\partial_rL_sz)^2\right)\;rdr\\
&~{}+\frac13\int_0^\infty\left((L_s\phi_r^L)z^2+\left(\phi_{rrr}^L-\frac{\phi_{rr}^L}{r}\right)z_r^2\right)\;rdr.
\end{aligned}\]
Combining this identity with \eqref{eq:Q main canonical scalar}, we obtain
\begin{equation}\label{eq:scalar positivity}
\begin{aligned}
&~{}\frac12\int_0^\infty\phi_r^L\left(U^2+\gamma U_r^2\right)\;rdr+(1-\gamma)\int_0^\infty\phi_r^LUL_sz\;rdr\\
=&~{}\frac12\int_0^\infty\phi_r^L\left(\left(z+\frac13L_sz\right)^2+\left(\partial_r\left(z+\frac13L_sz\right)\right)^2\right)\;rdr\\
&~{}+\frac{\gamma-1/9}{2}\int_0^\infty\phi_r^L\left(3z_r^2+4(L_sz)^2+(\partial_rL_sz)^2\right)\;rdr\\
&~{}-\frac{\gamma+1/3}{2}\int_0^\infty\left((L_s\phi_r^L)z^2+\left(\phi_{rrr}^L-\frac{\phi_{rr}^L}{r}\right)z_r^2\right)\;rdr.
\end{aligned}
\end{equation}

We next consider the vector case. By direct expansion,
\[\begin{aligned}
&~{}\int_0^\infty\phi_r^L\left(\left(\partial_r(z-L_vz)\right)^2+\frac{(z-L_vz)^2}{r^2}\right)\;rdr\\
=&~{}\int_0^\infty\phi_r^L\left(z_r^2+\frac{z^2}{r^2}+(\partial_rL_vz)^2+\frac{(L_vz)^2}{r^2}\right)\;rdr\\
&~{}-2\int_0^\infty\phi_r^L\left(z_r(\partial_rL_vz)+\frac{zL_vz}{r^2}\right)\;rdr.
\end{aligned}\]
On the other hand, \eqref{eq:canonical vector derivative expansion_2}, applied with $U=z-L_vz$, gives
\[\begin{aligned}
&~{}\int_0^\infty\phi_r^L\left(\left(\partial_r(z-L_vz)\right)^2+\frac{(z-L_vz)^2}{r^2}\right)\;rdr\\
=&~{}\int_0^\infty\phi_r^L\left(z_r^2+\frac{z^2}{r^2}+2(L_vz)^2+(\partial_rL_vz)^2+\frac{(L_vz)^2}{r^2}\right)\;rdr\\
&~{}-\int_0^\infty\left(\phi_{rrr}^L-\frac{\phi_{rr}^L}{r}\right)\left(z_r^2-\frac{z^2}{r^2}\right)\;rdr.
\end{aligned}\]
Comparing these two identities, we obtain
\[
\ba
& \int_0^\infty\phi_r^L\left(z_r(\partial_rL_vz)+\frac{zL_vz}{r^2}\right)\;rdr \\
&\quad =-\int_0^\infty\phi_r^L(L_vz)^2\;rdr+\frac12\int_0^\infty\left(\phi_{rrr}^L-\frac{\phi_{rr}^L}{r}\right)\left(z_r^2-\frac{z^2}{r^2}\right)\;rdr.
\ea
\]
Consequently,
\begin{equation}\label{eq:vector cal}
\begin{aligned}
&~{}\int_0^\infty\phi_r^L\left(\left(\partial_r\left(z+\frac13L_vz\right)\right)^2+\frac{\left(z+\frac13L_vz\right)^2}{r^2}\right)\;rdr\\
=&~{}\int_0^\infty\phi_r^L\left(z_r^2+\frac{z^2}{r^2}-\frac23(L_vz)^2+\frac19\left((\partial_rL_vz)^2+\frac{(L_vz)^2}{r^2}\right)\right)\;rdr\\
&~{}+\frac13\int_0^\infty\left(\phi_{rrr}^L-\frac{\phi_{rr}^L}{r}\right)\left(z_r^2-\frac{z^2}{r^2}\right)\;rdr.
\end{aligned}
\end{equation}
Also, by \eqref{eq:canonical L integration by parts},
\begin{equation}\label{eq:common cal vector}
\begin{aligned}
\int_0^\infty\phi_r^L\left(z+\frac13L_vz\right)^2\;rdr
=&~{}\int_0^\infty\phi_r^L\left(z^2-\frac23\left(z_r^2+\frac{z^2}{r^2}\right)+\frac19(L_vz)^2\right)\;rdr\\
&~{}+\frac13\int_0^\infty(L_s\phi_r^L)z^2\;rdr.
\end{aligned}
\end{equation}
Adding \eqref{eq:common cal vector} and \eqref{eq:vector cal}, we obtain
\[\begin{aligned}
&~{}\int_0^\infty\phi_r^L\left(\left(z+\frac13L_vz\right)^2+\left(\partial_r\left(z+\frac13L_vz\right)\right)^2+\frac{\left(z+\frac13L_vz\right)^2}{r^2}\right)\;rdr\\
&~{} = \int_0^\infty\phi_r^L\left(z^2+\frac13\left(z_r^2+\frac{z^2}{r^2}\right)-\frac59(L_vz)^2+\frac19\left((\partial_rL_vz)^2+\frac{(L_vz)^2}{r^2}\right)\right)\;rdr\\
&~{}\quad +\frac13\int_0^\infty\left((L_s\phi_r^L)z^2+\left(\phi_{rrr}^L-\frac{\phi_{rr}^L}{r}\right)\left(z_r^2-\frac{z^2}{r^2}\right)\right)\;rdr.
\end{aligned}\]
Combining this identity with \eqref{eq:Q main canonical vector}, we obtain
\begin{equation}\label{eq:vector positivity}
\begin{aligned}
&~{}\frac12\int_0^\infty\phi_r^L\left(U^2+\gamma\left(U_r^2+\frac{U^2}{r^2}\right)\right)\;rdr+(1-\gamma)\int_0^\infty\phi_r^LUL_vz\;rdr\\
&~{} = \frac12\int_0^\infty\phi_r^L\left(\left(z+\frac13L_vz\right)^2+\left(\partial_r\left(z+\frac13L_vz\right)\right)^2+\frac{\left(z+\frac13L_vz\right)^2}{r^2}\right)\;rdr\\
&~{}\quad +\frac{\gamma-1/9}{2}\int_0^\infty\phi_r^L\left(3\left(z_r^2+\frac{z^2}{r^2}\right)+4(L_vz)^2+(\partial_rL_vz)^2+\frac{(L_vz)^2}{r^2}\right)\;rdr\\
&~{}\quad -\frac{\gamma+1/3}{2}\int_0^\infty\left((L_s\phi_r^L)z^2+\left(\phi_{rrr}^L-\frac{\phi_{rr}^L}{r}\right)\left(z_r^2-\frac{z^2}{r^2}\right)\right)\;rdr.
\end{aligned}
\end{equation}

We now estimate the last terms in \eqref{eq:scalar positivity} and \eqref{eq:vector positivity}. By \eqref{eq:derivative bounds_2}, in the scalar case,
\[\left|\int_0^\infty\left((L_s\phi_r^L)z^2+\left(\phi_{rrr}^L-\frac{\phi_{rr}^L}{r}\right)z_r^2\right)\;rdr\right|\le\frac C{L^2}\int_0^\infty\phi_r^L\left(z^2+z_r^2\right)\;rdr,\]
while in the vector case,
\[\left|\int_0^\infty\left((L_s\phi_r^L)z^2+\left(\phi_{rrr}^L-\frac{\phi_{rr}^L}{r}\right)\left(z_r^2-\frac{z^2}{r^2}\right)\right)\;rdr\right|\le\frac C{L^2}\int_0^\infty\phi_r^L\left(z^2+z_r^2+\frac{z^2}{r^2}\right)\;rdr.\]
By Lemma \ref{lem:weighted radial resolvent}, both errors are bounded by
\[\frac{C_\gamma}{L^2}\int_0^\infty\phi_r^LU^2\;rdr.
\]
Since $\gamma>1/9$, the first two terms on the right hand sides of \eqref{eq:scalar positivity} and \eqref{eq:vector positivity} are nonnegative, and the second one has a strictly positive coefficient. By \eqref{eq:scalar U control}, their sum in the scalar case is bounded below by
\[c_\gamma\int_0^\infty\phi_r^L\left(U^2+U_r^2\right)\;rdr,\]
while by \eqref{eq:vector U control}, their sum in the vector case is bounded below by
\[c_\gamma\int_0^\infty\phi_r^L\left(U^2+U_r^2+\frac{U^2}{r^2}\right)\;rdr\]
for some $c_\gamma>0$. Therefore, choosing $L_2=L_2(\gamma) \ge 1$ sufficiently large so that $C_\gamma/L_2^2\le c_\gamma/2$, we obtain the desired coercive estimates.
\end{proof}

\subsection{Sharpness}
Now we provide an argument showing the sharpness of the hypothesis $\gamma>\frac19$ in Lemma~\ref{lem:canonical coercivity}.

\begin{lemma}\label{prop:sharp}
Let $\gamma>0$.  Let $\phi(r)=r$ in Lemmas~\ref{lem:phi_r canonical quadratic} and \ref{lem:canonical coercivity}. Then the following are satisfied:
\begin{enumerate}
\item If $U=(1-L_s)z$, one has
\[
\frac12\int_0^\infty\big(U^2+\gamma U_r^2\big)\,rdr+(1-\gamma)\int_0^\infty UL_sz\,rdr
=c_0 \int_{\R^2}p_\gamma(|\xi|^2)\,|\hat z(\xi)|^2\,d\xi,
\]
where $\hat z$ is the two-dimensional Fourier transform of the radial function $z$, and $p_\gamma$ is the polynomial
\[
p_\gamma(K)=1+3\gamma K+(4\gamma-1)K^2+\gamma K^3=(1+K)\big(\gamma K^2+(3\gamma-1)K+1\big).
\]
\item The same holds in the vector case with $L_v$ in place of $L_s$. 

\item The polynomial $p_\gamma$ is nonnegative on $[0,\infty)$ if and only if $\gamma\ge\frac19$.
\item In the case $\gamma=\frac19$ one has $p_{1/9}(K)=\frac19(K-3)^2(K+1)$, which vanishes at $K=|\xi|^2=3$. 
\item In the case $\gamma<\frac19$ the form with general $\phi^L_r$ is negative on a band of frequencies around $|\xi|=\sqrt3$. 
\end{enumerate}
\end{lemma}

\begin{remark}
Consequently the coercivity of Lemma \ref{lem:canonical coercivity} cannot hold for any $\gamma<\frac19$, and $\kappa_\gamma\to0$ as $\gamma\downarrow\frac19$.
\end{remark}

\begin{proof}
Under the assumption $\phi_r\equiv1$, the identities in \eqref{eq:canonical L integration by parts} have no terms containing derivatives of weights. Consequently, the operators $L_s$ \eqref{eq:scalar radial Laplacian} and $L_v$ \eqref{eq:vector radial Laplacian} are the radial Laplacians of the scalar and the irrotational radial vector fields. Therefore, using Fourier variables,
\[
\ba
& \int zL_sz\,rdr = - c_0 \int |\xi|^2|\hat z|^2 d\xi, \\
& \int (L_sz)^2rdr = c_0 \int |\xi|^4|\hat z|^2 d\xi, \quad \int(\partial_rL_sz)^2rdr  = c_0 \int |\xi|^6 |\hat z|^2 d\xi,
\ea
\]
and 
\[
\int \left( z_r^2+ \frac{z^2}{r^2} \right)rdr=c_0 \|D(ze_r)\|_{L^2}^2 =  \tilde c_0 \int |\xi|^2 |\widehat{(ze_r)}|^2 d\xi
\]
 in the vector case. Then \eqref{eq:phi_r canonical scalar quadratic} becomes $\frac12\int(1+3\gamma K+(4\gamma-1)K^2+\gamma K^3)|\hat z|^2$ with $K=|\xi|^2$. This proves (1) and (2). 
 
Proof of (3). The quadratic factor is nonnegative on $[0,\infty)$ iff $3\gamma-1\ge0$ or its discriminant $(3\gamma-1)^2-4\gamma=(9\gamma-1)(\gamma-1)$ is $\le0$, i.e.\ iff $\gamma\ge\frac19$. 
 
Proof of (4). Direct. 
 
Proof of (5).  For the localized weight $\phi^L_r$, which is of order $\simeq L^{-1}$ and essentially constant on the ball $r\ll L$, we choose a test profile $z$ which is concentrated at frequency $|\xi|\simeq\sqrt3$ and supported in the physical space $r\ll L$. With this choice, the weighted form becomes negative when the parameter $\gamma$ satisfies $\gamma<\frac19$. Moreover, the errors associated to weighted derivatives  are $O(L^{-2})$. We conclude that no such a $\kappa_\gamma>0$ can exist.
\end{proof}


\subsection{Local energies}

Let us define
\begin{equation}\label{eq:local energy_1}
\mathcal E_L(t):=\int_0^\infty\phi_r^L\left(q^2+q_r^2+\frac{q^2}{r^2}+h^2+h_r^2\right)(t,r)\;rdr.
\end{equation}
This quantity can be recast as a weighted local $H^1$ energy in the original variables. Indeed, for the radial profiles $u=qe_r$ and $\eta=h$,
\[|u|^2=q^2,\qquad |Du|^2=q_r^2+\frac{q^2}{r^2},\qquad |\nabla\eta|^2=h_r^2.\]
Since
\[\phi_r^L(r)=\frac1L\left(1+\frac{r^2}{L^2}\right)^{-3/4},\]
integration in polar coordinates gives
\begin{equation}\label{eq:local energy polar form}
\begin{aligned}
&~{}\int_{\R^2}\left(1+\frac{|x|^2}{L^2}\right)^{-3/4}\left(|u|^2+|Du|^2+\eta^2+|\nabla\eta|^2\right)(t,x)dx\\
&~{} \quad = 2\pi\int_0^\infty\left(1+\frac{r^2}{L^2}\right)^{-3/4}\left(q^2+q_r^2+\frac{q^2}{r^2}+h^2+h_r^2\right)(t,r)\;rdr\\
&~{} \quad =  2\pi L \mathcal E_L(t).
\end{aligned}
\end{equation}
In particular, since the weight in \eqref{eq:local energy polar form} is bounded by $1$, $\mathcal E_L(t)$ is finite for every radial $H^1(\R^2)\times H^1(\R^2)$ profile.

\begin{lemma}\label{lem:quadratic error}
There exists $C_{a,c} > 0$ such that, for every $L\ge L_1$,
\begin{equation}\label{eq:quadratic error}
\mathcal Q_{\phi^L}^{\mathrm{err}}\ge-\frac{C_{a,c}}L \mathcal E_L-\frac18\int_0^\infty\left(\left(\frac{\phi^L}{r}\right)_{rr}-\frac1r\left(\frac{\phi^L}{r}\right)_r\right)v^2\;rdr.
\end{equation}
The constant is independent of $L$ and the profiles. $L_1$ is the threshold from Lemma \ref{lem:weighted radial resolvent}.
\end{lemma}

\begin{proof}
Let $L\ge L_1$. By \eqref{eq:canonical quadratic error} and \eqref{eq:canonical radial variables},
\[\begin{aligned}
\mathcal Q_{\phi^L}^{\mathrm{err}}=&~{}-\frac{a+1}{2}\int_0^\infty\left(\frac{\phi^L}{r}\right)_rqD_rf\;rdr-\frac{c+1}{2}\int_0^\infty\left(\frac{\phi^L}{r}\right)_rhg_r\;rdr\\
&~{}-\frac c4\int_0^\infty L_s\left(\frac{\phi^L}{r}\right)h^2\;rdr\\
&~{}+\frac{c+1}{2}\int_0^\infty\left(\left(\frac{\phi^L}{r}\right)_{rr}-\frac1r\left(\frac{\phi^L}{r}\right)_r\right)vg_r\;rdr.
\end{aligned}\]
We first recall the lower bound of the canonical mass equivalence needed below. Applying Lemma \ref{lem:weighted radial resolvent} with $U=q=(1-L_v)f$, we obtain 
\begin{equation}\label{eq:vector norm upper}
\int_0^\infty\phi_r^L\left(f^2+f_r^2+\frac{f^2}{r^2}+(L_vf)^2\right)\;rdr\le C\int_0^\infty\phi_r^Lq^2\;rdr.
\end{equation}
Since
\[(D_rf)^2=\left(f_r+\frac fr\right)^2\le2\left(f_r^2+\frac{f^2}{r^2}\right),\]
\eqref{eq:vector norm upper} gives
\[\int_0^\infty\phi_r^L\left((D_rf)^2+(L_vf)^2\right)\;rdr\le C\int_0^\infty\phi_r^Lq^2\;rdr.\]
Similarly, applying Lemma \ref{lem:weighted radial resolvent} with $U=h=(1-L_s)g$ gives
\[\int_0^\infty\phi_r^L\left(g^2+g_r^2+(L_sg)^2\right)\;rdr\le C\int_0^\infty\phi_r^Lh^2\;rdr.\]
Combining the last two estimates, we obtain
\begin{equation}\label{eq:norm upper}
\int_0^\infty\phi_r^L\left((D_rf)^2+(L_vf)^2+g_r^2+(L_sg)^2\right)\;rdr\le C\int_0^\infty\phi_r^L(q^2+h^2)\;rdr\le C \mathcal E_L.
\end{equation}
We now estimate the four terms in $\mathcal Q_{\phi^L}^{\mathrm{err}}$. By \eqref{eq:derivative bounds_1} and \eqref{eq:norm upper},
\[
\ba
\left|\int_0^\infty\left(\frac{\phi^L}{r}\right)_rqD_rf\;rdr\right| \le &~{} \frac CL\left(\int_0^\infty\phi_r^Lq^2\;rdr\right)^{1/2}\left(\int_0^\infty\phi_r^L(D_rf)^2\;rdr\right)^{1/2} \\
\le &~{}  \frac CL \mathcal E_L,
\ea
\]
and similarly
\[\left|\int_0^\infty\left(\frac{\phi^L}{r}\right)_rhg_r\;rdr\right| \le\frac CL\left(\int_0^\infty\phi_r^Lh^2\;rdr\right)^{1/2}\left(\int_0^\infty\phi_r^Lg_r^2\;rdr\right)^{1/2} \le\frac CL \mathcal E_L.\]
By \eqref{eq:derivative bounds_2},
\[\left|\int_0^\infty L_s\left(\frac{\phi^L}{r}\right)h^2\;rdr\right|\le\frac C{L^2}\int_0^\infty\phi_r^Lh^2\;rdr\le\frac C{L^2} \mathcal E_L.\]
For the last term, \eqref{eq:size positivity}, Young's inequality, \eqref{eq:derivative bounds_2} and \eqref{eq:norm upper} give
\[\begin{aligned}
&~{}\left|\frac{c+1}{2}\int_0^\infty\left(\left(\frac{\phi^L}{r}\right)_{rr}-\frac1r\left(\frac{\phi^L}{r}\right)_r\right)vg_r\;rdr\right|\\
&~{} \quad \le \frac18\int_0^\infty\left(\left(\frac{\phi^L}{r}\right)_{rr}-\frac1r\left(\frac{\phi^L}{r}\right)_r\right)v^2\;rdr \\
&~{} \quad \quad +\frac{(c+1)^2}{2}\int_0^\infty\left(\left(\frac{\phi^L}{r}\right)_{rr}-\frac1r\left(\frac{\phi^L}{r}\right)_r\right)g_r^2\;rdr\\
&~{}  \quad \le \frac18\int_0^\infty\left(\left(\frac{\phi^L}{r}\right)_{rr}-\frac1r\left(\frac{\phi^L}{r}\right)_r\right)v^2\;rdr+\frac{C(c+1)^2}{L^2} \mathcal E_L.
\end{aligned}\]
Combining these estimates and using $L\ge L_1\ge1$ proves \eqref{eq:quadratic error}.
\end{proof}

\subsection{Control of cubic terms} We next estimate the cubic contribution. Assume
\begin{equation}\label{eq:radial small energy}
\int_0^\infty\left(q^2+q_r^2+\frac{q^2}{r^2}+h^2+h_r^2\right)\;rdr\le\varepsilon^2.
\end{equation}
For the radial profiles $u=qe_r$ and $\eta=h$, polar integration gives
\[\|u\|_{H^1(\R^2)}^2+\|\eta\|_{H^1(\R^2)}^2=2\pi\int_0^\infty\left(q^2+q_r^2+\frac{q^2}{r^2}+h^2+h_r^2\right)\;rdr.\]
Thus \eqref{eq:radial small energy} is precisely a smallness assumption on the radial $H^1(\R^2)\times H^1(\R^2)$ norm, up to the factor $2\pi$.

\begin{lemma}[Radial pointwise bounds]\label{l6p4}
Let $q$ and $h$ be smooth radial profiles with sufficient decay satisfying \eqref{eq:radial small energy}. Then
\begin{equation}\label{eq:radial pointwise bounds}
\sup_{r>0}|q(r)|\le\varepsilon,\qquad |q(r)|+|h(r)|\le C\varepsilon r^{-1/2},\qquad r>0.
\end{equation}
Moreover, uniformly for $L\ge1$,
\begin{equation}\label{eq:nonlinear_1}
\int_0^\infty\frac{\phi^L}{r}|h|q^2\;rdr\le C\varepsilon \mathcal E_L,\qquad \int_0^\infty\frac{\phi^L}{r}q^4\;rdr\le C\varepsilon^2 \mathcal E_L,
\end{equation}
and
\begin{equation}\label{eq:nonlinear_2}
\int_0^\infty\phi_r^Lq^2h^2\;rdr+\int_0^\infty\phi_r^Lq^4\;rdr\le C\varepsilon^2 \mathcal E_L.
\end{equation}
\end{lemma}

\begin{proof}
By Remark \ref{rem:behavior at origin}, $q(0)=0$. Therefore
\[\begin{aligned}
q(r)^2
=&~{}2\int_0^r q(s)q_r(s)ds\\
\le&~{}2\left(\int_0^\infty\frac{q(s)^2}{s}ds\right)^{1/2}\left(\int_0^\infty q_r(s)^2sds\right)^{1/2}\\
\le&~{}\int_0^\infty\left(\frac{q^2}{r^2}+q_r^2\right)\;rdr\le\varepsilon^2.
\end{aligned}\]
Similarly, for $z=q,h$ and $r>0$,
\[\begin{aligned}
|z(r)|^2
\le&~{}2\left(\int_r^\infty z(s)^2sds\right)^{1/2}\left(\int_r^\infty\frac{z_r(s)^2}{s}ds\right)^{1/2}\\
\le&~{}\frac2r\left(\int_r^\infty z(s)^2sds\right)^{1/2}\left(\int_r^\infty z_r(s)^2sds\right)^{1/2}
\le\frac{\varepsilon^2}{r}.
\end{aligned}\]
This proves \eqref{eq:radial pointwise bounds}.

\medskip

We prove the first estimate in \eqref{eq:nonlinear_1} separately on $r\le L$ and $r\ge L$. On $r\le L$, \eqref{eq:weight comparison} gives $\frac{\phi^L}{r}\le C\phi_r^L$. Hence
\[\int_0^L\frac{\phi^L}{r}|h|q^2\;rdr\le C\varepsilon\int_0^L\phi_r^L|hq|\;rdr\le C\varepsilon \mathcal E_L.\]
On $r\ge L$, \eqref{eq:weight comparison} and \eqref{eq:radial pointwise bounds} give
\[\frac{\phi^L}{r}|h|\le C\phi_r^L\left(\frac rL\right)^{1/2}\varepsilon r^{-1/2}=\frac{C\varepsilon}{\sqrt L}\phi_r^L.\]
Consequently,
\[\int_L^\infty\frac{\phi^L}{r}|h|q^2\;rdr\le\frac{C\varepsilon}{\sqrt L}\int_L^\infty\phi_r^Lq^2\;rdr\le C\varepsilon \mathcal E_L.\]
Combining the two regions proves the first estimate.

For the second estimate, \eqref{eq:weight comparison} and \eqref{eq:radial pointwise bounds} give, on $r\le L$,
\[\frac{\phi^L}{r}q^2\le C\varepsilon^2\phi_r^L,\]
while on $r\ge L$,
\[\frac{\phi^L}{r}q^2\le C\phi_r^L\left(\frac rL\right)^{1/2}\frac{\varepsilon^2}{r}=\frac{C\varepsilon^2}{\sqrt{Lr}}\phi_r^L\le C\varepsilon^2\phi_r^L.\]
Multiplying by $q^2$ and integrating proves the second estimate in \eqref{eq:nonlinear_1}. Finally, by \eqref{eq:radial pointwise bounds},
\[\int_0^\infty\phi_r^Lq^2h^2\;rdr+\int_0^\infty\phi_r^Lq^4\;rdr\le\varepsilon^2\int_0^\infty\phi_r^L(h^2+q^2)\;rdr\le\varepsilon^2 \mathcal E_L,\]
which proves \eqref{eq:nonlinear_2}.
\end{proof}

\begin{lemma}\label{lem:nonlinear error}
There exists $C>0$ such that, under \eqref{eq:radial small energy},
\begin{equation}\label{eq:nonlinear error}
|\mathcal N_{\phi^L}|\le C\varepsilon\left(\mathcal E_L+\int_0^\infty\left(\left(\frac{\phi^L}{r}\right)_{rr}-\frac1r\left(\frac{\phi^L}{r}\right)_r\right)v^2\;rdr\right)
\end{equation}
for every $L\ge L_1$, where $L_1$ is the threshold from Lemma \ref{lem:weighted radial resolvent}. The constant is independent of $L$ and the profiles.
\end{lemma}

\begin{proof}
Let $L\ge L_1$. By \eqref{eq:Virial cubic},
\[\begin{aligned}
\mathcal N_{\phi^L}:=&~{}\int_0^\infty\left(-\frac{\phi_r^L}{2}+\frac{\phi^L}{4r}\right)hq^2\;rdr+\int_0^\infty\phi_r^LqK_v(qh)\;rdr+\frac12\int_0^\infty\phi_r^LhK_s(q^2)\;rdr\\
&~{}-\frac12\int_0^\infty\left(\frac{\phi^L}{r}\right)_rh\partial_rK_s(q^2)\;rdr-\frac12\int_0^\infty\left(\frac{\phi^L}{r}\right)_rqD_rK_v(qh)\;rdr\\
&~{}-\frac14\int_0^\infty\left(\frac{\phi^L}{r}\right)_rvK_s(q^2)\;rdr.
\end{aligned}
\]
We first estimate the resolvent terms in $\mathcal N_{\phi^L}$. By \eqref{eq:derivative bounds_1},
\[\left|\left(\frac{\phi^L}{r}\right)_r\right|\le\frac CL\phi_r^L.\]
Applying Lemma \ref{lem:weighted radial resolvent} with the weight $\phi_r^L$ to $z=K_v(qh)$ and $U=qh$, we obtain
\[\int_0^\infty\phi_r^L\left(\left(K_v(qh)\right)^2+\left(\partial_rK_v(qh)\right)^2+\frac{\left(K_v(qh)\right)^2}{r^2}\right)\;rdr\le C\int_0^\infty\phi_r^Lq^2h^2\;rdr.\]
Since
\[\left(D_rK_v(qh)\right)^2\le2\left(\left(\partial_rK_v(qh)\right)^2+\frac{\left(K_v(qh)\right)^2}{r^2}\right),\]
it follows that
\[\int_0^\infty\phi_r^L\left(\left(K_v(qh)\right)^2+\left(D_rK_v(qh)\right)^2\right)\;rdr\le C\int_0^\infty\phi_r^Lq^2h^2\;rdr.\]

Similarly, applying Lemma \ref{lem:weighted radial resolvent} with the weight $\phi_r^L$ to $z=K_s(q^2)$ and $U=q^2$, we obtain
\[\int_0^\infty\phi_r^L\left(\left(K_s(q^2)\right)^2+\left(\partial_rK_s(q^2)\right)^2\right)\;rdr\le C\int_0^\infty\phi_r^Lq^4\;rdr.\]
Combining these estimates with \eqref{eq:nonlinear_2}, we obtain
\begin{equation}\label{eq:nonlinear error_1}
\begin{aligned}
&~{}\int_0^\infty\phi_r^L\left(\left(K_v(qh)\right)^2+\left(D_rK_v(qh)\right)^2+\left(K_s(q^2)\right)^2+\left(\partial_rK_s(q^2)\right)^2\right)\;rdr\\
&~{}\qquad \le C\int_0^\infty\phi_r^L\left(q^2h^2+q^4\right)\;rdr
\le C\varepsilon^2\mathcal E_L.
\end{aligned}
\end{equation}

The estimate \eqref{eq:nonlinear error_1} allows us to control the corresponding mixed cubic terms by Young's inequality. Indeed,
\[\left|qK_v(qh)\right|\le\frac{\varepsilon}{2}q^2+\frac1{2\varepsilon}\left(K_v(qh)\right)^2,\qquad \left|hK_s(q^2)\right|\le\frac{\varepsilon}{2}h^2+\frac1{2\varepsilon}\left(K_s(q^2)\right)^2.\]
Moreover, by \eqref{eq:derivative bounds_1},
\[\begin{aligned}
\left|\left(\frac{\phi^L}{r}\right)_rh\partial_rK_s(q^2)\right|
\le&~{}\frac CL\phi_r^L|h|\left|\partial_rK_s(q^2)\right|,\\
\left|\left(\frac{\phi^L}{r}\right)_rqD_rK_v(qh)\right|
\le&~{}\frac CL\phi_r^L|q|\left|D_rK_v(qh)\right|.
\end{aligned}\]
Applying Young's inequality with parameter $\varepsilon$, the four terms above are bounded by a constant multiple of
\[\varepsilon\phi_r^L(q^2+h^2)+\frac1\varepsilon\phi_r^L\left(\left(K_v(qh)\right)^2+\left(D_rK_v(qh)\right)^2+\left(K_s(q^2)\right)^2+\left(\partial_rK_s(q^2)\right)^2\right).\]
Therefore, by \eqref{eq:nonlinear error_1},
\[\left|\int_0^\infty\phi_r^LqK_v(qh)\;rdr\right|+\left|\int_0^\infty\phi_r^LhK_s(q^2)\;rdr\right|\le C\varepsilon\mathcal E_L,\]
and
\[
\left|\int_0^\infty\left(\frac{\phi^L}{r}\right)_rh\partial_rK_s(q^2)\;rdr\right|+\left|\int_0^\infty\left(\frac{\phi^L}{r}\right)_rqD_rK_v(qh)\;rdr\right|\le\frac{C\varepsilon}{L}\mathcal E_L.
\]
For the last term in $\mathcal N_{\phi^L}$, by \eqref{eq:weight comparison},
\[\left|\left(\frac{\phi^L}{r}\right)_r\right|\le\sqrt{\frac37}\left(\frac{\phi^L}{r}\right)^{1/2}\left(\left(\frac{\phi^L}{r}\right)_{rr}-\frac1r\left(\frac{\phi^L}{r}\right)_r\right)^{1/2}.\]
Hence, by Cauchy--Schwarz,
\[\begin{aligned}
&~{}\left|\int_0^\infty\left(\frac{\phi^L}{r}\right)_rvK_s(q^2)\;rdr\right|\\
\le&~{}\sqrt{\frac37}\left(\int_0^\infty\left(\left(\frac{\phi^L}{r}\right)_{rr}-\frac1r\left(\frac{\phi^L}{r}\right)_r\right)v^2\;rdr\right)^{1/2}\left(\int_0^\infty\frac{\phi^L}{r}\left(K_s(q^2)\right)^2\;rdr\right)^{1/2}.
\end{aligned}\]
Applying Lemma \ref{lem:weighted radial resolvent} with the weight $\phi^L/r$ to $z=K_s(q^2)$ and $U=q^2$, and using \eqref{eq:nonlinear_1}, we obtain
\[\int_0^\infty\frac{\phi^L}{r}\left(K_s(q^2)\right)^2\;rdr\le C\int_0^\infty\frac{\phi^L}{r}q^4\;rdr\le C\varepsilon^2\mathcal E_L.\]
Thus, together with Young's inequality, we obtain
\[
\begin{aligned}
\left|\int_0^\infty\left(\frac{\phi^L}{r}\right)_rvK_s(q^2)\;rdr\right|\le&~{}C\varepsilon\mathcal E_L^{1/2}\left(\int_0^\infty\left(\left(\frac{\phi^L}{r}\right)_{rr}-\frac1r\left(\frac{\phi^L}{r}\right)_r\right)v^2\;rdr\right)^{1/2}\\
\le&~{}C\varepsilon\left(\mathcal E_L+\int_0^\infty\left(\left(\frac{\phi^L}{r}\right)_{rr}-\frac1r\left(\frac{\phi^L}{r}\right)_r\right)v^2\;rdr\right).
\end{aligned}
\]
Now we consider the pure cubic term in $\mathcal N_{\phi^L}$. By \eqref{eq:size positivity} and \eqref{eq:nonlinear_1},
\[
\left|\int_0^\infty\left(-\frac{\phi_r^L}{2}+\frac{\phi^L}{4r}\right)hq^2\;rdr\right|\le\frac34\int_0^\infty\frac{\phi^L}{r}|h|q^2\;rdr\le C\varepsilon\mathcal E_L.
\]
Combining these estimates and using $L\ge L_1$ proves \eqref{eq:nonlinear error}.
\end{proof}

\begin{proposition}[Coercivity]\label{prop:coercivity}
Assume in the normalized system $b=d=1$
\begin{equation}\label{eq:parameter condition}
-a>\frac19,\qquad -c>\frac19.
\end{equation}
There exists $\varepsilon_0>0$, only depending on $a,c$, and $L_0\ge1$, $\kappa>0$, such that the following holds. Let $(u,\eta)$ be a global radial solution satisfying
\[\sup_{t\in\R}\left(\|u(t)\|_{H^1(\R^2)}+\|\eta(t)\|_{H^1(\R^2)}\right)\le\varepsilon\le\varepsilon_0.\]
Then, for every fixed $L\ge L_0$ and $t\in\R$,
\begin{equation}\label{eq:coercivity}
\frac1{2\pi}\frac{d}{dt}\mathcal H_{\phi^L}(t)\ge\kappa\mathcal E_L(t).
\end{equation}
\end{proposition}

\begin{proof}
Apply Lemma \ref{lem:canonical coercivity} with $\gamma=-a$, $z=f$ in the vector case and $\gamma=-c$, $z=g$ in the scalar case. Since $(1-L_v)f=q$ and $(1-L_s)g=h$, there exists $\kappa_1>0$ such that
\[\begin{aligned}
&~{}\frac12\int_0^\infty\phi_r^L\left(q^2-a\left(q_r^2+\frac{q^2}{r^2}\right)+h^2-ch_r^2\right)\;rdr\\
&~{}+(a+1)\int_0^\infty\phi_r^LqL_vf\;rdr+(c+1)\int_0^\infty\phi_r^LhL_sg\;rdr \ge\kappa_1\mathcal E_L
\end{aligned}\]
for every
\[L\ge\max\{L_2(-a),L_2(-c)\}.\] 
By \eqref{eq:Virial quadratic main} and \eqref{eq:canonical radial variables},
\[
\ba
\mathcal Q_{\phi^L}^{\mathrm{main}}\ge &~{} \kappa_1\mathcal E_L-\frac{3a}{2}\int_0^\infty\left(\frac{\phi^L}{r}-\phi_r^L\right)\frac{q^2}{r^2}\;rdr \\
&~{} +\frac14\int_0^\infty\left(\left(\frac{\phi^L}{r}\right)_{rr}-\frac1r\left(\frac{\phi^L}{r}\right)_r\right)v^2\;rdr.
\ea
\]
Combining this estimate with Lemma \ref{lem:quadratic error}, we obtain
\[\begin{aligned}
\mathcal Q_{\phi^L}^{\mathrm{main}}+\mathcal Q_{\phi^L}^{\mathrm{err}}\ge&~{}\left(\kappa_1-\frac{C_{a,c}}L\right)\mathcal E_L-\frac{3a}{2}\int_0^\infty\left(\frac{\phi^L}{r}-\phi_r^L\right)\frac{q^2}{r^2}\;rdr\\
&~{}+\frac18\int_0^\infty\left(\left(\frac{\phi^L}{r}\right)_{rr}-\frac1r\left(\frac{\phi^L}{r}\right)_r\right)v^2\;rdr.
\end{aligned}\]
Choose
\[L_0\ge\max\{L_1,L_2(-a),L_2(-c)\}\]
sufficiently large so that
\[\frac{C_{a,c}}{L_0}\le\frac{\kappa_1}{2}.\] 
Since $a<0$ and \eqref{eq:size positivity} gives $\phi^L/r-\phi_r^L\ge0$, we obtain
\[\mathcal Q_{\phi^L}^{\mathrm{main}}+\mathcal Q_{\phi^L}^{\mathrm{err}}\ge\frac{\kappa_1}{2}\mathcal E_L+\frac18\int_0^\infty\left(\left(\frac{\phi^L}{r}\right)_{rr}-\frac1r\left(\frac{\phi^L}{r}\right)_r\right)v^2\;rdr.\]
On the other hand, Lemma \ref{lem:nonlinear error} gives
\[|\mathcal N_{\phi^L}|\le C\varepsilon\left(\mathcal E_L+\int_0^\infty\left(\left(\frac{\phi^L}{r}\right)_{rr}-\frac1r\left(\frac{\phi^L}{r}\right)_r\right)v^2\;rdr\right).\]
Combining these estimates with \eqref{eq:Virial decomposition}, we obtain
\[\frac1{2\pi}\frac{d}{dt}\mathcal H_{\phi^L}\ge\left(\frac{\kappa_1}{2}-C\varepsilon\right)\mathcal E_L +\left(\frac18-C\varepsilon\right)\int_0^\infty\left(\left(\frac{\phi^L}{r}\right)_{rr}-\frac1r\left(\frac{\phi^L}{r}\right)_r\right)v^2\;rdr.\]
Choose $\varepsilon_0>0$ so that
\[C\varepsilon_0\le\min\left\{\frac{\kappa_1}{4},\frac1{16}\right\},\]
and set $\kappa=\kappa_1/4$. Hence, by \eqref{eq:size positivity},
\[\frac1{2\pi}\frac{d}{dt}\mathcal H_{\phi^L}(t)\ge\kappa\mathcal E_L(t),\]
which proves \eqref{eq:coercivity}.
\end{proof}

\begin{corollary}[Time integrability of the local energy]\label{cor:local energy integrability}
Assume \eqref{eq:parameter condition}. Let $\varepsilon_0>0$ and $L_0\ge1$ be as in Proposition \ref{prop:coercivity}. Let $(u,\eta)$ be a global radial solution satisfying
\[\sup_{t\in\R}\left(\|u(t)\|_{H^1(\R^2)}+\|\eta(t)\|_{H^1(\R^2)}\right)\le\varepsilon\le\varepsilon_0.\]
Then, there is $C,L_0>0$ such that, for every fixed $L\ge L_0$, 
\begin{equation}\label{eq:local energy integrability}
\int_{\R}\int_{\R^2}\left(1+\frac{|x|^2}{L^2}\right)^{-3/4}\left(|u|^2+|Du|^2+\eta^2+|\nabla\eta|^2\right)(t,x)dxdt\le CL\varepsilon^2.
\end{equation}
Moreover, there exist sequences $t_k^\pm\to\pm\infty$ such that
\begin{equation}\label{eq:local energy liminf}
\int_{\R^2}\left(1+\frac{|x|^2}{L^2}\right)^{-3/4}\left(|u|^2+|Du|^2+\eta^2+|\nabla\eta|^2\right)(t_k^\pm,x)dx\to0
\end{equation}
as $k\to\infty$.
\end{corollary}

\begin{proof}
Integrating \eqref{eq:coercivity} over $[-T,T]$ and using \eqref{eq:H_phi bound}, we obtain
\[\kappa\int_{-T}^T\mathcal E_L(t)dt\le\frac1{2\pi}\left(\mathcal H_{\phi^L}(T)-\mathcal H_{\phi^L}(-T)\right)\le C\varepsilon^2.\]
Since $\mathcal E_L\ge0$, letting $T\to\infty$ gives
\begin{equation}\label{eq:E_L integrable}
\int_{\R}\mathcal E_L(t)dt\le C\varepsilon^2.
\end{equation}
Multiplying by $2\pi L$ and using \eqref{eq:local energy polar form}, we obtain \eqref{eq:local energy integrability}. Since $\mathcal E_L$ is nonnegative and integrable on $\R$,
\[\liminf_{t\to\pm\infty}\mathcal E_L(t)=0.\]
Hence there exist sequences $t_k^\pm\to\pm\infty$ such that
\begin{equation}\label{eq:E_L liminf}
\mathcal E_L(t_k^\pm)\to0.
\end{equation}
Using \eqref{eq:local energy polar form} once more proves \eqref{eq:local energy liminf}.
\end{proof}

\subsection{Local energy decay}

For fixed $L\ge L_0$, define
\begin{equation}\label{eq:local energy_2}
\mathscr E_L(t):=\frac12\int_0^\infty\phi_r^L\left(-a\left(q^2+q_r^2+\frac{q^2}{r^2}\right)-c\left(h^2+h_r^2\right)\right)(t,r)\;rdr.
\end{equation}
By \eqref{eq:parameter condition} and \eqref{eq:local energy_1}, we have the bounds
\begin{equation}\label{eq:local energy comparison}
\frac12\min\{-a,-c\}\mathcal E_L(t)\le\mathscr E_L(t)\le\frac12\max\{-a,-c\}\mathcal E_L(t).
\end{equation}
Let us compute the time evolution of $\mathscr E_L(t)$.

\begin{lemma}\label{lem:local energy identity}
Let $f$ and $g$ be as in \eqref{eq:canonical radial variables}. For every fixed $L\ge1$,
\begin{equation}\label{eq:local energy identity}
\begin{aligned}
\frac{d}{dt}\mathscr E_L
=&~{}ac\int_0^\infty\phi_{rr}^L(qh+q_rh_r)\;rdr\\
&~{}+a(c+1)\int_0^\infty\phi_r^L\left(qg_r+q_rg_{rr}+\frac{qg_r}{r^2}\right)\;rdr\\
&~{}+c(a+1)\int_0^\infty\phi_r^L\left(hD_rf+h_rL_vf\right)\;rdr\\
&~{}+\frac a2\int_0^\infty\phi_r^L\left(q\partial_rK_s(q^2)+q_r\partial_{rr}K_s(q^2)+\frac q{r^2}\partial_rK_s(q^2)\right)\;rdr\\
&~{}+c\int_0^\infty\phi_r^L\left(hD_rK_v(qh)+h_rL_vK_v(qh)\right)\;rdr.
\end{aligned}
\end{equation}
\end{lemma}

\begin{proof}
By \eqref{eq:radial evolution}, \eqref{eq:radial derivative evolution}, and \eqref{eq:canonical radial variables}, we have for $q$
\begin{equation}\label{eq:q derivative}
\ba
& q_t=ch_r-(c+1)g_r-\frac12\partial_rK_s(q^2) \quad \text{and} \quad \\
& q_{tr}=ch_{rr}-(c+1)g_{rr}-\frac12\partial_{rr}K_s(q^2),
\ea
\end{equation}
and for $h$
\begin{equation}\label{eq:h derivative}
\ba
& h_t=aD_rq-(a+1)D_rf-D_rK_v(qh) \quad \text{and} \quad \\
& h_{tr}=aL_vq-(a+1)L_vf-L_vK_v(qh).
\ea
\end{equation}
Differentiating \eqref{eq:local energy_2} gives
\[\begin{aligned}
\frac{d}{dt}\mathscr E_L
=&~{}-a\int_0^\infty\phi_r^Lqq_t\;rdr-a\int_0^\infty\phi_r^Lq_rq_{tr}\;rdr-a\int_0^\infty\phi_r^L\frac{q}{r^2}q_t\;rdr\\
&~{}-c\int_0^\infty\phi_r^Lhh_t\;rdr-c\int_0^\infty\phi_r^Lh_rh_{tr}\;rdr.
\end{aligned}\]
Substituting \eqref{eq:q derivative} into the first three terms yields
\[\begin{aligned}
&~{}-ac\int_0^\infty\phi_r^L\left(qh_r+q_rh_{rr}+\frac{q}{r^2}h_r\right)\;rdr\\
&~{}+a(c+1)\int_0^\infty\phi_r^L\left(qg_r+q_rg_{rr}+\frac{qg_r}{r^2}\right)\;rdr\\
&~{}+\frac a2\int_0^\infty\phi_r^L\left(q\partial_rK_s(q^2)+q_r\partial_{rr}K_s(q^2)+\frac q{r^2}\partial_rK_s(q^2)\right)\;rdr.
\end{aligned}\]
Similarly, substituting \eqref{eq:h derivative} into the last two terms gives
\[\begin{aligned}
=&~{}-ac\int_0^\infty\phi_r^L\left(hD_rq+h_rL_vq\right)\;rdr\\
&~{}+c(a+1)\int_0^\infty\phi_r^L\left(hD_rf+h_rL_vf\right)\;rdr\\
&~{}+c\int_0^\infty\phi_r^L\left(hD_rK_v(qh)+h_rL_vK_v(qh)\right)\;rdr.
\end{aligned}\]
Hence, collecting all terms containing no resolvent profiles, we obtain
\[-ac\int_0^\infty\phi_r^L\left(qh_r+hD_rq+q_rh_{rr}+\frac{q}{r^2}h_r+h_rL_vq\right)\;rdr.\]
The first two terms have the divergence form
\[qh_r+hD_rq=qh_r+hq_r+\frac{qh}{r}=D_r(qh).\]
For the remaining three terms, \eqref{eq:vector radial Laplacian} gives
\[\begin{aligned}
q_rh_{rr}+\frac{q}{r^2}h_r+h_rL_vq
=&~{}q_rh_{rr}+\frac{q}{r^2}h_r+h_r\left(q_{rr}+\frac{q_r}{r}-\frac{q}{r^2}\right)\\
=&~{}q_rh_{rr}+q_{rr}h_r+\frac{q_rh_r}{r}
=D_r(q_rh_r).
\end{aligned}\]
Therefore, by \eqref{eq:radial integration by parts}, the whole local contribution equals
\[-ac\int_0^\infty\phi_r^LD_r(qh+q_rh_r)\;rdr=ac\int_0^\infty\phi_{rr}^L(qh+q_rh_r)\;rdr.\]
Combining this identity with the four remaining resolvent contributions proves \eqref{eq:local energy identity}.
\end{proof}

We are now ready to prove Theorem \ref{thm:MT}, stated as follows in the normalized case where $b=d=1$. The original case in Theorem \ref{thm:MT} is obtained by direct scaling.

\begin{theorem}[Local energy decay]\label{thm:local energy decay}
Assume \eqref{eq:parameter condition}. Let $\varepsilon_0>0$ and $L_0\ge1$ be as in Proposition \ref{prop:coercivity}. Let $(u,\eta)$ be a global radial solution satisfying
\[\sup_{t\in\R}\left(\|u(t)\|_{H^1(\R^2)}+\|\eta(t)\|_{H^1(\R^2)}\right)\le\varepsilon\le\varepsilon_0.\]
Then, for every fixed $L\ge L_0$,
\[\lim_{t\to\pm\infty}\int_{\R^2}\left(1+\frac{|x|^2}{L^2}\right)^{-3/4}\left(|u|^2+|Du|^2+\eta^2+|\nabla\eta|^2\right)(t,x)dx=0.\]
Consequently, for every $R >0$, with $B_R=\{x\in\R^2:|x|<R\}$,
\begin{equation}\label{eq:local H1 decay}
\lim_{t\to\pm\infty}\left(\|u(t)\|_{H^1(B_R)}+\|\eta(t)\|_{H^1(B_R)}\right)=0.
\end{equation}
\end{theorem}

\begin{proof}
Fix $L\ge L_0$. By Lemma \ref{lem:local energy identity}, it suffices to estimate the five terms in \eqref{eq:local energy identity}.

By \eqref{eq:derivative bounds_1} and Cauchy--Schwarz,
\[\begin{aligned}
\left|\int_0^\infty\phi_{rr}^L(qh+q_rh_r)\;rdr\right| \le&~{}\frac CL\int_0^\infty\phi_r^L\left(|qh|+|q_rh_r|\right)\;rdr\\
\le&~{}\frac CL
\left(\int_0^\infty\phi_r^L(q^2+q_r^2)\;rdr\right)^{1/2}
\left(\int_0^\infty\phi_r^L(h^2+h_r^2)\;rdr\right)^{1/2}\\
\le&~{}\frac CL\mathcal E_L.
\end{aligned}\]

For the second and fourth terms in \eqref{eq:local energy identity}, let $z$ be a radial scalar profile. By \eqref{eq:scalar radial Laplacian},
\[\begin{aligned}
\int_0^\infty\phi_r^L(L_sz)^2\;rdr
=&~{}\int_0^\infty\phi_r^L\left(z_{rr}^2+\frac{z_r^2}{r^2}\right)\;rdr+2\int_0^\infty\phi_r^Lz_{rr}z_r\;dr\\
=&~{}\int_0^\infty\phi_r^L\left(z_{rr}^2+\frac{z_r^2}{r^2}\right)\;rdr-\int_0^\infty\phi_{rr}^Lz_r^2\;dr.
\end{aligned}\]
Since $\phi_{rr}^L\le0$ by \eqref{eq:sign of phi_rr},
\begin{equation}\label{eq:z_rr bound}
\int_0^\infty\phi_r^L\left(z_{rr}^2+\frac{z_r^2}{r^2}\right)\;rdr\le\int_0^\infty\phi_r^L(L_sz)^2\;rdr.
\end{equation}
For the second term, applying \eqref{eq:z_rr bound} to $z=g$ and using \eqref{eq:norm upper}, we obtain
\[\int_0^\infty\phi_r^L\left(g_r^2+g_{rr}^2+\frac{g_r^2}{r^2}\right)\;rdr\le C\mathcal E_L.\]
Therefore, by Cauchy--Schwarz,
\[\begin{aligned}
&~{}\left|\int_0^\infty\phi_r^L\left(qg_r+q_rg_{rr}+\frac{qg_r}{r^2}\right)\;rdr\right|\\
&~{}\quad \le \left(\int_0^\infty\phi_r^L\left(q^2+q_r^2+\frac{q^2}{r^2}\right)\;rdr\right)^{1/2} \left(\int_0^\infty\phi_r^L\left(g_r^2+g_{rr}^2+\frac{g_r^2}{r^2}\right)\;rdr\right)^{1/2}\\
&~{}\quad \le C\mathcal E_L.
\end{aligned}
\]
For the third term, by Cauchy--Schwarz and \eqref{eq:norm upper},
\[\begin{aligned}
&~{}\left|\int_0^\infty\phi_r^L\left(hD_rf+h_rL_vf\right)\;rdr\right|\\
&~{}\quad \le \left(\int_0^\infty\phi_r^L\left(h^2+h_r^2\right)\;rdr\right)^{1/2}\left(\int_0^\infty\phi_r^L\left((D_rf)^2+(L_vf)^2\right)\;rdr\right)^{1/2}\\
&~{}\quad \le C\mathcal E_L.
\end{aligned}
\]
For the fourth term, since
\begin{equation}\label{eq:young ineq}
\left(K_s(q^2)-q^2\right)^2\le2\left(K_s(q^2)\right)^2+2q^4,
\end{equation}
applying \eqref{eq:z_rr bound} to $z=K_s(q^2)$ and using \eqref{eq:resolvent commutation}, \eqref{eq:nonlinear error_1}, and \eqref{eq:nonlinear_2}, we obtain
\[
\begin{aligned}
&~{}\int_0^\infty\phi_r^L\left(\left(\partial_{rr}K_s(q^2)\right)^2+\frac{\left(\partial_rK_s(q^2)\right)^2}{r^2}\right)\;rdr\\
&~{}\quad \le\int_0^\infty\phi_r^L\left(L_sK_s(q^2)\right)^2\;rdr\\
&~{}\quad =\int_0^\infty\phi_r^L\left(K_s(q^2)-q^2\right)^2\;rdr\\
&~{}\quad \le C\varepsilon^2\mathcal E_L.
\end{aligned}
\]
Together with \eqref{eq:nonlinear error_1}, this yields
\[\int_0^\infty\phi_r^L\left(\left(\partial_rK_s(q^2)\right)^2+\left(\partial_{rr}K_s(q^2)\right)^2+\frac{\left(\partial_rK_s(q^2)\right)^2}{r^2}\right)\;rdr\le C\varepsilon^2\mathcal E_L.\]
Hence, by Cauchy--Schwarz,
\[\begin{aligned}
&~{}\left|\int_0^\infty\phi_r^L\left(q\partial_rK_s(q^2)+q_r\partial_{rr}K_s(q^2)+\frac q{r^2}\partial_rK_s(q^2)\right)\;rdr\right|\\
\le&~{}\left(\int_0^\infty\phi_r^L\left(q^2+q_r^2+\frac{q^2}{r^2}\right)\;rdr\right)^{1/2}\\
&~{}\times\left(\int_0^\infty\phi_r^L\left(\left(\partial_rK_s(q^2)\right)^2+\left(\partial_{rr}K_s(q^2)\right)^2+\frac{\left(\partial_rK_s(q^2)\right)^2}{r^2}\right)\;rdr\right)^{1/2}\\
\le&~{}C\varepsilon\mathcal E_L.
\end{aligned}\]

Finally, for the fifth term, by \eqref{eq:resolvent commutation} and the same estimate as in \eqref{eq:young ineq},
\[\left(L_vK_v(qh)\right)^2=\left(K_v(qh)-qh\right)^2\le2\left(K_v(qh)\right)^2+2q^2h^2.\]
Therefore, using \eqref{eq:nonlinear error_1} and \eqref{eq:nonlinear_2}, we obtain
\[\int_0^\infty\phi_r^L\left(\left(D_rK_v(qh)\right)^2+\left(L_vK_v(qh)\right)^2\right)\;rdr\le C\varepsilon^2\mathcal E_L.\]
Then, by Cauchy--Schwarz,
\[\begin{aligned}
&~{}\left|\int_0^\infty\phi_r^L\left(hD_rK_v(qh)+h_rL_vK_v(qh)\right)\;rdr\right|\\
\le&~{}\left(\int_0^\infty\phi_r^L\left(h^2+h_r^2\right)\;rdr\right)^{1/2}
\left(\int_0^\infty\phi_r^L\left(\left(D_rK_v(qh)\right)^2+\left(L_vK_v(qh)\right)^2\right)\;rdr\right)^{1/2}\\
\le&~{}C\varepsilon\mathcal E_L.
\end{aligned}\]
Substituting all these estimates into \eqref{eq:local energy identity} and using $L\ge1$ and $\varepsilon\le\varepsilon_0 \le 1$, we obtain
\begin{equation}\label{eq:local energy bound}
\left|\frac{d}{dt}\mathscr E_L(t)\right|\le C_{a,c}\mathcal E_L(t).
\end{equation}

By \eqref{eq:E_L liminf}, there exist sequences $t_k^\pm\to\pm\infty$ such that
\[\mathcal E_L(t_k^\pm)\to0.\]
Hence, by \eqref{eq:local energy comparison},
\[\mathscr E_L(t_k^\pm)\to0.\]
We first consider the limit as $t\to\infty$. Fix $t\ge0$ and choose $k$ sufficiently large so that $t_k^+>t$. Integrating \eqref{eq:local energy bound} over $[t,t_k^+]$, we obtain
\[\begin{aligned}
0\le\mathscr E_L(t)
=&~{}\mathscr E_L(t_k^+)-\int_t^{t_k^+}\frac{d}{ds}\mathscr E_L(s)\;ds\\
\le&~{}\mathscr E_L(t_k^+)+C_{a,c}\int_t^{t_k^+}\mathcal E_L(s)\;ds.
\end{aligned}\]
By \eqref{eq:E_L integrable},
\[\int_t^\infty\mathcal E_L(s)\;ds<\infty.\]
Therefore, letting $k\to\infty$ and using $\mathscr E_L(t_k^+)\to0$, we obtain
\[0\le\mathscr E_L(t)\le C_{a,c}\int_t^\infty\mathcal E_L(s)\;ds.\]
Using the lower bound in \eqref{eq:local energy comparison}, we conclude that
\[0 \le \mathcal E_L(t)\le C_{a,c}\int_t^\infty\mathcal E_L(s)\;ds.\]
Since \eqref{eq:E_L integrable} implies
\[\int_t^\infty\mathcal E_L(s)\;ds\to 0 \qquad \text{as} \;\; t\to\infty,\]
we obtain
\[\mathcal E_L(t)\to 0 \qquad \text{as} \;\; t\to\infty.\]
The limit as $t\to-\infty$ follows analogously by using the sequence $t_k^-\to-\infty$.

Finally, fix $R>0$. For $|x|\le R$,
\[\left(1+\frac{|x|^2}{L^2}\right)^{3/4}\le\left(1+\frac{R^2}{L^2}\right)^{3/4}.\]
Hence, by \eqref{eq:local energy polar form},
\[\begin{aligned}
&~{}\|u(t)\|_{H^1(B_R)}^2+\|\eta(t)\|_{H^1(B_R)}^2\\
&~{}\quad \le \left(1+\frac{R^2}{L^2}\right)^{3/4}\int_{\R^2}\left(1+\frac{|x|^2}{L^2}\right)^{-3/4}\left(|u|^2+|Du|^2+\eta^2+|\nabla\eta|^2\right)(t,x)\;dx\\
&~{}\quad  = 2\pi L\left(1+\frac{R^2}{L^2}\right)^{3/4}\mathcal E_L(t).
\end{aligned}\]
Since $\mathcal E_L(t)\to0$ as $t\to\pm\infty$, this proves \eqref{eq:local H1 decay}.
\end{proof}

\subsection{Decay on expanding regions}\label{sec:expanding}

We now prove Corollary~\ref{cor:expanding}. By the scaling of Section~\ref{sec:4}, it suffices to consider the normalized system $b=d=1$. This scaling changes the ball $B_{L(t)}$ into $B_{\sqrt b\,L(t)}$. The $H^1$ norms are modified by constants depending on $b$ only, which does not affect the required statement. We take $\varepsilon_0$ and $L_0$ as in Proposition~\ref{prop:coercivity}. The new ingredient will be the following estimate on the dependence of $\mathcal H_{\phi^L}$ on the scale $L$.

\begin{lemma}\label{lem:L derivative}
Let $L\ge1$ be a parameter. Let $(u,\eta)(t)\in \mathbb H^1(\R^2)$ be radial profiles, under the representation by $q$ and $h$ given by $u=qe_r$, $\eta=h$. Then the mapping 
\[
L \mapsto\mathcal H_{\phi^L}
\]
is differentiable for $L>1$ and we have the bound
\begin{equation}\label{eq:L derivative bound}
\left|\partial_L\mathcal H_{\phi^L}\right|\le\frac CL\left(\|u\|_{H^1(\R^2)}^2+\|\eta\|_{H^1(\R^2)}^2\right),
\end{equation}
with $C>0$ an absolute constant.
\end{lemma}

\begin{proof}
Recall \eqref{eq:weight choice} and the scaled weight $\phi^L(r)=\phi(r/L)$. We have
\[\partial_L\phi^L(r)=-\frac r{L^2}\,\phi'\!\left(\frac rL\right)=-\frac rL\,\phi_r^L(r),\qquad \partial_L\left(\frac{\phi^L}r\right)=-\frac1L\,\phi_r^L(r).\]
Since $\phi'(s)=(1+s^2)^{-3/4}$, we have the quick estimate $\sup_{s>0}s\,\phi'(s)=2^{1/2}3^{-3/4}<1$. Additionally, $\phi'\le1$, so that
\begin{equation}\label{eq:L derivative weight}
\left|\partial_L\phi^L\right|\le\frac1L,\qquad \left|\partial_L\left(\frac{\phi^L}r\right)\right|\le\frac1{L^2}.
\end{equation}
Recall that $\mathcal H_\phi=I_\phi-\frac12\mathcal J_\phi-\frac12\mathcal K_\phi$, where $I_\phi$, $\mathcal J_\phi$ and $\mathcal K_\phi$ are given by \eqref{eq:principal virial functional}, \eqref{eq:tangential derivative correction} and \eqref{eq:inverse tangential correction}, and that each of them depends linearly on the weight. Differentiating in $L$ and using \eqref{eq:L derivative weight} together with Cauchy--Schwarz in $L^2(rdr)$ and the bound $(q_r+q/r)^2\le2(q_r^2+q^2/r^2)=2|Du|^2$, we get
\[
\left|\partial_LI_{\phi^L}\right|\le\frac{2\pi}L\int_0^\infty\left(|qh|+|q_rh_r|\right)rdr\le\frac1L\left(\|u\|_{H^1}^2+\|\eta\|_{H^1}^2\right),
\]
\[
\left|\partial_L\mathcal J_{\phi^L}\right|\le\frac{2\pi}{L^2}\int_0^\infty|h|\,|D_rq|\,rdr\le\frac{\sqrt2}{L^2}\,\|\eta\|_{L^2}\|Du\|_{L^2}.
\]
Now we deal with the term $\mathcal K_{\phi^L}$. Let $v=D_r^{-1}h$, then
\[\partial_L\mathcal K_{\phi^L}=2\pi\int_0^\infty\partial_L\left(\frac{\phi^L}{r}\right)q\,v\,rdr=2\pi\int_0^\infty\left(\partial_L\phi^L\right)q\,\frac vr\,rdr,\]
so that, using \eqref{eq:L derivative weight} and the Hardy inequality \eqref{eq:radial inverse Hardy},
\[\left|\partial_L\mathcal K_{\phi^L}\right|\le\frac{2\pi}L\|q\|_{L^2(rdr)}\left\|\frac vr\right\|_{L^2(rdr)}\le\frac1L\|u\|_{L^2}\|\eta\|_{L^2}.\]
Gathering the three previous bounds and using that $L\ge1$ lead to \eqref{eq:L derivative bound}.
\end{proof}

\begin{proof}[Proof of Corollary~\ref{cor:expanding}]
As in Remark~\ref{rem:behavior at origin} and the final remark below, it suffices to consider smooth solutions, the general case following by density. Let $L:[T_0,\infty)\to[L_0,\infty)$ be as in the statement and set
\[\mathscr H(t):=\mathcal H_{\phi^{L(t)}}(t),\qquad t\ge T_0.\]
Since $\mathcal H_{\phi^L}(t)$ is a $C^1$ function of $(L,t)$, the chain rule gives
\[\mathscr H'(t)=\left.\frac{d}{dt}\mathcal H_{\phi^{L}}(t)\right|_{L=L(t)}+L'(t)\left.\partial_L\mathcal H_{\phi^L}(t)\right|_{L=L(t)}.\]
The first term is estimated from below by Proposition~\ref{prop:coercivity}, which holds for every fixed $L\ge L_0$. The second term is bounded using Lemma~\ref{lem:L derivative}, the property $L'\ge0$ and the uniform in time bound $\sup_t(\|u(t)\|_{H^1}+\|\eta(t)\|_{H^1})\le\varepsilon$:
\[\mathscr H'(t)\ge2\pi\kappa\,\mathcal E_{L(t)}(t)-C\varepsilon^2\,\frac{L'(t)}{L(t)}.\]
Integrating the previous inequality on $[T_0,T]$ and using \eqref{eq:H_phi bound}, and the direct bounds $\|\phi^L\|_{L^\infty}=C_*$ and $\|\phi^L/r\|_{L^\infty}\le L^{-1}\le1$, we obtain
\begin{equation}\label{eq:expanding integral}
\begin{aligned}
2\pi\kappa\int_{T_0}^T\mathcal E_{L(t)}(t)\,dt\le &~{}\mathscr H(T)-\mathscr H(T_0)+C\varepsilon^2\log\frac{L(T)}{L(T_0)} \\
\le &~{} C\varepsilon^2\left(1+\log\frac{L(T)}{L(T_0)}\right).
\end{aligned}
\end{equation}
On the other hand, by \eqref{eq:local energy polar form} and the fact that $(1+|x|^2/L^2)^{-3/4}\ge2^{-3/4}$ on $B_L$,
\[\|u(t)\|_{H^1(B_{L(t)})}^2+\|\eta(t)\|_{H^1(B_{L(t)})}^2\le2^{3/4}\,2\pi L(t)\,\mathcal E_{L(t)}(t).\]
Suppose, by contradiction, that there exist $\delta>0$ and $T_1\ge T_0$ such that 
\[
\|u(t)\|_{H^1(B_{L(t)})}^2+\|\eta(t)\|_{H^1(B_{L(t)})}^2\ge\delta
\]
 for all $t\ge T_1$. Then, by \eqref{eq:expanding integral},
\[\frac{\delta}{2^{3/4}}\int_{T_1}^T\frac{dt}{L(t)}\le C\varepsilon^2\left(1+\log\frac{L(T)}{L(T_0)}\right)\qquad\text{for all }T\ge T_1.\]
Since $L(T)\to\infty$, condition \eqref{eq:expanding condition} forces $\int_{T_0}^T dt/L(t)\to\infty$ and then 
\[
\log L(T)=o\big(\int_{T_0}^Tdt/L(t)\big).
\]
 This fact contradicts the last inequality as $T\to\infty$. This proves the claim as $t\to+\infty$. The case $t\to-\infty$ follows by applying the previous argument to $(-u,\eta)(-t,x)$, which is again a radial solution of \eqref{eq:normalized abcd} satisfying the same bounds.

Finally, for $L(t)=t^\alpha$ with $0<\alpha<1$ one has $\int_{T_0}^Tdt/L(t)\sim T^{1-\alpha}/(1-\alpha)$ while $\log L(T)=\alpha\log T$; for $L(t)=t/(\log t)^\beta$, $\beta>0$, which is nondecreasing for $\log t>\beta$, one has $\int_{T_0}^Tdt/L(t)\sim(\log T)^{\beta+1}/(\beta+1)$ while $\log L(T)\le\log T$. In both cases \eqref{eq:expanding condition} holds.
\end{proof}

We finish this work with a final remark.

\begin{remark}
As the reader may have noticed, all decay computations performed in this paper and leading to Theorem~\ref{thm:MT} were obtained by assuming smooth radial solutions. By standard arguments, they also apply to the $\mathbb H^1$ solutions of Theorem~\ref{thm:GWP}. We mention the standard procedure to make this rigorous:

\medskip

\noindent
\emph{Radial symmetry.} If $(u_0,\eta_0)$ is the initial data and it is radial, then by uniqueness of the solution, and rotation invariance of the model \eqref{2D ABCD}, the solution remains radial for all time. Therefore, by Theorem~\ref{thm:GWP}, small radial initial data evolves into global radially symmetric solutions. Consequently, Theorem~\ref{thm:MT} applies to all small radial data.

\medskip

\noindent
\emph{Persistence of regularity.} 
By Theorem \ref{thm:GWP} (see also Remark \ref{rem:persistence regularity}), the higher Sobolev regularity of the solution persists as long as the $\mathbb H^1$ norm remains bounded in time. Therefore, smooth $\mathbb H^3$ small radial data produce globally defined solutions in the class $C(\R;\mathbb H^3)$, enough for all the virial computations, which require at most two derivatives of $q,h$ (i.e.\ the radial representatives of the solution), after the resolvents are taken into account, and the boundary terms at $r=0$ vanish (Remark~\ref{rem:behavior at origin}). Finally, decay at $r \to \infty$ can be handled with Schwartz data.

\medskip

\noindent
\emph{Density.} In order to obtain a decay conclusion for radial data only in $\mathbb H^1$, we approximate $(u_0,\eta_0)\in\mathbb H^1$ radial by smooth radial data $(u_0^n,\eta_0^n)$ in $\mathbb H^1$. The integrated inequalities \eqref{eq:local energy integrability} and the Gronwall form $\mathscr E_L(t)\le\mathscr E_L(t')+C\int_t^{t'}\mathcal E_L$, obtained by integrating \eqref{eq:local energy bound}, involve only $\mathbb H^1$ quantities on compact time intervals, so they pass to the limit by the Lipschitz continuity of the flow in $\mathbb H^1$ on $[-T,T]$, obtained from Theorem \ref{thm:LWP} and the uniform bound of Theorem \ref{thm:GWP}. The limiting argument of Theorem \ref{thm:local energy decay} then applies verbatim.

A similar remark applies to Theorem \ref{thm:GWP}, which uses the conservation of the Hamiltonian/energy $E(u,\eta)(t)=E(u_0,\eta_0)$ for $\mathbb H^1$ solutions. Notice that $E$ is continuous on $\mathbb H^1$ by \eqref{Sob}. 
\end{remark}

\subsection*{AI Declaration} For this paper, no AI tools of any kind have been used except for improving the writing. The authors take full responsibility for the mathematical arguments of the final manuscript.


\begin{thebibliography}{99}
\bibitem{ASL} B. Alvarez-Samaniego and D. Lannes, \emph{Large time existence for 3D water-waves and asymptotics}, Invent. Math. 171 (2008), no. 3, 485--541.

\bibitem{BCSM2026} J. Bae, S. Cuccagna, D. Del Santo, and M. Maeda, \emph{The asymptotic stability of the stable solitons of the generalized good Boussinesq equation}, preprint arXiv:2610.25939, 2026.

\bibitem{BCM2026} J. Bae, S. Cuccagna, and M. Maeda, \emph{Conditional asymptotic stability of solitary waves of the Euler--Poisson system on the line}, preprint arXiv:2603.05334, 2026.

\bibitem{BonaChen1998} J. L. Bona and M. Chen, \emph{A Boussinesq system for two-way propagation of nonlinear dispersive waves}, Phys. D 116 (1998), no. 1-2, 191--224.

\bibitem{BCS1} J. L. Bona, M. Chen, and J.-C. Saut, \emph{Boussinesq equations and other systems for small-amplitude long waves in nonlinear dispersive media. I: Derivation and linear theory}, J. Nonlinear. Sci. Vol. 12: pp. 283--318 (2002).

\bibitem{BCS2} J. L. Bona, M. Chen, and J.-C. Saut, \emph{Boussinesq equations and other systems for small-amplitude long waves in nonlinear dispersive media. II: The nonlinear theory}, Nonlinearity 17 (2004) 925--952.

\bibitem{BCL} J. L. Bona, T. Colin and D. Lannes, \emph{Long wave approximations for water waves}, Arch. Ration. Mech. Anal. 178 (2005), no. 3, 373--410.

\bibitem{BLS} J. L. Bona, D. Lannes, and J.-C. Saut, \emph{Asymptotic models for internal waves}, J. Math. Pures Appl. (9) 89 (2008), no. 6, 538--566.

\bibitem{Bous} J. Boussinesq, \emph{Th\'eorie des ondes et des remous qui se propagent le long d'{}un canal rectangulaire horizontal, en communiquant au liquide contenu dans ce canal des vitesses sensiblement pareilles de la surface au fond}, J. Math. Pure Appl. (2) 17 (1872), 55--108.

\bibitem{Burtea2016} C. Burtea, \emph{New long time existence results for a class of Boussinesq-type systems}, J. Math. Pures Appl. (9) 106 (2016), no. 2, 203--236.

\bibitem{CLY2026} G. Chen, Y. Li, and Z. Yang, \emph{Asymptotic stability of moving kinks with an internal mode}, preprint arXiv:2611.13279, 2026.

\bibitem{CSM2026} S. Cuccagna, D. Del Santo, and M. Maeda, \emph{The virial--Kato method and the asymptotic stability of solitons of the generalized Benjamin, Bona and Mahony equation}, preprint arXiv:2610.41549, 2026.

\bibitem{CM2023} S. Cuccagna and M. Maeda, \emph{Asymptotic stability of kink with internal modes under odd perturbation}, Nonlinear Differential Equations Appl. NoDEA \textbf{30} (2023), Paper No. 1.

\bibitem{CM2027} S. Cuccagna and M. Maeda, \emph{On stabilization at a soliton for generalized Korteweg--de Vries pure power equation for any power $p\in(1,5)$}, Nonlinear Anal. \textbf{274} (2027), 114253.

\bibitem{Chen1998} M. Chen, \emph{Exact solutions of various Boussinesq systems}, Appl. Math. Lett. 11 (1998), no. 5, 45--49.

\bibitem{CL2013} R. M. Chen and Y. Liu, \emph{On the ill-posedness of a weakly dispersive one-dimensional Boussinesq system}, J. Anal. Math. 121 (2013), 299--316.

\bibitem{ChenNguyenSun2010} M. Chen, N. V. Nguyen and S.-M. Sun, \emph{Solitary-wave solutions to Boussinesq systems with large surface tension}, Discrete Contin. Dyn.
  Syst. 26 (2010), no. 4, 1153--1184.

\bibitem{ChenNguyenSun2011} M. Chen, N. V. Nguyen and S.-M. Sun, \emph{Existence of traveling-wave solutions to Boussinesq systems}, Differential Integral Equations
  24 (2011), no. 9-10, 895--908.

\bibitem{CMPS2016} R. C\^ote, C. Mu\~noz, D. Pilod and G. Simpson, \emph{Asymptotic stability of high-dimensional Zakharov--Kuznetsov solitons}, Arch. Ration. Mech. Anal. 220 (2016), no. 2, 639--710.

\bibitem{DD} P. Daripa, and R. K. Dash, \emph{On a class of Boussinesq equations for shallow water waves}, Computational science and its applications--ICCSA 2003. Part II, 764--773, Lecture Notes in Comput. Sci., 2668, Springer, Berlin, 2003.

\bibitem{DMS2007} V.~A. Dougalis, D.~E. Mitsotakis and J.-C. Saut, \emph{On some Boussinesq systems in two space dimensions: theory and numerical analysis} M2AN Math. Model. Numer. Anal. 41 (2007), no. 5, 825--854.

\bibitem{DMS2009} V.~A. Dougalis, D.~E. Mitsotakis and J.-C. Saut, \emph{On initial-boundary value problems for a Boussinesq system of BBM-BBM type in a plane domain}, Discrete Contin. Dyn. Syst. 23 (2009), no. 4, 1191--1204.   

\bibitem{Grisvard} P. Grisvard, \emph{Quelques propri\'et\'es des espaces de Sobolev, utiles dans l'\'etude des \'equations de Navier-Stokes}, Probl\`emes d'{}\'evolution, non lin\'eaires, S\'eminaire de Nice (1974-1976).

\bibitem{Hadamard} J. Hadamard, \emph{Sur les probl\`emes aux d\'eriv\'ees partielles et leur signification physique}, Princeton University Bulletin (1902) Vol. 13  49--52.

\bibitem{HMM1} S.~Herr, C.~Maul\'en, and C.~Mu\~noz, \emph{Decay of solutions of nonlinear Dirac equations}. Comm. Math. Phys. 407 (2026), no. 5, Paper No. 94, 52 pp.

\bibitem{HMM2} S.~Herr, C.~Maul\'en, and C.~Mu\~noz, \emph{Decay of solutions of nonlinear Dirac equations: the 2D case}. preprint 2026 to appear in SIAM J. Math. Anal.

\bibitem{HTT} S. Herr, D. Tataru, and N. Tzvetkov, \emph{Global well-posedness of the energy-critical nonlinear Schr\"odinger equation with small initial data in $H^1(\mathbb T^3)$}, Duke Math. J. 159 (2011), no. 2, 329--349.

\bibitem{KimKwak2026} S. Kim and C. Kwak, \emph{Well-posedness issues for the generalized Benjamin--Bona--Mahony equation}, preprint arXiv:2603.21060 [math.AP] (2026).

\bibitem{KMM2017} M. Kowalczyk, Y. Martel and C. Mu\~noz, \emph{Kink dynamics in the $\phi^4$ model: asymptotic stability for odd perturbations in the energy space},
  J. Amer. Math. Soc. 30 (2017), no. 3, 769--798.

\bibitem{KMM2017b} M. Kowalczyk, Y. Martel and C. Mu\~noz, \emph{Nonexistence of small, odd breathers for a class of nonlinear wave equations}, Lett. Math. Phys.
  107 (2017), no. 5, 921--931.

\bibitem{KMMV2021} M. Kowalczyk, Y. Martel, C. Mu\~noz and H. Van Den Bosch, \emph{A sufficient condition for asymptotic stability of kinks in general $(1+1)$-scalar field models}, Ann. PDE 7 (2021), no. 1, Paper No. 10.

\bibitem{KwakMaulen} C. Kwak, and C. Maul\'en, \emph{Ill-Posedness Issues on (abcd)-Boussinesq System}, J. Dyn. Diff. Eqns. 36 (2), 1123--1152 (2024).

\bibitem{KwakMunoz2019} C. Kwak and C. Mu\~noz, \emph{Extended decay properties for generalized BBM equation}, in: Nonlinear dispersive partial differential equations
  and inverse scattering, Fields Inst. Commun. 83, Springer, New York, 2019, 397--411.

\bibitem{KM2020} C. Kwak and C. Mu\~noz, \emph{Asymptotic dynamics for the small data weakly dispersive one-dimensional Hamiltonian ABCD system}, Trans. Amer. Math. Soc. 373 (2020), no. 2, 1043--1107.

\bibitem{KMPP2019} C. Kwak, C. Mu\~noz, F. Poblete and J. C. Pozo, \emph{The scattering problem for Hamiltonian ABCD Boussinesq systems in the energy space},
  J. Math. Pures Appl. (9) 127 (2019), 121--159.     

\bibitem{Lannes2013} D. Lannes, \emph{The water waves problem. Mathematical analysis and asymptotics}, Mathematical Surveys and Monographs 188, American Mathematical
  Society, Providence, RI, 2013.

\bibitem{LPS} F. Linares, D. Pilod, and J.~C. Saut, \emph{Well-posedness of strongly dispersive two-dimensional surface wave Boussinesq systems},  SIAM J. Math. Anal. 44 (2012), no. 6, 4195--4221.

\bibitem{MartelMerle2000} Y. Martel and F. Merle, \emph{A Liouville theorem for the critical generalized Korteweg-de Vries equation}, J. Math. Pures Appl. (9) 79 (2000), no. 4, 339--425.
 
 \bibitem{MM0} Y. Martel, and F. Merle, \emph{Instability of solitons for the critical generalized Korteweg de Vries equation}, GAFA 11 (2001), 74-123.

\bibitem{MM1} Y. Martel and F. Merle, \emph{Asymptotic stability of solitons for subcritical generalized {K}d{V} equations}, Arch. Ration. Mech. Anal. {157} (2001), no.~3, 219--254.

\bibitem{MMnon} Y. Martel and F. Merle, \emph{Asymptotic stability of solitons of the subcritical gKdV equations revisited}, Nonlinearity {18} (2005) 55--80.

\bibitem{MMas1} Y. Martel and F. Merle, \emph{Asymptotic stability of solitons of the gKdV equations with general nonlinearity}, Math. Ann. {34} (2008), no. 2, 391--427.

\bibitem{MMas2} Y. Martel and F. Merle, \emph{Refined asymptotics around solitons for gKdV equations}, Disc. Cont. Dyn. Syst. {20} (2008), no. 2, 177--218.

\bibitem{Mau} C. Maul\'en, \emph{Asymptotic stability manifolds for solitons in the generalized good Boussinesq equation}, J. Math. Pures Appl. (9) \textbf{177} (2023), 260--328.

\bibitem{MaMu0} C. Maul\'en, and C. Mu\~noz, \emph{Decay in the one dimensional generalized Improved Boussinesq equation}, SN Partial Differ. Equ. Appl. \textbf{1} (2020), Paper No. 1.

\bibitem{MaMu} C. Maul\'en, and C. Mu\~noz, \emph{Asymptotic stability of the fourth order $\phi^4$ kink for general perturbations in the energy space}, Ann. Inst. H. Poincar\'e C Anal. Non Lin\'eaire 42 (2025), no. 3, 647-714.

\bibitem{MaMu2} C. Maul\'en, and C. Mu\~noz, \emph{On asymptotic stability of stable Good Boussinesq solitary waves}, preprint, arXiv:2501.13649 (2025).

\bibitem{MMP} C. Maul\'en, C. Mu\~noz, and F. Poblete, \emph{Decay of small energy ABCD Boussinesq waves under the influence of an uneven bottom}, preprint, 2025.

\bibitem{MS} V. Maz'ya and T. Shaposhnikova, \emph{Jacques Hadamard, A Universal Mathematician}.  Amer. Math. Soc., Providence, RI, 1998.

\bibitem{MMPP2021} A. J. Mendez, C. Mu\~noz, F. Poblete, and J. C. Pozo, \emph{On local energy decay for large solutions of the Zakharov-Kuznetsov equation}, Comm. Partial Differential Equations \textbf{46} (2021), no.~8, 1440--1487.

\bibitem{MMPP2} A. J. Mendez, C. Mu\~noz, F. Poblete, and J. C. Pozo, \emph{Long time asymptotics of large data in the Kadomtsev-Petviashvili models}, Nonlinearity \textbf{37} (2024), no.~5, Paper No. 055017, 28 pp.

\bibitem{MSZ2012} M. Ming, J.-C. Saut and P. Zhang, \emph{Long-time existence of solutions to Boussinesq systems}, SIAM. J. Math. Anal. 44 (2012) 4078--4100.

\bibitem{MZZ1} M. Ming, P. Zhang and Z. Zhang, \emph{Large time well-posedness of the three-dimensional capillary-gravity waves in the long wave regime}, Arch. Ration. Mech. Anal. 204 (2012) no. 2 387--444.

\bibitem{MZZ2} M. Ming, P. Zhang and Z. Zhang, \emph{Long-wave approximation to the 3-D capillary-gravity waves}, SIAM J. Math. Anal. 44 (2012) no. 4 2920--2948.

\bibitem{MunozPonce2019} C. Mu\~noz and G. Ponce, \emph{Breathers and the dynamics of solutions in KdV type equations}, Comm. Math. Phys. 367 (2019), no. 2, 581--598.

\bibitem{MPS2021} C. Mu\~noz, G. Ponce and J.-C. Saut, \emph{On the long time behavior of solutions to the intermediate long wave equation}, SIAM J. Math.  Anal. 53 (2021), no. 1, 1029--1048.

\bibitem{MPP2018} C. Mu\~noz, F. Poblete and J. C. Pozo, \emph{Scattering in the energy space for Boussinesq equations}, Comm. Math. Phys. 361 (2018), no. 1,  127--141.     

\bibitem{SX} J.-C. Saut and L. Xu, \emph{The Cauchy problem on large time for surface waves Boussinesq systems}, J. Math. Pures Appl. (9) 97 (2012) 635--662.

\bibitem{SautXu} J.-C. Saut and L. Xu, \emph{Long time existence for the Boussinesq--Full Dispersion systems}, J. Differential Equations 269 (2020), 2627--2663.

\bibitem{SWX} J.-C. Saut, C. Wang and L. Xu, \emph{The Cauchy problem on large time for surface-waves-type Boussinesq Systems II}, SIAM. J. Math. Anal. 49 (2017) 2321--2386.

\bibitem{Strauss1977} W. A. Strauss, \emph{Existence of solitary waves in higher dimensions}, Comm. Math. Phys. 55 (1977), no. 2, 149--162.

\bibitem{Tao2001} T. Tao, \emph{Multilinear weighted convolution of $L^2$ functions and applications to nonlinear dispersive equations}, Amer. J. Math. 123(5) (2001) 839--908.

\end{thebibliography}
\end{document}